%% file: charac.tex
\documentclass[12pt]{amsart}

\usepackage[T1]{fontenc}
\usepackage{times,mathptm}
\usepackage{amssymb,epsfig,verbatim,xypic}
\theoremstyle{plain}

\newtheorem{thm}{Theorem}[section]

\newtheorem{cor}[thm]{Corollary}
\newtheorem{pro}[thm]{Proposition}
\newtheorem{lem}[thm]{Lemma}
\newtheorem{proposition-principale}[thm]{Proposition principale}
\newtheorem{thm-principal}{Th\'eor\`eme principal}[section]

\theoremstyle{definition}

\newtheorem{eg}[thm]{Example}
\newtheorem{rem}[thm]{Remark}

\newenvironment{thm-A}
{{\vv \noindent \bf Theorem A.$\,$}\it}{\vv}

\newenvironment{thm-B}
{{\vv \noindent \bf Theorem B.$\,$}\it}{\vv}

\newenvironment{thm-C}
{{\vv \noindent \bf Theorem C.$\,$}\it}{\vv}

\newenvironment{thm-D}
{{\vv \noindent \bf Theorem D.$\,$}\it}{\vv}

\newenvironment{thm-E}
{{\vv \noindent \bf Theorem E.$\,$}\it}{\vv}

\def\vv{\vspace{0.2cm}}

\def\C{\mathbf{C}}
\def\R{\mathbf{R}}
\def\Q{\mathbf{Q}}
\def\Z{\mathbf{Z}}

\def\P{\mathbb{P}}
\def\Sphere{\mathbb{S}^2}

\def\T{\mathbb{T}^2}
\def\M{{S_M}}

\def\SS{{\sf{Fam}}}
\def\Sing{{\sf{Sing}}}

\def\Rep{{\sf{Rep}}}

\def\H{\mathbb{H}}
\def\B{\mathbb{B}}

\def\hol{{\mathbf{Mon}}}
\def\dev{{\mathbf{dev}}}
\def\Lin{{\mathbf{Lin}}}
\def\Af{{\mathbf{Aff}}\,}

\def\A{\mathcal{A}}
\def\Out{{\sf{Out}}}
\def\Aut{{\sf{Aut}}}
\def\Bir{{\sf{Bir}}}
\def\Sym{{\sf{Sym}}}

\def\Orb{{\text{Orb}}}

\def\Ok{{\text{Ok}}}
\def\Quad{{\text{Quad}}}

\def\MCG{{\sf{MCG}}}
\def\PGL{{\sf{PGL}}\,}
\def\PSL{{\sf{PSL}}\,}
\def\GL{{\sf{GL}}\,}
\def\SL{{\sf{SL}}\,}
\def\Aff{{\sf{Aff}}\,}
\def\SU{{\sf{SU}}\,}
\def\Diff{{\sf{Diff}}\,}
\def\SO{{\sf{SO}}\,}

\def\tr{{\sf{tr}}}
\def\Tr{{\sf{Tr}}}

\def\Ind{{\text{Ind}}}

\addtocounter{section}{0}             % Start with section 1
\numberwithin{equation}{section}       % Number formulas within sections
\begin{document}
\setlength{\baselineskip}{0.47cm}        % Previous 0.57
%
%%%%%%%%%%%%%%%%%%%%%%%%%%%%%%%%%%%%%%%%%%%%%%%%%%%%%%%%%%%%%%%%%%
%
\title[Dynamics, Painlev\'e VI and Character Varieties.]
{Holomorphic dynamics, Painlev\'e VI equation and Character Varieties.}
\date{2007}
\author{Serge Cantat, Frank Loray}
\address{D\'epartement de math\'ematiques\\
         Universit\'e de Rennes\\
         Rennes\\
         France}
\email{serge.cantat@univ-rennes1.fr} 
\email{frank.loray@univ-rennes1.fr}
%
%%%%%%%%%%%%%%%%%%%%%%%%%%%%%%%%%%%%%%%%%%%%%%%%%%%%%%%%%%%%%%%%%%
%
%\begin{abstract} . \end{abstract}
%
%%%%%%%%%%%%%%%%%%%%%%%%%%%%%%%%%%%%%%%%%%%%%%%%%%%%%%%%%%%%%%%%%%
%

\maketitle
\tableofcontents

%\newpage

%
%%%%%%%%%%%%%%%%%%%%%%%%%%%%%%%%%%%%%%%%%%%%%%%%%%%%%%%%%%%%%%%%%%
%

\section{Introduction}

%%%%%%%%%%%%%%%%%%%%%%%%%%%%%%%%%%%%%%%%

This is the first part of a series of two papers (see \cite{Cantat:BHPS}), the aim of which is to describe the dynamics of a polynomial  action of the group 
\begin{equation}
\Gamma_2^*=\{ M\in \PGL(2,\Z)\, \vert \quad M={\text{Id}} \, \,{\text{mod}}(2)\}
\end{equation}
on the family of affine cubic surfaces 
\begin{equation}
x^2+y^2+z^2+xyz=Ax+By+Cz+D,
\end{equation}
where $A$, $B$, $C$, and $D$ are complex parameters.
This dynamical system appears in several
different mathematical areas,  like the monodromy of the sixth Painlev\'e differential equation,
the geometry of hyperbolic threefolds, and  the spectral properties of certain discrete Schr\"odinger operators. One of our main goals here is to classify parameters $(A,B,C,D)$ for which $\Gamma_2^*$ preserves
a holomorphic geometric structure, and to apply this classification  to
provide a galoisian proof of the irreducibility of the sixth Painlev\'e equation. 

%%%%%%%%%%%%%%%%%%%%%%%

\subsection{Character variety}\label{par:char_var_nota}

%%%%%%%%%%%%%%%%%%%%%%%%%%%%%%%%%%%%%%%%

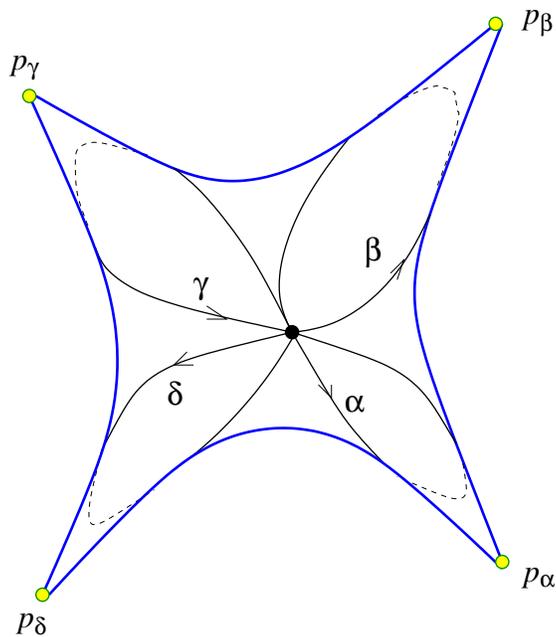
\begin{figure}[h]\label{fig:four_punctured_sphere}
\input{sphere-ter.pstex_t}
\caption{ {\sf{The four punctured sphere. }}}
\end{figure}

 Let $\Sphere_4$ be the four punctured sphere. Its fundamental group
is isomorphic to a free group of rank $3$; if $\alpha$, $\beta$, $\gamma$
and $\delta$ are the four loops which are depicted on figure 
\ref{fig:four_punctured_sphere},
then 
$$
\pi_1(\Sphere_4)= \langle \alpha,\beta,\gamma,\delta \, \vert \, \alpha\beta\gamma\delta = 1 \rangle.
$$

Let $\Rep(\Sphere_4)$ be the set of representations of $\pi_1(\Sphere_4)$ into $\SL(2,\C)$. Such 
a representation $\rho$ is uniquely determined by the $3$ matrices $\rho(\alpha)$, $\rho(\beta)$, and 
$\rho(\gamma)$, so that $\Rep(\Sphere_4)$ can be identified with the affine algebraic variety $(\SL(2,\C))^3$.
Let us associate the $7$ following traces to any element $\rho$ of $\Rep(\Sphere_4)$:
\begin{equation*}
\begin{array}{clclclc}
 a = \tr(\rho(\alpha)) & ; & b  =  \tr(\rho(\beta)) & ; & c = \tr(\rho(\gamma)) & ; &  d  =  \tr(\rho(\delta)) \\
 x  = \tr(\rho(\alpha \beta)) & ; &  y  =  \tr(\rho(\beta \gamma))   & ; & z  =  \tr( \rho(\gamma \alpha)) . & \, & \,
 \end{array}
\end{equation*}
The polynomial map $\chi:\Rep(\Sphere_4)\to \C^7$ defined by 
\begin{equation}
\chi(\rho)=(a,b,c,d,x,y,z)
\end{equation}
is invariant under conjugation, by which we mean that $\chi(\rho')=\chi(\rho)$ 
if $\rho'$ is conjugate to $\rho$ by an element of $\SL(2,\C).$ Moreover,
\begin{enumerate}
\item the algebra of polynomial functions on $\Rep(\Sphere_4)$ which are invariant
under conjugation is generated by the components of $\chi$;
\item the components of $\chi$ satisfy the quartic equation 
\begin{equation}\label{eq:surface}
x^2+y^2+z^2+x y z = A x + By +C z + D,
\end{equation}
in which the variables $A$, $B$, $C$, and $D$ are given by
\begin{equation}\label{eq:parameters}
\begin{array}{c} A  =  ab + cd, \quad B  = bc + ad,\quad  C  =  ac + bd,  \\
 {\text{and}} \quad  D  =  4 - a^2 - b^2- c^2 - d^2 - abcd .  \end{array} 
\end{equation}
\item the algebraic quotient $\Rep(\Sphere_4)/\! /\SL(2,\C)$ of  $\Rep(\Sphere_4)$
by the action of $\SL(2,\C)$ by conjugation is isomorphic to the six-dimensional
quartic hypersurface of $\C^7$ defined by equation (\ref{eq:surface}). 
\end{enumerate}
The affine algebraic variety $\Rep(\Sphere_4)/\! /\SL(2,\C)$ will be denoted  $\chi(\Sphere_4)$ and called the {\sl{character
variety of $\Sphere_4$}}. For each choice of four complex parameters $A$, $B$, $C$, and $D$, $S_{(A,B,C,D)}$
(or $S$ if there is no obvious possible confusion) will denote the cubic surface of $\C^3$ defined by the equation (\ref{eq:surface}).
The family of these surfaces $S_{(A,B,C,D)}$ will be denoted $\SS.$

\begin{rem}
The map $\C^4\to\C^4;(a,b,c,d)\mapsto(A,B,C,D)$ defined by (\ref{eq:parameters})
is a non Galois ramified cover of degree $24$. Fibers are studied
in Appendix B. It is important to notice that a point 
$m\in S_{(A,B,C,D)}$ will give
rise to representations of very different nature depending on 
the choice of $(a,b,c,d)$ in the fiber, e.g. reducible or irreducible,
finite or infinite image.
\end{rem}

\begin{rem} As we shall see in section \ref{par:OncePuncturedTorus}, if we replace the four puntured sphere
by the once puntured torus, the character variety is naturally fibered by the family of cubic surfaces $S_{(0,0,0,D}.$
\end{rem}

%%%%%%%%%%%%%%%%%%%%%%%%%%%%%%%%%%%%%%%%%%%%%%%%%%

\subsection{Automorphisms and modular groups}\label{par:automorphisms_modular} 

%%%%%%%%%%%%%%%%%%%%%%%%%%%%%%%%%%%%%%%%

The group of automorphisms $\Aut(\pi_1(\Sphere_4))$ acts on $\Rep(\Sphere_4)$
by composition: $(\Phi,\rho)\mapsto \rho \circ \Phi^{-1}$.
Since inner automorphisms act trivially on $\chi(\Sphere_4)$, we get a morphism from the
group of outer automorphisms $\Out(\pi_1(\Sphere_4))$ into the group of polynomial diffeomorphisms
of $\chi(\Sphere_4)$:
\begin{equation}
\left\{ \begin{array}{ccc} \Out(\pi_1(\Sphere_4)) & \to & \Aut[\chi(\Sphere_4)] \\
 \Phi & \mapsto & f_\Phi
\end{array} \right.
\end{equation}
such that $f_\Phi ( \chi (\rho) ) = \chi (\rho\circ \Phi^{-1})$ for any representation $\rho$.

\vv

The group $\Out(\pi_1(\Sphere_4))$ 
is isomorphic to the extended mapping class group
$\MCG^*(\Sphere_4)$, {\sl{i.e.}} to the group of isotopy classes of 
homeomorphisms of $\Sphere_4,$ that preserve or reverse the orientation.
It contains a copy of $\PGL(2,\Z)$ which is obtained
as follows.
Let $\T=\R^2/\Z^2$ be a torus of dimension $2$ and $\sigma$ be the
involution of $\T$ defined by $\sigma(x,y)=(-x,-y)$. The fixed point set of $\sigma$ is 
the $2$-torsion subgroup $H\subset\T,$ isomorphic to $\Z/2\Z\times\Z/2\Z$:
\begin{equation}
H=\{(0,0);\, (0,1/2);\, (1/2,0);\, (1/2,1/2)\}.
\end{equation}
The quotient $\T/\sigma$ is homeomorphic to the 
sphere, $\Sphere$, and the quotient map $\pi:\T\to \T/\sigma=\Sphere$  has four ramification
points, corresponding to the four fixed points of $\sigma$. 
The affine group $\GL(2,\Z)\ltimes H$ acts linearly on $\T$ and commutes with $\sigma$. This yields an action of
$\PGL(2,\Z)\ltimes H$  on the sphere $\Sphere,$ that permutes  ramification points
of $\pi$. Taking these four ramification points as the punctures of $\Sphere_4$, we get a morphism 
\begin{equation}
\PGL(2,\Z)\ltimes H\to \MCG^*(\Sphere_4),
\end{equation}
which, in fact, is an isomorphism (see \cite{Birman:book}, section 4.4).
The image of $\PGL(2,\Z)$ is the stabilizer of $\pi(0,0)$, freely permuting the three other points.
As a consequence, $\PGL(2,\Z)$ acts by polynomial 
transformations on $\chi(\Sphere_4)$. 
The image of $H$ permutes the $4$ punctures
by products of disjoint transpositions and acts trivially on $\chi(\Sphere_4),$ so that  the action of the whole mapping class group $\MCG^*(\Sphere_4)$ on $\chi(\Sphere_4)$ actually reduces to that 
of $\PGL(2,\Z)$ (see section \ref{par:MappingClassGroup}).

\vv

Let $\Gamma_2^*$ be the subgroup of $\PGL(2,\Z)$ whose elements coincide with 
the identity  modulo $2$. This group coincides with the (image in $\PGL(2,\Z)$ of the)
stabilizer of the fixed points of $\sigma$, so that $\Gamma_2^*$ acts on $\Sphere_4$ and fixes its four punctures.
Consequently, $\Gamma_2^*$ acts polynomially on $\chi(\Sphere_4)$ and preserves the fibers of the projection
$$
(a,b,c,d,x,y,z)\mapsto (a,b,c,d).
$$
From this we obtain, for any choice of four complex parameters $(A,B,C,D)$, a morphism from 
$\Gamma_2^*$ to the group of polynomial diffeomorphisms of the surface $S_{(A,B,C,D)}$.
The following result is essentially due to \`El'-Huti (see \cite{El-Huti:1974}, and \S \ref{par:triangle_automorphisms_area}). 

\begin{thm-A}
For any choice of the parameters $A$, $B$, $C$, and $D$, the morphism 
$$
\Gamma_2^*\to \Aut[S_{(A,B,C,D)}]  
$$
is injective and the index of its image is bounded by $24$. For a generic 
choice of the parameters, this morphism is an isomorphism. 
\end{thm-A}

As a consequence of this result, it suffices to understand the action of  $\Gamma^*_2$ 
on the surfaces $S_{(A,B,C,D)}$ in order to get a full understanding of the action of 
$\MCG^*(\Sphere_4)$ on $\chi(\Sphere_4)$. (see also remark \ref{rem:artin} 
for the case of orientation preserving transformations and an action of the pure 
braid group on three strings).

\begin{rem}
If the parameters $A$, $B$, $C$, and $D$ belong to a ring ${\mathbb{K}}$, the group $\Gamma^*_2$ 
acts on  $S_{(A,B,C,D)}({\mathbb{K}})$, i.e. on the set of points of the surface with coordinates 
in ${\mathbb{K}}$. In particular, when the parameters are real numbers, $\Gamma_2^*$ acts
on the real surface   $S_{(A,B,C,D)}(\R)$. 
\end{rem}

\vv

There are useful symmetries of the parameter space, as well as covering between 
distinct surfaces $S_{(A,B,C,D)}$ and $S_{(A',B',C',D')},$  that can be used to relate
dynamical properties of $\Gamma_2^*$ on different surfaces of our family.
These symmetries and covering will be described in section \ref{par:geometry-surfaces} and appendix B.

%%%%%%%%%%%%%%%%%%%%%%%%%%%%%%%%%%%%%%%%

\subsection{Projective structures}

%%%%%%%%%%%%%%%%%%%%%%%%%%%%%%%%%%%%%%%%

Once $\Sphere_4$ is endowed with a complex projective structure, which means
that we have an atlas on $\Sphere_4$ made of charts into $\P^1(\C)$ with transition
functions in the group of homographic transformations of $\P^1(\C)$, the holonomy defines
a morphism from $\pi_1(\Sphere_4)$ to $\PSL(2,\C)$. Since $\pi_1(\Sphere_4)$ is
a free group, the holonomy can be lifted to a morphism
$$
\rho : \pi_1(\Sphere_4)\to \SL(2,\C).
$$
Properties of the holonomy such as discreteness, finiteness, or the presence of parabolic elements
in $\rho(\pi_1(\Sphere_4))$, are invariant by conjugation and by the action of the mapping
class group $\MCG^*(\Sphere_4)$. This kind of properties may be used to construct invariant
subsets of $S_{(A,B,C,D)}$ for the action of $\Gamma_2^*$, and the dynamics of this action
may be used to understand those invariant sets and the space of projective structures. 
This approach has been popularized by Goldman
(see \cite{Goldman:1997}, \cite{Goldman:survey} for example).

%%%%%%%%%%%%%%%%%%%%%%%%%%%%%%%%%%%%%%%%

\subsection{ Painlev\'e VI equation }

%%%%%%%%%%%%%%%%%%%%%%%%%%%%%%%%%%%%%%%%

The dynamics of $\Gamma_2^*$ on the varieties $S_{(A,B,C,D)}$ is also related to the
monodromy of a famous ordinary differential equation. The sixth Painlev\'e
 equation $P_{VI}=P_{VI}(\theta_\alpha,\theta_\beta,\theta_\gamma,\theta_\delta)$
is the second order non linear O.D.E.
\begin{equation*}
P_{VI}\ \ \ \left\{
\begin{matrix}{ 
{d^2 q\over dt^2}}&=&
{{1\over 2}\left({1\over q}+{1\over q-1}+{1\over q-t}\right)
\left({dq\over dt}\right)^2
-\left({1\over t}+{1\over t-1}+{1\over q-t}\right)
\left({dq\over dt}\right)}\\
&&
{+{q(q-1)(q-t)\over t^2(t-1)^2}
\left(\frac{(\theta_\delta-1)^2}{2}-\frac{\theta_\alpha^2}{2}{t\over q^2}
+\frac{\theta_\beta^2}{2}{t-1\over (q-1)^2}+\frac{1-\theta_\gamma^2}{2}{t(t-1)\over
(q-t)^2}\right)  }.
\end{matrix}\right.
\end{equation*}
the coefficients of which depend on complex parameters 
$$
\theta=(\theta_\alpha,\theta_\beta,\theta_\gamma,\theta_\delta).
$$
The main property of this equation is the absence 
of movable singular points, the so-called Painlev\'e property: 
All essential singularities 
of all solutions $q(t)$ of the equation
only appear when $t\in\{0,1,\infty\}$; in other words, 
any solution $q(t)$ extends analytically 
as a meromorphic function on the universal cover of
$\P^1(\C)\setminus\{0,1,\infty\}$. 

Another important property,
expected by Painlev\'e himself, is the irreducibility.
Roughly speaking, the general solution is more transcendental
than solutions of linear, or first order non linear, ordinary differential
equations with polynomial coefficients. Painlev\'e proved
that any irreducible second order polynomial differential
equation without movable singular point falls after reduction
into the $4$-parameter family $P_{VI}$ or one of its degenerations $P_I,\ldots,P_V$.
The fact that Painlev\'e equations are actually irreducible
was proved by Nishioka and Umemura for $P_I$ 
(see \cite{Nishioka:1988,Umemura:1990}) and 
by Watanabe in \cite{Watanabe:1998} for $P_{VI}$. 
Another notion of irreducibility, related with transcendence of 
first integrals, was developped by Malgrange and Casale
in \cite{Malgrange:2001,Casale:2006} and then applied to the 
first of Painlev\'e equations (see \ref{par:Painleve} for more details).
 
A third important property, discovered by R. Fuchs,
is that solutions of $P_{VI}$
parametrize isomonodromic deformations of rank $2$ 
meromorphic connections over the Riemann sphere 
having simple poles at $\{0,t,1,\infty\}$,
with respective set of local exponents 
$(\pm\frac{\theta_\alpha}{2},\pm\frac{\theta_\beta}{2},
\pm\frac{\theta_\gamma}{2},\pm\frac{\theta_\delta}{2})$. 
From this point of view, 
the good space of initial conditions at, say, $t_0$, 
is the moduli space ${\mathcal{M}}_{t_0}(\theta)$ of those connections for $t=t_0$  
(see \cite{IIS:2005}); it turns to be a convenient 
semi-compactification of the naive space of initial conditions 
$\C^2\ni(q(t_0),q'(t_0))$ (compare \cite{Okamoto:1979}).
Via the Riemann-Hilbert correspondence, ${\mathcal{M}}_{t_0}(\theta)$
is analytically isomorphic to the moduli space of corresponding
monodromy representations, namely to (a desingularization of) $S_{(A,B,C,D)}$ with 
parameters
\begin{equation}\label{theta/abcd}
a=2\cos(\pi\theta_\alpha),\ b= 2\cos(\pi\theta_\beta),\  c= 2\cos(\pi\theta_\gamma),\  d= 2\cos(\pi\theta_\delta).
\end{equation}
The (non linear) monodromy of  $P_{VI}$, obtained after analytic continuation around $0$ and $1$ 
of local $P_{VI}$ solutions at $t=t_0$, induces a representation 
$$\pi_1(\P^1(\C)\setminus\{0,1,\infty\},t_0)\to \Aut[S_{(A,B,C,D)}]$$
whose image coincides with the action of $\Gamma_2\subset\PSL(2,\Z)$ (see \cite{Dubrovin-Mazzocco:2000,IIS:2005}).

%%%%%%%%%%%%%%%%%%%%%%%%%%%%%%%%%%%%%%%%%%%%%%%%%%%%%%%%%%%%%%%

\subsection{The Cayley cubic}\label{par:intro_cayley}

%%%%%%%%%%%%%%%%%%%%%%%%%%%%%%%%%%%%%%%%%%%%%%%%%%%%%%%%%%%%%%%

One very specific choice of the parameters will play a central role in this paper. 
The parameters are $(0,0,0,4)$, and the surface $S_{(0,0,0,4)}$ is the unique 
surface in our family with four singularities. Four is the maximal possible number
of isolated singularities for a cubic surface, and $S_{(0,0,0,4)}$ is therefore isomorphic to the
well known Cayley cubic. From the point of view of character varieties, this surface 
appears in the very special case $(a,b,c,d)=(0,0,0,0)$ 
consisting only of solvable representations (dihedral or reducible).

\begin{figure}[t]\label{fig:examples_surfaces}
\input{ch-surf.pstex_t}
\caption{
{\sf{Four examples.}}
{\sf{I.}} The Cayley cubic $S_C$ ; 
{\sf{II.}} $S_{(-0.2,-0.2,-0.2,4.39)}$ ;
 {\sf{III.}} $S_{(0,0,0,3)}$ ;
 {\sf{IV.}}  $S_{(0,0,0,4.1)}.$ 
  }
\end{figure}
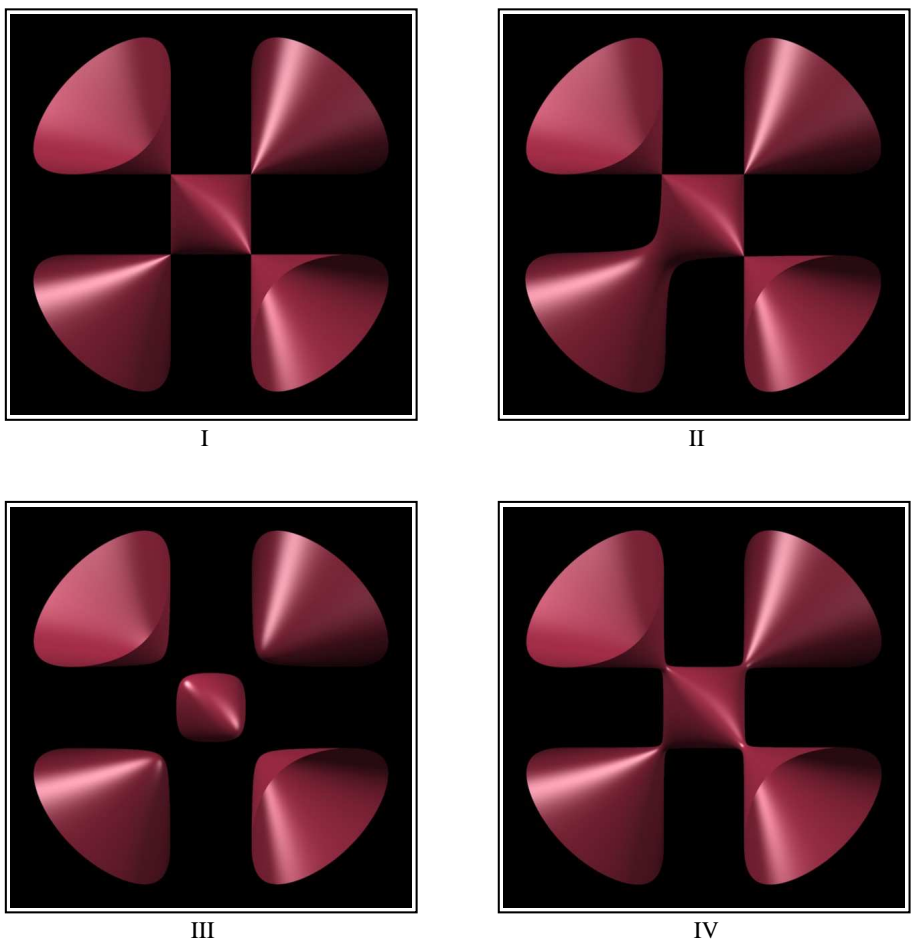

{\sf{I.}} The Cayley cubic $S_C$ ; 
{\sf{II.}} $S_{(-0.2,-0.2,-0.2,4.39)}$ ;
 {\sf{III.}} $S_{(0,0,0,3)}$ ;
 {\sf{IV.}}  $S_{(0,0,0,4.1)}.$

From the Painlev\'e point of view, it corresponds to the
Picard parameter $(\theta_1,\theta_2,\theta_3,\theta_4)=(0,0,0,1).$
The singular foliation which is defined by the corresponding Painlev\'e equation $P_{VI}(0,0,0,1)$
is transversely affine (see \cite{Casale:2006}) and, as was shown by Picard himself, admits
explicit first integrals by means of elliptic functions (see \ref{par:Painleve}). 
Moreover, this specific equation has countably many agebraic solutions, that
are given by finite order points on the Legendre family of elliptic curves  (see \ref{par:Painleve}).

The Cayley cubic has also the ``maximal number of automorphisms'':
The whole group $\PGL(2,\Z),$ in which $\Gamma_2^*$ has index $6,$
stabilizes  the Cayley cubic, and there are extra symmetries coming from 
the permutation of coordinates (see section \ref{par:triangle_automorphisms_area}), so that
the maximal index $24$ of theorem A is obtained in the case of the Cayley cubic. 

Moreover, the degree $2$ orbifold cover
\begin{equation}
\pi_C:\C^*\times\C^*\to S_{(0,0,0,4)}
\end{equation}
 semi-conjugates the action of 
$\PGL(2,\Z)$ on  the character surface 
$S_{(0,0,0,4)}$ to the monomial action of $\GL(2,\Z)$ on $\C^*\times\C^*,$ which is
defined by 
\begin{equation}\label{eq:monomial}
M\left( \begin{array}{c} 
         u \\ v
        \end{array}
\right)
=
\left(\begin{array}{c} 
         u^{m_{11}}v^{m_{12}} \\ u^{m_{21}}v^{m_{22}}
        \end{array}
\right),
\end{equation}
for any element $M$ of $\GL(2,\Z).$ 
On the universal cover $\C\times\CÊ\to\C^*\times\C^*$,
the lifted dynamics is the usual affine action of the group $\GL(2,\Z)\ltimes\Z^2$ on 
the complex plane $\C^2.$  

%%%%%%%%%%%%%%%%%%%%%%%%%%%%%%%%%%%%%%%%%%%%%%%%%%%%%%%%%%%%%%%%

\subsection{Compactification and entropy}

%%%%%%%%%%%%%%%%%%%%%%%%%%%%%%%%%%%%%%%%%%%%%%%%%%%%%%%%%%%%%%%

Our first goal is to classify automorphisms of surfaces $S_{(A,B,C,D)}$ in three
types, elliptic, parabolic and hyperbolic, and to describe the main properties of
the dynamics of each type of automorphisms. This classification is compatible
with the description of mapping classes, Dehn twists corresponding to parabolic
transformations, and pseudo-Anosov mappings to hyperbolic automorphisms. 
The most striking result in that direction is summarized in the following theorem.

\begin{thm-B}
Let $A,$ $B,$ $C,$ and $D$ be four complex numbers. Let $M$ be an element of 
$\Gamma_2^*,$ and $f_M$ be the automorphism of $S_{(A,B,C,D)}$ which is
determined by $M.$ The topological entropy of $f_M:S_{(A,B,C,D)}(\C)\to 
S_{(A,B,C,D)}(\C)$ is equal to the logarithm of the spectral radius of $M.$ 
\end{thm-B}
  
The proof is obtained by a deformation argument: We 
shall show that the topological entropy does not depend on the parameters 
$(A,B,C,D),$ and then compute it in the case of the Cayley cubic. To do so, 
we first describe the geometry of surfaces $S \in \SS$ (section 
\ref{par:geometry-surfaces}), their groups of automorphisms
(section \ref{par:automorphisms-surfaces}), and the action of automorphisms 
by birational transformations on the Zariski closure ${\overline{S}}$ of $S$ in $\P^3(\C)$ 
(section \ref{par:bira-extension-dyna}). 

Another algorithm to compute the topological entropy has been obtained by 
Iwasaki and Uehara for non singular cubics $S$ in $\SS$ (see \cite{Iwasaki-Uehara:2006}). The case of singular cubics is crucial for the study of the set of quasi-fuchsian 
deformations of fuchsian representations, in connection with Bers embedding of
Teichm\"uller spaces (see \cite{Goldman:survey} and \cite{Cantat:BHPS}).

%%%%%%%%%%%%%%%%%%%%%%%%%%%%%%%%%%%%%%%%%%%%%%%%%%%%%%%%%%%%%%%%

\subsection{Bounded orbits}

%%%%%%%%%%%%%%%%%%%%%%%%%%%%%%%%%%%%%%%%%%%%%%%%%%%%%%%%%%%%%%%

Section \ref{par:sectionbounded} is devoted to the study  of parabolic elements 
(or Dehn twists), and bounded or periodic (i.e. finite) orbits 
of $\Gamma_2^*$. 
For instance, given a representation $\rho:\Sphere_4\to\SU(2)\subset\SL(2,\C)$, 
the $\Gamma_2^*$-orbit of the correponding point 
$\chi(\rho)$ will be bounded, contained in the cube $[-2,2]^3$.
If moreover the image of $\rho$ is finite, then so will be the 
corresponding orbit. Though, there are  
periodic orbits with complex coordinates.

First of all, fixed points of $\Gamma_2^*$ are precisely  
the singular points of $S_{(A,B,C,D)}$ and have been 
extensively studied (see \cite{IIS:2005}).
Singular points arise
from semi-stable points of $\Rep(\Sphere_4)$, that is to say
either from reducible representations, or from those representations 
for which one of the matrices $\rho(\alpha)$, $\rho(\beta)$, 
$\rho(\gamma)$ or $\rho(\delta)$ is $\pm I$. Both type
of degeneracy occur at each singular point of $S_{(A,B,C,D)}$ depending on the choice of parameters $(a,b,c,d)$ fitting 
to $(A,B,C,D)$. The Riemann-Hilbert
correspondance ${\mathcal{M}}_{t_0}(\theta)\to S_{(A,B,C,D)}$
is a minimal resolution of singularities and $P_{VI}$ equation
restricts to the exceptional divisor as a Riccati equation :
this is the locus of Riccati-type solutions. We note that 
any point $(x,y,z)$ is the singular point of one member
$S_{(A,B,C,D)}$.

Periodic orbits of length $\ge2$ correspond to algebraic
solutions of $P_{VI}$ equation (see \cite{Iwasaki:2007}).
In Proposition \ref{Prop:SmallOrbits}, we classify orbits
of length $\le4$ : we find one $2$-parameter family of 
length $2$ orbits and two $1$-parameter families of 
length $3$ and $4$ orbits. They correspond to 
well-known\footnote{Although we usually find $4$ families 
of algebraic solutions of $P_{VI}$ in the litterature (see \cite{Boalch:2007,BHG}),
there are actually $3$ up to Okamoto symmetries :
degree $4$ solutions 3B and 4C in \cite{BHG} are conjugated
by the symmetry $s_1s_2s_1$ (with notations of \cite{Noumi-Yamada}).}
algebraic solutions of $P_{VI}$ equation (see \cite{BHG}).
For instance, the length $2$ orbit arise when $A=C=0$ ;
the corresponding $P_{VI}$-solution is $q(t)=\sqrt{t}$.

The following result shows that infinite bounded orbits are
real and contained in the cube $[-2,2]^3$. 

\begin{thm-C}
Let $m$ be a point of $S_{(A,B,C,D)}$ with a bounded  
$\Gamma_2^*$-orbit of length $>4$. 
Then, the parameters $(A,B,C,D)$
are real numbers and the orbit is contained in the real part
$S_{(A,B,C,D)}(\R)$ of the surface. 

If the orbit of $m$ is finite, then both the surface and the orbit
are actually defined over a (real) number field.

If the orbit of $m$ is infinite, then it corresponds to a 
$\SU(2)$-representation for a convenient choice 
of parameters $(a,b,c,d)$, and the orbit is contained
and dense in the unique bounded connected component 
of the smooth part of $S_{(A,B,C,D)}(\R)$.
\end{thm-C}

As a corollary, periodic orbits of length $>4$ are rigid
and we recover the main result of \cite{BHG}.
Recall that Cayley cubic contains infinitely many
periodic orbits, of arbitrary large order. It is conjectured
that there are finitely many periodic orbits apart
from the Cayley member, but this is still an open problem.
A classification of known periodic orbits can be found in 
\cite{Boalch:2007}.

About infinite orbits, Theorem C should be compared with
results of Goldman and Previte and Xia, 
concerning the dynamics on the character variety
for representations into $\SU(2)$ \cite{Previte-Xia:2005}\footnote{After reading the first version of this paper, professor Iwasaki informed us that part of theorem C was already announced in \cite{Iwasaki:bounded}.}.
We note that an infinite bounded orbit may also
correspond to $\SL(2,\R)$-representation for an alternate
choice of parameters $(a,b,c,d)$.  

This theorem stresses the particular role played by the real case, when all the parameters $A$, $B$,
$C$, and $D$ are real numbers; in that case,  $\Gamma_2^*$ preserves the real part of the surface
and we have two different, but closely related, dynamical systems: The action on the complex
surface $S_{(A,B,C,D)}(\C)$ and the action on the real surface $S_{(A,B,C,D)}(\R)$. The link
between those two dynamical systems will be studied in \cite{Cantat:BHPS}.

%%%%%%%%%%%%%%%%%%%%%%%%%%%%%%%%%%%%%%%%%%%%%%%%%%%%%%%%%%%%%%%%

\subsection{Dynamics, affine structures, and the irreducibility of $P_{VI}$}

%%%%%%%%%%%%%%%%%%%%%%%%%%%%%%%%%%%%%%%%%%%%%%%%%%%%%%%%%%%%%%%

The last main result that we shall prove concerns the
classification of parameters $(A,B,C,D)$ for which $S_{(A,B,C,D)}$ admits a 
$\Gamma_2^*$-invariant holomorphic geometric structure.

\begin{thm-D}\label{thm:main_affine}
The group $\Gamma_2^*$ does not preserve any holomorphic curve of finite type, any singular holomorphic
foliation, or any singular holomorphic web. 
The group $\Gamma_2^*$  does not preserve any meromorphic affine structure, except in the case of
the Cayley cubic, {\sl{i.e.}} when $(A,B,C,D)=(0,0,0,4),$ or equivalently when
$$
(a,b,c,d)=  (0,0,0,0) \ \text{or} \ (2,2,2,-2),$$
up to multiplication by $-1$ and permutation of the parameters.
\end{thm-D}

Following \cite{Casale:2005}, the same strategy shows that the Galois groupoid
is the whole symplectic pseudo-group except in the Cayley case (see section \ref{par:Painleve}),  and
we get

\begin{thm-E}
The sixth Painlev\'e equation is irreducible in the sense of Malgrange and 
Casale except when 
$(A,B,C,D)=(0,0,0,4)$, {\sl{i.e.}} except in one of the
following cases:
\begin{itemize}
\item $\theta_\omega\in\frac{1}{2}+\Z$, $\forall \, \omega=\alpha,\beta,\gamma,\delta$,
\item $\theta_\omega\in \Z$, $\forall \,\omega=\alpha,\beta,\gamma,\delta$, and $\sum_\omega \theta_\omega$ is even.
\end{itemize}
\end{thm-E}

Following \cite{Casale:2007}, Malgrange-Casale irreducibility also
implies Nishioka-\-Ume\-mura irreducibility, so that theorem \ref{thm:main_affine} indeed
provides a galoisian proof of the irreducibility in the spirit of Drach and Painlev\'e.

%%%%%%%%%%%%%%%%%%%%%%%%%%%%%%%%%%%%%%%%%%%%%%%%%%%%%%%%%%%%%%%%

\subsection{Aknowledgement}

%%%%%%%%%%%%%%%%%%%%%%%%%%%%%%%%%%%%%%%%%%%%%%%%%%%%%%%%%%%%%%%
\vv

This article has been written while the first author was visiting Cornell 
University: Thanks to Cliff Earle, John Smillie and Karen Vogtmann
for nice discussions concerning this paper, and to the DREI for 
travel expenses between Rennes and Ithaca.

We would like to kindly thank Marta Mazzocco who introduced 
us to Painlev\'e VI equation,  its geometry and dynamics.
The talk she gave in Rennes was the starting point of our work.
Many thanks also to Guy Casale who taught us about irreducibility, 
and to Yair Minsky who kindly explained some aspects of character
varieties to the first author.
 
Part of this paper was the subject of
a conference held in Rennes in 2005, which was funded by the ANR 
project "Syst\`emes dynamiques polynomiaux", and both authors are now
taking part to the ANR project "Symplexe."

%%%%%%%%%%%%%%%%%%%%%%%%%%%%%%
%%%%%%%%%%%%%%%%%%%%%%%%%%%%%%
\section{The family of surfaces}\label{par:geometry-surfaces}

As explained in \ref{par:char_var_nota}, we shall consider the family $\SS$ of complex 
affine surfaces which are defined by the following type of cubic equations 
$$
x^2+y^2+z^2+x y z = A x + By +C z + D,
$$
in which $A$, $B$, $C$, and $D$ are four complex parameters. 
Each choice of $(A,B,C,D)$
gives rise to one surface $S$ in our family; if necessary, 
$S$ will also be denoted $S_{(A,B,C,D)}$.
When the parameters are real numbers, $S(\R)$ will denote the real part of $S$.
Figure \ref{fig:examples_surfaces}  presents a few pictures of $S(\R)$ for various choices of the parameters. 

This section contains preliminary results on the geometry of the surfaces 
$S_{(A,B,C,D)}$, and
the automorphisms of these surfaces. Most of these results are well known to 
algebraic geometers and specialists of Painlev\'e VI equations. 
%%%%%%%%%%%%%%%%%%%%%%%%%%%%
\subsection{The Cayley cubic}\label{par:Cayley} 

In 1869, Cayley proved that,  up to  projective transformations,
there is a unique cubic surface in $\P^3(\C)$ with four isolated singularities. 
One of the nicest models of the Cayley cubic is the surface $S_{(0,0,0,4)},$
whose equation is
$$
x^2+y^2+z^2+xyz=4.
$$ 
The four singular points of $S_C$ are rational nodes located at
$$
(-2,-2,-2),\ (-2,2,2),\ (2,-2,2)\ \ \ \text{and}\ \ \ (2,2,-2),
$$
 and can be seen on figure \ref{fig:examples_surfaces}. This specific member of our family of surfaces 
 will be called {\sl{the Cayley cubic}} and denoted $S_C.$ 
 This is justified by the following theorem (see Appendix A).
 
\begin{thm}[Cayley]\label{thm:carac_cayley}
If $S$ is a member of the family $\SS$ with four singular points, then 
$S$ coincides with the Cayley cubic $S_C$. 
\end{thm}

The  Cayley cubic is isomorphic to the quotient of $\C^*\times \C^*$ by the involution 
$\eta(u,v)=(u^{-1},v^{-1})$. The map 
$$
\pi_C(x,y) = \left( -u-\frac{1}{u},-v-\frac{1}{v}, -uv - \frac{1}{uv}\right)
$$
gives an explicit isomorphism between $(\C^*\times \C^*)/\eta$ and $S_C.$ The four
fixed points 
$$
(1,1),\ \ \ (1,-1),\ \ \ (-1,1)\ \ \ \text{and}\ \ \ (-1,-1)
$$ of $\eta$ respectively correspond to the singular points of $S_C$ above.

The real surface $S_C(\R)$ contains the four singularities of $S_C,$ and the smooth 
locus $S_C(\R)\setminus \Sing(S_C)$ is made of five components : A bounded
one, the closure of which coincides with the image of 
$\T=\mathbb S^1\times\mathbb S^1\subset\C^*\times\C^*$ by $\pi_C,$ 
and four unbounded ones, corresponding to images of $\R^+\times \R^+,$
$\R^+\times \R^-,$ $\R^-\times \R^+,$ and $\R^-\times \R^-$ (see figure \ref{fig:examples_surfaces}).

As explained in section \ref{par:intro_cayley}, the group $\GL(2,\Z)$ acts on $\C^*\times \C^*$
by monomial transformations, and this action commutes with the involution $\eta$,
permuting its fixed points. As a consequence, $\PGL(2,\Z)$ acts on the quotient $S_C$. 
Precisely, the generators 
$$
\begin{pmatrix}1&0\\-1&1\end{pmatrix},\ \ \ \begin{pmatrix}1&1\\0&1\end{pmatrix}\ \ \ 
\text{and}\ \ \ \begin{pmatrix}1&0\\0&-1\end{pmatrix}
$$
of $\PGL(2,\Z)$ respectively send the triple $(x,y,z)$ to
$$
(x,-z-xy,y),\ \ \ (z,y,-x-yz)\ \ \ \text{and}\ \ \ (x,y,-z-xy).
$$
As we shall see, the induced action of $\PGL(2,\Z)$ on $S_C$ coincides with the action 
of the extended mapping class group of $\Sphere_4$ considered
in \S \ref{par:automorphisms_modular}.

The group $\PGL(2,\Z)$ preserves the real part of $S_C$ ;
for example, the product $\C^*\times\C^*$ retracts by deformation 
on the real $2$-torus 
$\T=\mathbb S^1\times\mathbb S^1,$
and the monomial action of $\GL(2,\Z)$  preserves this torus
(it is the standard one under the parametrization
$(s,t)\mapsto(e^{2i\pi s},e^{2i\pi t})$).
 
%%%%%%%%%%%%%%%%%%%%%%%%%%%%%%%%%%%%%%%%%%

\subsection{Mapping class group action}\label{par:MappingClassGroup}

%%%%%%%%%%%%%%%%%%%%%%%%%%%%%%%%%%%%%%%%

First, let us detail section \ref{par:automorphisms_modular}.
The extended mapping class group $\MCG^*(\Sphere_4)$ 
is the group of isotopy classes of homeomorphisms
of the four punctured sphere $\Sphere_4$; 
the usual mapping class group $\MCG(\Sphere_4)$ 
is the index $2$ subgroup consisting only in orientation
preserving homeomorphisms. Those groups embed 
in the group of outer automorphisms of $\pi_1(\Sphere_4)$ 
in the following way. Fix a base
point $p_0\in\Sphere_4$.  In any isotopy class, 
one can find a homeomorphism $h$ fixing $p_0$ and thus
inducing an automorphism of the fundamental group
$$
h_*:\pi_1(\Sphere_4,p_0)\to\pi_1(\Sphere_4,p_0)\ ;\ 
\gamma\mapsto h\circ\gamma.
$$
The class of $h_*$ modulo inner automorphisms does not depend 
on the choice of the representative $h$ in the homotopy class
and we get a morphism 
\begin{equation}
\MCG^*(\Sphere_4)\to\Out(\pi_1(\Sphere_4))
\end{equation} 
which turns out to be an isomorphism.

Now, the action of $\Out(\pi_1(\Sphere_4))$ on $\chi(\Sphere_4)$ gives rise
to a morphism
\begin{equation}
\left\{\begin{matrix}\MCG^*(\Sphere_4)&\to&\Aut[\chi(\Sphere_4)]\\
[h]&\mapsto&\left\{\chi(\rho)\mapsto\chi(\rho\circ h^{-1})\right\}
\end{matrix} \right.
\end{equation}
into the group of polynomial diffeomorphisms of $\chi(\Sphere_4).$
(here, we use that $\rho\circ(h_*)^{-1}=\rho\circ h^*=\rho\circ h^{-1}$).
Our goal in this section is to give explicit formulae for this action of $\MCG^*(\Sphere_4)$
on $\chi(\Sphere_4),$ and to describe the subgroup  of $\MCG(\Sphere_4)$ which 
stabilizes each surface $S_{(A,B,C,D)}.$ 

%%%%%%%%%
\subsubsection{Torus cover}\label{par:toruscover}
%%%%%%%%%

Consider the two-fold ramified cover 
\begin{equation}
\pi_T:\T=\R^2/\Z^2\to\Sphere
\end{equation}
with Galois involution $\sigma:(x,y)\mapsto(-x,-y)$ sending its
ramification points $(1/2,0)$,
$(0,1/2)$, $(1/2,1/2)$ and $(0, 0)$
respectively to the four punctures $p_\alpha$, $p_\beta$, $p_\gamma$ and $p_\delta$
(see figure \ref{fig:four_punctured_sphere}).

\begin{figure}[h]\label{fig:torus}
\input{torus.pstex_t}
\caption{ {\sf{The torus cover. }}}
\end{figure}
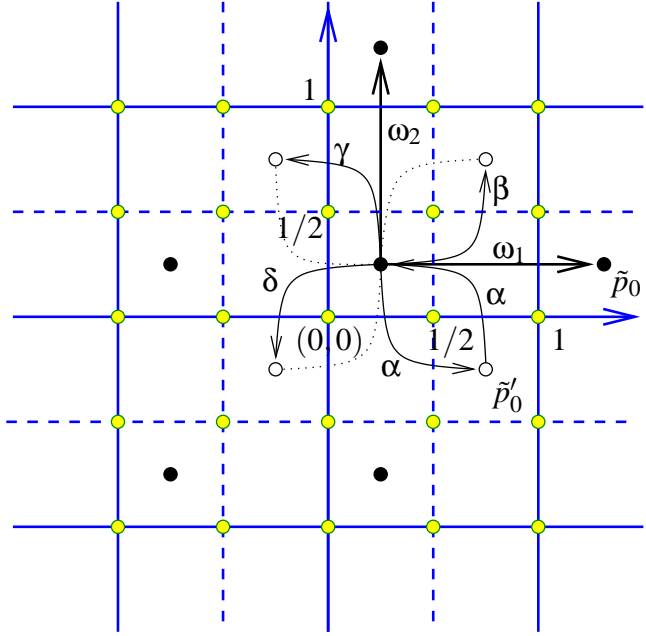

The mapping class group of the torus, and also of the once punctured torus $\T_1=\T\setminus\{(0,0)\},$ 
is isomorphic to $\GL(2,\Z).$ This group acts by linear homeomorphisms on the torus, fixing $(0, 0),$ 
and permuting the other three ramification points of $\pi_T.$ This action provides a section of the
projection $\Diff(\T)\to \MCG^*(\T).$ Since this action commutes with the involution $\sigma$
(which generates the center of $\GL(2,\Z)$), we get a morphism from $\PGL(2,\Z)$ to 
$\MCG^*(\Sphere_4).$  This morphism is one to one and its image is contained in 
the stabilizer of $p_\delta$ in $\MCG^*(\Sphere_4)$. 

The subset  $H\subset\T$ of ramification points of $\pi$ 
coincides with the $2$-torsion subgroup of $(\T,+)$ ; $H$ 
acts by translation on $\T$ and commutes with the involution 
$\sigma$ as well. This provides an isomorphism (see section 4.4 in \cite{Birman:book})
\begin{equation}
\PGL(2,\Z)\ltimes H\to\MCG^*(\Sphere_4).
\end{equation}

\begin{lem}\label{lem:formulae-mcg} The subgroup of $\Aut(\chi(\Sphere_4))$ obtained by the action of the subgroup $\PGL(2,\Z)$ of $\MCG^*(\Sphere_4)$
is generated by the three polynomial automorphisms $B_1,$ $B_2$ and $T_3$ of equations \ref{eq:b1}, \ref{eq:b2}, 
and \ref{eq:t3} below.  The $4$-order translation group $H$ acts trivially on parameters $(A,B,C,D,x,y,z),$ permuting parameters
$(a,b,c,d)$ as follows 
\begin{eqnarray}\label{eq:p1p2}
P_1=(1/2,0) :  (a,b,c,d)\mapsto (d,c,b,a) \\
P_2=(0,1/2) :  (a,b,c,d)\mapsto (b,a,d,c) 
\end{eqnarray}
\end{lem}

\begin{proof}
Let $\tilde p_0$ and $\tilde p_0'$ be the lifts of the base point
$p_0\in\Sphere_4$. Still denote by $\alpha$, $\beta$, $\gamma$ 
and $\delta$ the two lifts of those loops, with respective initial points 
$\tilde p_0$ and $\tilde p_0'$. The fundamental group of the four punctured torus 
$\T_4=\T\setminus H$ based at $\tilde p_0$ may
be viewed as the set of even words in $\alpha$, $\beta$, $\gamma$ 
and $\delta$, or equivalently of words in $\omega_1$, $\omega_2$
and $\delta$ that are even in $\delta$ where
$$
\omega_1=\beta\gamma=\alpha^{-1}\delta^{-1}\ \ \ \text{and}\ \ \ \omega_2=\gamma\delta=\beta^{-1}\alpha^{-1}.
$$
(see Figure \ref{fig:torus}). The action of the linear homeomorphism 
$$
B_1=\begin{pmatrix}1&0\\-1&1\end{pmatrix}:\T_4\to\T_4,
$$ 
or we should say, of a convenient isotopic homeomorphism $h$
fixing $\tilde p_0$, on the fundamental groups $\pi_1(\T_4,\tilde p_0)$
 and $\pi_1(\Sphere_4,p_0)$ is given by :
$$
h_*:\ \ \ 
\left\{\begin{matrix}
\omega_1 & \mapsto & \delta^{-1}\omega_1^{-1}\omega_2\delta^{-1}\\
\omega_2 & \mapsto & \omega_2\\
\delta & \mapsto & \delta
\end{matrix}\right.
\ \ \ \text{i.e.}\ \ \ 
\left\{\begin{matrix}
\alpha & \mapsto & \alpha\beta\alpha^{-1}\\
\beta & \mapsto & \alpha\\
\gamma & \mapsto & \gamma\\
\delta & \mapsto & \delta
\end{matrix}\right.
$$
This automorphism of  $\pi_1(\Sphere_4,p_0)$, which depends on the choice of $h$ in 
the isotopy class of $B_1$, induces an automorphism
$$
\Rep(\Sphere_4)\to\Rep(\Sphere_4)\ ;\ \rho\mapsto\rho\circ(h_*)^{-1}.
$$
The corresponding action on the character variety, 
{\sl{i.e.}} on the corresponding $7$-uples $(a,b,c,d,x,y,z)\in\C^7$,
is independant of that choice. In order to compute it, note
that 
$$
(h_*)^{-1}=h^*:\ \ \ 
\left\{\begin{matrix}
\alpha & \mapsto & \beta\\
\beta & \mapsto & \beta^{-1}\alpha\beta\\
\gamma & \mapsto & \gamma\\
\delta & \mapsto & \delta
\end{matrix}\right.
$$
We therefore obtain
\begin{equation}\label{eq:b1}
B_1=\begin{pmatrix}1&0\\-1&1\end{pmatrix}\ :\ \ \ 
\left\{\begin{matrix}
a & \mapsto & b\\
b & \mapsto & a\\
c & \mapsto & c\\
d & \mapsto & d
\end{matrix}\right.
\ \ \ \text{and}\ \ \ 
\left\{\begin{matrix}
x & \mapsto & x\hfill\\
y & \mapsto & -z-xy+ac+bd\\
z & \mapsto & y\hfill
\end{matrix}\right.
\end{equation}
For instance, the coordinate $y'$ of the image is given by $y'=\tr(\rho\circ h^*(\beta\gamma))=\tr(\rho(\beta^{-1}\alpha\beta\gamma)),$ 
and its value is  easily computed using standard Fricke-Klein formulae, like
$$
\tr(M_1)=\tr(M_1^{-1}),\ \ \ \tr(M_1M_2)=\tr(M_2M_1),
$$
$$
\tr(M_1M_2^{-1})+\tr(M_1M_2)=\tr(M_1)\tr(M_2)
$$
$$
\text{and}\ \ \ \tr(M_1M_2M_3)+\tr(M_1M_3M_2)+\tr(M_1)\tr(M_2)\tr(M_3)
$$
$$
=\tr(M_1)\tr(M_2M_3)+\tr(M_2)\tr(M_1M_3)+\tr(M_3)\tr(M_1M_2)
$$
for any $M_1,M_2,M_3\in\SL(2,\C)$.

A similar computation yields 
\begin{equation}\label{eq:b2}
B_2=\begin{pmatrix}1&1\\0&1\end{pmatrix}\ :\ \ \ 
\left\{\begin{matrix}
a & \mapsto & a\\
b & \mapsto & c\\
c & \mapsto & b\\
d & \mapsto & d
\end{matrix}\right.
\ \ \ \text{and}\ \ \ 
\left\{\begin{matrix}
x & \mapsto & z\hfill\\
y & \mapsto & y\\
z & \mapsto & -x-yz+ab+cd\hfill
\end{matrix}\right.
\end{equation}
which, together with $B_1$, provide a system of generators
for the $\PSL(2,\Z)$-action.
In order to generate $\PGL(2,\Z)$, we have to add the involution
\begin{equation}\label{eq:t3}
T_3=\begin{pmatrix}0&1\\1&0\end{pmatrix}\ :\ \ \ 
\left\{\begin{matrix}
a & \mapsto & c\\
b & \mapsto & b\\
c & \mapsto & a\\
d & \mapsto & d
\end{matrix}\right.
\ \ \ \text{and}\ \ \ 
\left\{\begin{matrix}
x & \mapsto & y\\
y & \mapsto & x\\
z & \mapsto & z
\end{matrix}\right.
\end{equation}
The formulae for the action of $H$ are obtained in the same way.
\end{proof}

\begin{rem}
The formulae \ref{eq:b1}, \ref{eq:b2}, 
and \ref{eq:t3} for $B_1,$ $B_2$ and $T_3$ specialize to the
formulae of section \ref{par:Cayley} when $(A,B,C,D)=(0,0,0,4).$
\end{rem}

\begin{rem}\label{rem:artin}
The Artin Braid Group 
${\mathcal{B}}_3=\langle \beta_1,\beta_2\ \vert\ \beta_1\beta_2\beta_1=\beta_2\beta_1\beta_2\rangle$
is isomorphic to the group of isotopy classes of the thrice punctured disk fixing its boundary. There is
therefore a morphism from ${\mathcal{B}}$ into the subgroup of $\MCG(\Sphere_4)$ that stabilizes $p_\delta.$ This
morphism gives rise to the following well known exact sequence 
$$
I\to\langle (\beta_1\beta_2)^3\rangle\to \mathcal B_3\to \PSL(2,\Z)\to 1,
$$
where generators $\beta_1$ and $\beta_2$ are respectively sent
to $B_1$ and $B_2,$ and the group $\langle (\beta_1\beta_2)^3\rangle$  coincides
with the center of ${\mathcal{B}}_3.$ In particular, the action of ${\mathcal{B}}_3$ on $\chi(\Sphere_4)$
coincides with the action of $\PSL(2,\Z).$
We note that $\PSL(2,\Z)$ is the free product of the trivolution 
$B_1B_2$ and the involution $B_1B_2B_1$. In $\PGL(2,\Z)$,
we also have relations $T_3^2=I$, $T_3B_1T_3=B_2^{-1}$ and $T_3B_2T_3=B_1^{-1}$.
\end{rem}

%%%%%%%%%
\subsubsection{The modular groups $\Gamma_2^*$ and $\Gamma_2$}\label{par:modularformulae}
%%%%%%%%%

Since the action of $M\in\GL(2,\Z)$ on the set $H$  of  points of order $2$
depends only on the equivalence class of $M$ modulo $2$, we get
an exact sequence
$$
I\to\Gamma_2^*\to \PGL(2,\Z)\ltimes H\to \Sym_4\to 1
$$
where $\Gamma_2^*\subset\PGL(2,\Z)$ is the subgroup
defined by those matrices $M\equiv I$ modulo $2$.
This group acts on the character variety, and since it preserves
the punctures, it fixes $a,$ $b,$ $c,$ and $d.$ The group $\Gamma_2^*$ is the free 
product of  $3$ involutions, $s_x,$ $s_y,$ and
$s_z,$ acting on the character variety as follows.
\begin{equation}\label{eq:formula_involution_x}
s_x=\begin{pmatrix}-1&2\\ 0&1\end{pmatrix}\ :\ \ \ 
\left\{\begin{matrix}
x & \mapsto & -x-yz+ab+cd\\
y & \mapsto & y\hfill\\
z & \mapsto & z\hfill
\end{matrix}\right.
\end{equation}
\begin{equation}\label{eq:formula_involution_y}
s_y=\begin{pmatrix}1&0\\ 2&-1\end{pmatrix}\ :\ \ \ 
\left\{\begin{matrix}
x & \mapsto & x\hfill\\
y & \mapsto & -y-xz+bc+ad\\
z & \mapsto & z\hfill
\end{matrix}\right.
\end{equation}
\begin{equation}\label{eq:formula_involution_z}
s_z=\begin{pmatrix}1&0\\0&-1\end{pmatrix}\ :\ \ \ 
\left\{\begin{matrix}
x & \mapsto & x\hfill\\
y & \mapsto & y\hfill\\
z & \mapsto & -z-xy+ac+bd
\end{matrix}\right.
\end{equation}
We note that $s_x=B_2B_1^{-1}B_2^{-1}T_3$, 
$s_y=B_2B_1B_2^{-1}T_3$ and $s_z=B_2B_1B_2T_3$.
The standard modular group $\Gamma_2\subset\PSL(2,\Z)$
is generated by 
$$
\left\{\begin{matrix}
g_x&=&s_z s_y&=&B_1^{2}&=&\begin{pmatrix}1&0\\-2&1\end{pmatrix}\\
g_y&=&s_x s_z&=&B_2^{2}&=&\begin{pmatrix}1&2\\0&1\end{pmatrix}\\
g_z&=&s_y s_x&=&B_1^{-2}B_2^{-2}&=&\begin{pmatrix}1&-2\\2&-3\end{pmatrix}
\end{matrix}\right.
$$
(we have $g_z g_y g_x=I$); as we shall see, this corresponds to Painlev\'e VI monodromy (see \cite{IIS:2005}
and section \ref{par:Painleve}). The following proposition is now a direct consequence of lemma \ref{lem:formulae-mcg}.

\begin{pro}\label{pro:mcg}
Let $\MCG_0^*(\Sphere_4)$ (resp. $\MCG_0(\Sphere_4)$) be the subgroup of $\MCG^*(\Sphere_4)$
(resp. $\MCG(\Sphere_4)$) which stabilizes the four punctures of $\Sphere_4.$ This group coincides
with the stabilizer of the projection  $\pi:\chi(\Sphere_4)\to \C^4$ which is defined by
$$
\pi(a,b,c,d,x,y,z)=(a,b,c,d).
$$
Its image in $\Aut(\chi(\Sphere_4))$ coincides with the image of $\Gamma_2^*$ (resp. $\Gamma_2$)
and is therefore generated by the three involutions $s_x,$ $s_y$ and $s_z$ (resp. the three automorphisms
$g_x,$ $g_y,$ $g_z$).
\end{pro}

As we shall see in sections \ref{par:triangle_automorphisms_area} and
\ref{par:consequences_notations}, this group is of finite index in 
$\Aut(\chi(\Sphere_4)$. 

\begin{rem}
Let us consider the exact sequence
$$
I\to\Gamma_2^*\to \PGL(2,\Z)\to \Sym_3\to 1,
$$
where $\Sym_3\subset \Sym_4$ is the stabilizer of $p_\delta$, or equivalently of $d$, or $D.$
A splitting $\Sym_3\hookrightarrow\PGL(2,\Z)$ is generated by the  transpositions
$T_1=T_3B_1B_2$ and $T_2=B_1B_2T_3.$ 
They act as follows on the 
character variety. 
$$
T_1=\begin{pmatrix}-1&0\\1&1\end{pmatrix}\ :\ \ \ 
\left\{\begin{matrix}
a & \mapsto & b\\
b & \mapsto & a\\
c & \mapsto & c\\
d & \mapsto & d
\end{matrix}\right.
\ \ \ \text{and}\ \ \ 
\left\{\begin{matrix}
x & \mapsto & x\\
y & \mapsto & z\\
z & \mapsto & y
\end{matrix}\right.
$$
and
$$
T_2=\begin{pmatrix}1&1\\0&-1\end{pmatrix}\ :\ \ \ 
\left\{\begin{matrix}
a & \mapsto & a\\
b & \mapsto & c\\
c & \mapsto & b\\
d & \mapsto & d
\end{matrix}\right.
\ \ \ \text{and}\ \ \ 
\left\{\begin{matrix}
x & \mapsto & z\\
y & \mapsto & y\\
z & \mapsto & x
\end{matrix}\right.
$$
\end{rem}

%%%%%%%%%%%%%%%%%%%%%%%%%%%%%%%%%%%%%%%%

\subsection{Twists}\label{par:twists}

%%%%%%%%%%%%%%%%%%%%%%%%%%%%%%%%%%%%%%%%

There are other symmetries between surfaces $S_{(A,B,C,D)}$
that do not arise from the action of the mapping class group.
Indeed, given any $4$-uple $\epsilon=(\epsilon_1,\epsilon_2,\epsilon_3,\epsilon_4)\in\{\pm1\}^4$ with $\prod_{i=1}^4\epsilon_i=1$,
the $\epsilon$-twist of a representation $\rho\in\Rep(\Sphere_4)$
is the new representation $\otimes_\epsilon\rho$
generated by 
$$
\left\{\begin{matrix}
\tilde\rho(\alpha) & = & \epsilon_1\rho(\alpha)\\
\tilde\rho(\beta) & = & \epsilon_2\rho(\beta)\\
\tilde\rho(\gamma) & = & \epsilon_3\rho(\gamma)\\
\tilde\rho(\delta) & = & \epsilon_4\rho(\delta)
\end{matrix}\right.
$$
This provides an action of $\Z/2\Z\times \Z/2\Z\times \Z/2\Z$
on the character variety given by
$$
\otimes_\epsilon:\ \ \ 
\left\{\begin{matrix}
a & \mapsto & \epsilon_1a\\
b & \mapsto & \epsilon_2c\\
c & \mapsto & \epsilon_3b\\
d & \mapsto & \epsilon_4d
\end{matrix}\right.
\ \ \ 
\left\{\begin{matrix}
A & \mapsto & \epsilon_1\epsilon_2A\\
B & \mapsto & \epsilon_2\epsilon_3B\\
C & \mapsto & \epsilon_1\epsilon_3C\\
D & \mapsto & \hfill D
\end{matrix}\right.
\ \ \ \text{and}\ \ \ 
\left\{\begin{matrix}
x & \mapsto & \epsilon_1\epsilon_2x\\
y & \mapsto & \epsilon_2\epsilon_3y\\
z & \mapsto & \epsilon_1\epsilon_3z
\end{matrix}\right.
$$
The action on $(A,B,C,D,x,y,z)$ is trivial iff $\epsilon=\pm(1,1,1,1)$.
The "Benedetto-Goldman symmetry group" of order 192 acting on $(a,b,c,d,x,y,z)$
which is described in \cite{Benedetto-Goldman:1999} (\S 3C) 
is precisely the group generated by $\epsilon$-twists
and the symmetric group $\Sym_4=\langle T_1,T_2,P_1,P_2\rangle$.
The subgroup $Q$ acting trivially on $(A,B,C,D,x,y,z)$ is 
of order $8$ generated by
\begin{equation}
Q=\langle P_1,P_2,\otimes_{(-1,-1,-1,-1)}\rangle.
\end{equation} 

%%%%%%%%%%
\subsection{Character variety of the once-punctured torus}\label{par:OncePuncturedTorus}
%%%%%%%%%%
Our family of surfaces $S_{(A,B,C,D)}$ also provides,
for $(A,B,C,D)=(0,0,0,D)$, the moduli
space of representations of the torus $\T=\R^2/\Z^2$ with
one puncture at $(0,0)$. Precisely, if we go back to the notations 
of \S \ref{par:MappingClassGroup} (see figure \ref{fig:torus}),
the fundamental group $\pi_1(\T_1)$, $\T_1=\T\setminus\{(0,0)\}$,
is the free group generated by $\omega_1$ and $\omega_2$.
The algebraic quotient $\chi(\T_1)=\Rep(\T_1)/\! /\SL(2,\C)$
is given by the map 
$$
\left\{\begin{matrix}\Rep(\T_1)&\to&\chi(\T_1)\simeq\C^3\hfill\\
\rho&\mapsto&(X,Y,Z)=(\tr(\rho(\omega_1)),\tr(\rho(\omega_2)),-\tr(\rho(\omega_1\omega_2)))
\end{matrix}\right.
$$
(see \cite{Benedetto-Goldman:1999}). Using that
$$
\tr([M_1,M_2])=\tr(M_1)^2+\tr(M_2)^2+\tr(M_1M_2)^2-\tr(M_1)\tr(M_2)\tr(M_1M_2)-2,
$$
for all $M_1,M_2\in\SL(2,\C)$,
we note that those representations with given trace $\tilde d=\tr(\rho([\omega_1,\omega_2]))$ are parametrized by the affine cubic
$$
X^2+Y^2+Z^2+XYZ=\tilde d+2
$$
which is precisely $S_{(0,0,0,D)}$ with $D=\tilde d+2$. The reason is given
by the two-fold ramified cover $\pi:\T\to\Sphere$ used in \S \ref{par:MappingClassGroup}.
Consider a representation $\rho\in\Rep(\Sphere_4)$
corresponding to some point $(x,y,z)\in S_{(0,0,0,D)}$,
with local traces given by $(a,b,c,d)=(0,0,0,d)$, $D=4-d^2$.
One can lift the representation on the $4$-punctured torus,
where punctures are given by the set $H$ of $2$-torsion points.
Since $a=b=c=0$, we have 
$$
\rho(\alpha),\ \rho(\beta),\ \rho(\gamma)\sim\begin{pmatrix}i&0\\0&-i\end{pmatrix}
$$
and the lifted representation $\rho\circ\pi$ has local monodromy 
$-I$ around the corresponding punctures $(1/2,0)$, $(0,1/2)$ and $(1/2,1/2)$. After twisting $\rho\circ\pi$ by 
$-I$ at each of the punctures, we finally deduce a representation 
$\widetilde\rho\in\Rep(\T_1)$. Since $\pi_*\omega_1=\beta\gamma$ and $\pi_*\omega_2=\beta^{-1}\alpha^{-1}$ 
(see \S \ref{par:MappingClassGroup}), 
the character associated to the lifted representation $\widetilde{\rho}$ is given by
$$
\left\{\begin{matrix}
X&=&\tr(\widetilde{\rho}(\omega_1))&=&y\\
Y&=&\tr(\widetilde{\rho}(\omega_2))&=&x\\
Z&=&-\tr(\widetilde{\rho}(\omega_1\omega_2)&=&-z-xy
\end{matrix}\right.
$$
which satisfies $X^2+Y^2+Z^2+XYZ=4-d^2$. Note that
the local monodromy of $\widetilde{\rho}$ at $(0,0)$ is 
$-\rho(\delta^2)$ and we indeed find $\tilde d=2-d^2$.
We can now reverse the formulae and deduce that 
any representation $(X,Y,Z)\in\chi(\T_1)$ is the lifting
of a representation $(x,y,z)\in\chi(\Sphere_4)$.
This is due to the hyperelliptic nature of the once punctured torus.

%%%%%%%%%%%%%%%%%%%%%%%%%%%%%%%%%%%%%%%%%%

\section{Geometry and Automorphisms}\label{par:automorphisms-surfaces}

%%%%%%%%%%%%%%%%%%%%%%%%%%%%%%%%%%%%%%%%%%

This section is devoted to a geometric study of the family of surfaces $S_{(A,B,C,D)},$
and to the description of the groups of polynomial automorphisms $\Aut[S_{(A,B,C,D)}].$

In section  \ref{par:markov}, we describe a special case that is famous in Teichm\"uller 
theory.  Section \ref{par:classification-eph} introduces the concept of elliptic, parabolic,
and hyperbolic automorphisms of $S_{(A,B,C,D)}$.

%%%%%%%%%%%%%%%%%%%%%%%%%%%%%%%%%%%%%%%%%%

\subsection{The triangle at infinity and automorphisms}\label{par:triangle_automorphisms_area}

%%%%%%%%%%%%%%%%%%%%%%%%%%%%%%%%%%%%%%%%%%

 Let $S$ be any member of the
family $\SS$. The closure $\overline{S}$ of $S$ in $\P^3(\C)$ is given by a cubic homogeneous equation 
$$
w(x^2+y^2+z^2)+ xyz = w^2(Ax+By+Cz) + Dw^3.
$$
The intersection of $\overline{S}$ with the plane at infinity does not depend 
on the parameters and coincides with the triangle $\Delta$ given by the equation 
$$
\Delta\ \ \ :\ \ \ xyz=0;
$$
moreover, one easily checks that the surface $\overline{S}$ is smooth in a neighborhood of $\Delta$
(all the singularities of $\overline{S}$ are contained in $S$).

Since the equation defining $S$ is of degree $2$ with respect to the $x$ variable, 
each point $(x,y,z)$ of $S$ gives rise to a unique second point $(x',y,z)$. This procedure
determines a holomorphic involution of $S$, namely
$$
s_x(x,y,z)=(A-x-yz,y,z).
$$
This automorphism coincides with the automorphism of $S$ determined 
by the involution $s_x$ of $\Gamma_2^*$ (see equation \ref{eq:formula_involution_x}, \S \ref{par:modularformulae}).
Geometrically, the involution $s_x$ corresponds to the following: If $m$ is a point of
$\overline{S}$, the projective line which joins $m$ and the vertex $v_x=[1;0;0;0]$ of the triangle
$\Delta$ intersects ${\overline{S}}$ on a third point; this point is $s_x(m)$. 
The same construction provides two more involutions 
$$s_y(x,y,z)=(x,B-y-xz,z)\ \ \ \text{and}\ \ \ s_z(x,y,z)=(x,y,C-z-xy),$$ 
and therefore a subgroup 
$$
\A = \langle s_x,s_y,s_z\rangle
$$
of the group $\Aut[S]$ of polynomial automorphisms of the surface $S$. 

From section \ref{par:modularformulae}, we deduce that {\sl{ for any member $S$ of the 
family $\SS,$  the group $\A$ coincides with the image of $\Gamma_2^*$ into $\Aut[S],$ 
which is obtained by the action of $\Gamma_2^*\subset \MCG^*(\Sphere_4)$ on $\chi(\Sphere_4)$}}
 (see \S \ref{par:automorphisms_modular}). 

\begin{thm}\label{thm:elhuti} Let $S=S_{(A,B,C,D)}$ be any  member of the family of surfaces $\SS.$ Then 
\begin{itemize} 
 \item there is no non-trivial relation between the three involutions $s_x$, $s_y$ and $s_z$, and  $\A$ is 
therefore isomorphic to the free product $(\Z/2\Z)\star (\Z/2\Z)\star (\Z/2\Z)$ ;
\item the index of $\A$ in $\Aut[S]$ is bounded by $24$ ;
\item $\A$ coincides with the image of $\Gamma_2^*$ in $\Aut[S]$.
\end{itemize} 
Moreover, for a generic choice of the parameters $(A,B,C,D)$, $\A$ coincides with $\Aut[S]$.
\end{thm}

This result is almost contained in \`El'-Huti's article \cite{El-Huti:1974} and is more precise
than Horowitz's main theorem (see \cite{Horowitz:1972}, \cite{Horowitz:1975}).

\begin{proof}
Since ${\overline{S}}$ is smooth in a neighborhood of the triangle at infinity and the three
involutions are the reflexions with respect to the vertices of that triangle, we can
apply the main theorems of  \`El'-Huti's article: 
\begin{itemize}
 \item $\A$ is isomorphic to the free product 
$$
(\Z/2\Z)\star (\Z/2\Z)\star (\Z/2\Z)= \langle s_x \rangle \star \langle s_y \rangle \star \langle s_z \rangle ;
$$
\item $\A$ is of finite index in $\Aut[S]$ ;
\item $\Aut[S]$ is generated by $\A$ and the group of projective transformations of $\P^3(\C)$ 
which preserve $\overline{S}$ and $\Delta$ (i.e., by affine transformations of $\C^3$ that
preserve $S$). 
\end{itemize}
We already know that $\A$ and the image of $\Gamma_2^*$ in $\Aut[S]$ coincide. We now
need to study the index of $\A$ in $\Aut[S]$. Let $f$ be an affine invertible transformation of $\C^3,$ 
that we decompose as the composition of a linear part $M$ and a translation of vector $T.$
Let $S$ be any member of $\SS$. If $f$ preserves
$S,$  then the equation of $S$ is multiplied by a non zero complex number when we apply
$f.$ Looking at the cubic terms, this means that $M$ is a diagonal matrix composed with a permutation
of the coordinates. Looking at the quadratic terms, this implies that $T$ is the nul vector, so that
$f=M$ is linear. Coming back to the equation of $S,$ we now see that $M$ is one of the  $24$ linear transformations 
of the type $\sigma\circ \epsilon$ where $\epsilon$ either is the identity or changes the sign of two coordinates,
and $\sigma$ permutes the coordinates.
If $(A,B,C,D)$ are generic, $S_{(A,B,C,D)}$ is not invariant by any of these linear maps.
Moreover, one easily verifies that the subgroup $\A$ is a normal subgroup of $\Aut[S]$: If such a 
linear transformation $M=\sigma\circ \epsilon$ preserves $S,$ then it normalizes $\A.$ This shows
that $\A$ is a normal subgroup of $\Aut[S],$ the index of which is bounded by $24.$
\end{proof}

%%%%%%%%%%%%%%%%%%%%%%%%%%%%%%%%%%%%%%%%%%

\subsection{Consequences and notations}\label{par:consequences_notations}
As a corollary of theorem \ref{thm:elhuti} and proposition \ref{pro:mcg}, we get the following result:
{\sl{The mapping class group $\MCG_0^*(\Sphere_4)$ acts on the character variety $\chi(\Sphere_4),$ preserving each
surface $S_{(A,B,C,D)},$ and its  image in $\Aut[S_{(A,B,C,D)}]$ coincides with the image of $\Gamma_2^*,$ and therefore
with the finite index subgroup $\A$ of  $\Aut[S_{(A,B,C,D)}].$}} In other words, up to finite index subgroups, describing 
the dynamics of $\MCG^*(\Sphere_4)$ on the character variety $\chi(\Sphere_4)$ or of the group $\Aut[S]$ on $S$ for any
member $S$ of the family $\SS$ is one and the same problem.

\vv

\begin{figure}[t]\label{fig:markov}
\input{conj_c.pstex_t}
\caption{ {\sf{Conjugation for the Markov example.}} {{The right hand part of this figure 
depicts the dynamics of $\Gamma_2^*$ on $\M_+(\R)$, but viewed in $\P^2(\R)$ after the birational 
change of variables $[x:y:z:w]=[XQ:YQ:ZQ:XYZ]$, with $Q=X^2+Y^2+Z^2$. 
This change of variables sends the interior of the triangle $\{X \geq 0, Y\geq 0, Z\geq 0\}$ 
onto $\M_+(\R)$. }}}
\end{figure}
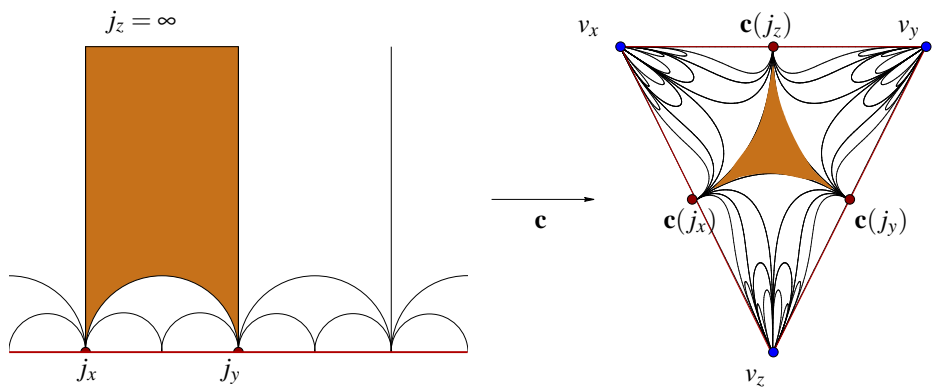

Let $\H$ be the Poincar\'e half plane. The group of isometries of $\H$ is isomorphic to 
$\PGL(2,\R)$: If $M$ is an element of $\GL(2,\R)$, its action on $\H$ is defined by
$$
M(z)=\frac{m_{11} z + m_{12}}{m_{21}z+ m_{22}} 
$$
if the determinant of $M$ is positive, and by the same formula but with $z$ replaced by ${\overline{z}}$
if the determinant is negative. In particular, $\Gamma_2^*$ acts isometrically on $\H$. 
Let $j_x$, $j_y$ and $j_z$ be the three points on the boundary 
of $\H$ with coordinates $-1$, $0$, and $\infty$ respectively. Let $r_x$ (resp. $r_y$ 
and $r_z$) be the reflection of  $\H$ around the geodesic between $j_y$ and $j_z$ 
(resp. $j_z$ and $j_x$, resp. $j_x$,   and $j_y$). These isometries are respectively induced by the three matrices $s_x,$ $s_y,$ and $s_z$ given in section \ref{par:modularformulae}. 
As a consequence, $\Gamma_2^*$ coincides with the group of  symmetries of the tesselation 
of $\H$ by ideal triangles, one of which has vertices $j_x$, $j_y$ and $j_z$ (see the left part 
of figure \ref{fig:markov}). 

\vv

In the following, we shall identify the subgroup 
$\Gamma_2^*$ of $\PGL(2,\Z)$ and the subgroup $\A$ of 
$\Aut[S_{(A,B,C,D)}]$ : If $f$ is an element of $\A$, $M_f$ will denote
the associated element of $\Gamma_2^*$ (either viewed as a matrix or an isometry of $\H$),
and if $M$ is an element of $\Gamma_2^*$, $f_M$ will denote the automorphism associated 
to $M$ (for any surface $S$ of the family $\SS$). If $f$ is one of 
the three involutions $s_x,$ $s_y,$ or $s_z$ (resp. the three elements $g_x,$ $g_y,$ or $g_z$), we shall use exactly the
same  letters to denote the element $f$ of $\Gamma_2^*$ or the corresponding automorphism $f\in \A.$  The only place where this rule is not followed is when we study
the action of $\Gamma_2^*$ on the PoincarŽ disk: We then use the notation $r_x,$ $r_y,$ and $r_z$ to denote the involutive isometrys induced by $s_x,$ $s_y,$ and $s_z.$

%%%%%%%%%%%%%%%%%%%%%%%%%%%%%%%%%%%%%%%%%%

\subsection{Elliptic, Parabolic, Hyperbolic}\label{par:classification-eph}
Non trivial isometries of $\H$ are classified into three different species. 
Let $M$ be an element of $\PGL(2,\R)\setminus\{{\text{Id}}\}$, viewed as  an isometry of $\H$. Then,

\begin{itemize}
 \item $M$ is elliptic if $M$ has a fixed point in the interior of $\H$. 
Ellipticity is equivalent to $\det(M)=1$ and $\vert\tr(M)\vert < 2$ (in which case
$M$ is a rotation around a unique fixed point)  or $\det(M)=-1$ and $\tr(M)=0$
(in which case $M$ is a reflexion around a geodesic of fixed points). 
\item $M$ is parabolic if $M$ has  a unique fixed point, which is located
on the boundary of $\H$; $M$ is parabolic if and only if  $\det(M)=1$ and $\tr(M)=2$ or $-2$;  
\item $M$ is hyperbolic if it has exactly two fixed points which are on the boundary of $\H$; this
occurs if and only if $\det(M)=1$ and $\vert \tr(M) \vert >2$,
or $\det(M)=-1$ and $\tr(M)\neq 0$.
\end{itemize}

An element $f$  of $\A\setminus \{{\text{Id}}\}$ will be termed elliptic, parabolic, or hyperbolic, according 
to the type of $M_f$. Examples of elliptic elements are given by the three involutions
$s_x$, $s_y$ and $s_z$. Examples of parabolic elements are given by 
the three automorphisms $g_x,$ $g_y$ and $g_z$ (see  section \ref{par:modularformulae}).
The dynamics of these automorphisms will be described in details in \S \ref{par:dynamics_parabolic}.
Let us just mention the fact that $g_x$ (resp. $g_y$, $g_z$) preserves the 
conic fibration $\{x=c^{ste}\}$ (resp. $\{y=c^{ste}\}$, $\{z=c^{ste}\}$) of any member $S$ of $\SS$. 

\begin{pro}\label{pro:classification_isometry}
Let $S$ be one of the surfaces in the family $\SS$ ($S$ may be singular). 
An element $f$ of $\A$ is
\begin{itemize}
  \item  elliptic if and only if  $f$ is conjugate to one of the 
involutions $s_x$, $s_y$ or $s_z$, if and only if $f$ is periodic; 
\item parabolic if and only if $f$ is conjugate to a non trivial power of one of the 
automorphisms $g_x$, $g_y$ or $g_z$;
\item hyperbolic if and only if $f$ is conjugate to a cyclically reduced composition
which involves the three involutions $s_x$, $s_y$, and $s_z$. 
\end{itemize}
\end{pro}

\begin{proof}
Since $\Gamma_2^*$ and the image of $\A$ in $\Aut[S]$ are isomorphic for any
$S$ in $\SS$, we just need to prove the same statement for $\Gamma_2^*$.
The group $\Gamma_2^*$ is a subgroup of $\PGL(2,\Z)$. As a consequence, any 
elliptic element of $\Gamma_2^*$ is periodic. Since 
$$
\Gamma_2^*= (\Z/2\Z)\star (\Z/2\Z)\star (\Z/2\Z),
$$
any periodic element of $\Gamma_2^*$ is conjugate to one of the involutions
$r_x$, $r_y$, $r_z$ (see for example \cite{Serre:book}), and the
first property is proved.

\vv

If $M$ is a parabolic element of $\Gamma_2^*$, its unique fixed point on the
boundary $\R\cup\{\infty\}$ of $\H$ is a rational number. The action of $\Gamma_2^*$
on the set $\Q\cup\{\infty\}$ of rational numbers has three distinct orbits: The orbits
of $j_x=-1$, $j_y=0$ and $j_z=\infty$. This implies that there exists an element $F$
of $\Gamma_2^*$ such the $FMF^{-1}$ is parabolic and fixes one of these three points, 
say $j_z$. Any parabolic element $G$ of $\Gamma_2^*$ that fixes $\infty$ is of the type
$$
\pm \left( \begin{array}{cc}
 1 & 2 k \\ 0 & 1
\end{array}
\right)
$$ 
where $k$ is an integer. This fact shows that $M$ is conjugate to a power
of $g_z$ (see  section \ref{par:modularformulae}) and concludes the proof of the second point.

\vv

Let $M$ be a hyperbolic element of $\Gamma_2^*$. After conjugation, we can write
$M$ as a cyclically reduced word in the involutive generators $r_x$, $r_y$ and 
$r_z$. If the number of involutions that appear in this composition is equal to 
$1$ or $2$, then $M$ is an involution or a power of $g_x$, $g_y$ or $g_z$. The third
property follows from this remark.
\end{proof}

\begin{rem}\label{rem:cyclic}
The three vertices $j_x$, $j_y$ and $j_z$ disconnect $\partial \H$ in three
segments $[j_y,j_z]$, $[j_z,j_x]$ and $[j_x,j_y]$.
Let $M$ be a hyperbolic element of $\Gamma_2^*$. Let $\alpha_M$ be the repulsive
and attrating fixed points of $M$ on the boundary of $\H$.
The Fricke-Klein ping-pong lemma, as described in \cite{delaHarpe:book}, page 25,
shows that $M$ is a cyclically reduced composition of $r_x$, $r_y$, and $r_z$ if and
only if the fixed points of $M$ are contained in two distinct connected components of
$\partial \H\setminus \{ j_x, j_y, j_z\}$. 
\end{rem}

%%%%%%%%%%%%%%%%%%%%%%%%%%%%%%%%%%%%%%%%%%%%%%%%%%%%%%%%%%%%%%

\subsection{The Markov surface}\label{par:markov} 

Let $\M$ be the element of $\SS$ corresponding to the parameter $(A,B,C,D)=(0,0,0,0)$.
After a simultaneous multiplication of each coordinate by $-3$, the equation of $\M$ 
is 
$$
x^2+y^2+z^2=3xyz.
$$
This surface has been studied by Markov in 1880 in his  papers
concerning diophantine approximation.  The real part $\M(\R)$ of the Markov
surface has an isolated singular point at the origin and four other connected components, each of which 
is homeomorphic to a disk. One of these components is 
$$
\M_+(\R)=\M(\R)\cap (\R^+)^3.
$$

\begin{pro}[Markov, \cite{Cassels:book}]
 The action of $\A=\Gamma_2^*$ on the Markov surface $\M$ preserves each connected
component of $\M(\R)$. There exists a diffeomorphism ${\mathbf{c}}:\H\to \M_+(\R)$
such that $(i)$ the image of the (closed) ideal triangle with vertices $j_x$, $j_y$ and $j_z$ is the subset
of $\M_+(\R)$ defined by the three inequalities 
$$
xy\leq 2z, \quad  yz\leq 2 x, \quad {\text{and}} \quad zx \leq 2y,
$$ 
and (ii) ${\mathbf{c}}$ conjugates
the action of $\Gamma_2^*$ on $\H$ with the action of $\Gamma_2^*$ on $\M_+(\R)$ in such a way that 
$$
{\mathbf{c}}\circ r_x = s_x\circ {\mathbf{c}},  \quad {\mathbf{c}}\circ r_y = s_y\circ {\mathbf{c}}, \quad {\text{and}} \quad{\mathbf{c}}\circ r_z = s_z\circ {\mathbf{c}}.
$$
\end{pro}

\begin{rem} We refer the reader to \cite{Cassels:book} or 
\cite{Goldman:2003} for a proof (see figure \ref{fig:markov} for a visual argument).
This result is not surprising if one notices that $\M_+(\R)$ is a model of 
the Teichm\"uller space of the once 
punctured torus  with a cusp at the puncture, and finite 
area $2\pi$. \end{rem}

%\begin{rem} Notice that this proposition implies that the dynamics of $\Gamma_2^*$ on $\M(\R)$ is totally
%discontinuous. In particular, if $f$ is an element of $\Gamma_2^*$, $f:\M(\R)\to \M(\R)$
%has a unique peridic point, the origin. \end{rem}

%%%%%%%%%%%%%%%%%%%%%%%%%%%%%%%%%%%%%%%%%%%%%%%%%%

\subsection{An (almost) invariant area form}

%%%%%%%%%%%%%%%%%%%%%%%%%%%%%%%%%%%%%%%%%%%%%%%%%%

The monomial action of the group $\GL(2,\Z)$ on $\C^*\times \C^*$ almost preserves the holomorphic 
$2$-form 
$$
\Omega = \frac{dx}{x}\wedge \frac{dy}{y}.
$$
More precisely, $M^*\Omega=\pm  \Omega$ for any element $M$ of $\GL(2,\Z)$. This form 
is invariant under the action of $\eta$ and determines a holomorphic volume form on 
the Cayley cubic, that is almost $\Aut[S_C]$ invariant.This property
is shared by all the members of $\SS$ (the proof is straightforward).

\begin{pro}\label{pro:AreaForm}
Let $S\in \SS$ be the surface corresponding to the parameters $(A,B,C,D)$. 
The volume form $\Omega$, which is globally defined by the formulas 
$$
\Omega = \frac{dx\wedge dy}{2z+xy-C} = \frac{dy\wedge dz}{2x + yz -A} = \frac{dz\wedge dx}{2y+zx-B}
$$
on $S\setminus\Sing(S)$, is almost invariant under the action of $\Aut[S]$, by which we mean 
that $f^*\Omega = \pm \Omega$  for any $f$ in $\Aut[S]$.
\end{pro}

%%%%%%%%%%%%%%%%%%%%%%%%%%%%%%%%%%%%%%%%%%%%%%%%%%%%%%%%%%%%%%%

\subsection{Singularities, fixed points, and an orbifold structure}\label{par:sing_fp_orb}

The singularities of the elements of $\SS$ will play an important role in this article. 
In this section, we collect a few results regarding these singularities.

\begin{lem} Let $S$ be a member of $\SS$. 
A point $m$ of $S$ is singular if and only if $m$ is a fixed point of the group $\A$. 
\end{lem}

\begin{proof} This is a direct consequence of the fact that $m$ is a fixed point
 of $s_x$ if and only if   $2x+yz=Ax$, if and only if  the partial derivative of the equation of $S$ with respect to the $x$-variable vanishes.   \end{proof}

\begin{eg}\label{eg:singu}
The family of surfaces with parameters 
$(4+2d,4+2d,4+2d, -(8+ 8d + d^2))$ with $d\in \C$
is a deformation of the Cayley cubic, that corresponds to $d=-2$, and any of these
surfaces has $3$ singular points (counted with multiplicity). 
\end{eg}

\begin{lem}
If $m$ is a singular point of $S$, there exists a neighborhood of $m$ which is isomorphic
to the quotient of the unit ball in $\C^2$ by a finite subgroup of $\SU(2)$.  
\end{lem}

\begin{proof} 
Any singularity of a cubic surface is a quotient singularity, except when the singularity
is isomorphic to $x^3+y^3+z^3+\lambda xyz=0$,
for at least one parameter $\lambda$ (see \cite{Bruce-Wall:1979}).
Since the second jet of the equation of $S$ never vanishes when $S$ is a member of $\SS$,
the singularities of $S$ are quotient singularities. Since $S$ admits a global
volume form $\Omega$, the finite group is conjugate to a subgroup of $\SU(2,\C)$. 
\end{proof}

As a consequence, any member $S$ of $\SS$ is endowed with a well defined orbifold 
structure. If $S$ is singular, the group $\A$ fixes each of the singular
points and preserves the orbifold structure. We shall consider this action in 
the orbifold category, but we could as well  extend the action of $A$ to a smooth desingularization of $S.$

\begin{lem}\label{Lem:TopologyCubic}
The complex affine surface $S$ is simply connected. 
When $S$ is singular, the fundamental group of the complex surface $S\setminus {\sf{Sing}}(S)$ 
is normally generated by the local finite fundamental groups around the singular points.
\end{lem}

\begin{proof} First of all, recall that a smooth cubic surface in $\P^3(\C)$
may be viewed as the blowing-up of $\P^2(\C)$ at $6$-points
in general position. Let us be concrete.
After a projective change of coordinates, one can assume that
those $6$ points lie on the triangle $XYZ=0$ 
and are labelled as follows
$$
p_i=[0:1:u_i],\ \ \ q_i=[v_i:0:1]\ \ \ \text{et}\ \ \ r_i=[1:w_i:0],\ \ \ i=1,2
$$
where $[X:Y:Z]$ are projective coordinates of $\P^2.$ 
One can moreover assume that the three following products
take the same value $\lambda$:
$$u_1u_2=v_1v_2=w_1w_2=:\lambda.$$
Now, consider the map 
$$\Phi:\mathbb P^2\dashrightarrow \C^3\ ;\ (X:Y:Z)
\mapsto \left(\frac{P}{YZ},\frac{Q}{XZ},\frac{R}{XY}\right)$$
where $P$, $Q$ and $R$ are degree $2$ homogeneous polynomials given by
$$\left\{\begin{matrix}
P&=&-X^2-\frac{1}{\lambda}Y^2-\lambda Z^2+(\frac{1}{w_1}+\frac{1}{w_2})XY+(v_1+v_2)XZ\\ 
Q&=&-\lambda X^2-Y^2-\frac{1}{\lambda}Z^2+(w_1+w_2)XY+(\frac{1}{u_1}+\frac{1}{u_2})YZ\\ 
R&=&-\frac{1}{\lambda}X^2-\lambda Y^2-Z^2+(\frac{1}{v_1}+\frac{1}{v_2})XZ+(u_1+u_2)YZ
\end{matrix}\right.$$
For $u_i$, $v_i$ and $w_i$ generic, 
the map $\Phi$ sends the triangle $XYZ=0$ to the triangle at infinity
$xyz=0$ of $\P^3\supset\C^3$ and has simple indeterminacy points 
exactly at $p_i,$ $q_i,$ and $r_i,$ $i=1,2.$ Let $\tilde S$ be the surface obtained by blowing-up the $6$ 
indeterminacy points of $\Phi.$ One can check that the image
of $\Phi:{\tilde{S}} \to \P^3(\C)$ is exactly the cubic surface 
$S=S_{(A,B,C,D)},$ with parameters
$$\left\{\begin{matrix}
A&=&\left(\frac{v_1}{w_1}+\frac{v_2}{w_2}+\frac{v_1}{w_2}+\frac{v_2}{w_1}\right)
-\left(u_1\lambda+\frac{1}{u_1\lambda}+u_2\lambda+\frac{1}{u_2\lambda}\right)\\ 
B&=&\left(\frac{u_1}{w_1}+\frac{u_2}{w_2}+\frac{u_1}{w_2}+\frac{u_2}{w_1}\right)
-\left(v_1\lambda+\frac{1}{v_1\lambda}+v_2\lambda+\frac{1}{v_2\lambda}\right)\\ 
C&=&\left(\frac{u_1}{v_1}+\frac{u_2}{v_2}+\frac{u_1}{v_2}+\frac{u_2}{v_1}\right)
-\left(w_1\lambda+\frac{1}{w_1\lambda}+w_2\lambda+\frac{1}{w_2\lambda}\right)\\
D&=&\sum_{i,j,k\in\{1,2\}}\left(u_iv_jw_k+\frac{1}{u_iv_jw_k}\right)\\
&&-\left(\frac{u_1}{u_2}+\frac{u_2}{u_1}+\frac{v_1}{v_2}+\frac{v_2}{v_1}+\frac{w_1}{w_2}+\frac{w_2}{w_1}+\lambda^3+\frac{1}{\lambda^3}+4\right)
\end{matrix}\right.$$
Singular cubics arise when $3$ of the $6$ points lie on 
a line, or all of them lie on a conic. In this case, the corresponding
line(s) and/or conic have negative self-intersection in 
$\tilde S,$  
and are blown-down by $\Phi$ to singular point(s) of  
$S.$ A smooth resolution of $S$ is therefore given by $\tilde S$.

Our claim is that the quasi-projective surface $\tilde S'$ 
obtained by deleting the strict transform of
the triangle $XYZ=0$ from ${\tilde{S}}$ is simply connected. Indeed, the 
fundamental group of $\P^2-\{XYZ=0\}$ is isomorphic to 
$\Z^2$, generated by two loops, say one turning around $X=0$,
and the other one around $Y=0$. After blowing-up one 
point lying on $X=0$, and adding the exceptional divisor 
(minus $X=0$), the first loop becomes homotopic to $0$; 
after blowing-up the $6$ points and adding all exceptional divisors, the two generators 
become trivial and the resulting surface $\tilde S'$ is simply connected.
The affine surface $S$ is obtained after blowing-down 
some rational curves in $\tilde S$ and is therefore simply connected
as well. 

The second assertion of the lemma directly follows from Van Kampen
Theorem.
\end{proof}

\section{Birational Extension and Dynamics}\label{par:bira-extension-dyna}
%%%%%%%%%%%%%%%%%%%%%%%%%%%%%%%%%%%%%

\subsection{Birational transformations of surfaces}\label{par:bira_nota} Let $f$ be a birational 
transformation of a complex projective surface $X$ and $\Ind(f)$ be its 
indeterminacy set. The critical set of $f$ is the union of all the curves
$C$ in $S$ such that $f(C\setminus \Ind(f))$ is a point (in fact a point
of $\Ind(f^{-1})$). 
One says that $f$ is not algebraically stable if there is a curve $C$ 
in the critical set and a positive integer $k$ such that $f^k(C\setminus \Ind(f))$
is contained in $\Ind(f).$ Otherwise, $f$ is said to be algebraically stable
(see \cite{Diller-Favre:2001}).
Let $H^2(X,\Z)$ be the second cohomology group of $X$ and $f^*:H^2(X,\Z)\to H^2(X,\Z)$
be the linear transformation induced by $f.$ It turns out that $f$ is algebraically stable
if and only if $(f^k)^*=(f^*)^k$ for any positive integer $k$ (see \cite{Diller-Favre:2001}).
More generally, $(f\circ g)^*=g^*\circ f^*$ if and only if $g$ does not blow down any curve
onto an element of $\Ind(f).$ 

\vv

The (first) dynamical degree $\lambda(f)$ of $f$ is the spectral radius of the sequence
of linear operators $(f^k)^*.$ If $f$ is algebraically stable, $\lambda(f)$ is therefore
the largest eigenvalue of $f^*.$ It follows from  Hodge theory that
$$
\limsup \frac{1}{k}  \log \Vert (f^k)^* [v]\Vert = \log \lambda(f).
$$
for any class $[v]$ that is obtained through a hyperplane section of $X.$
The dynamical degree of $f$ is invariant under birational conjugation   (see \cite{Diller-Favre:2001,Guedj:etats}),
and provides an upper bound for the 
topological entropy of $f$ (see \cite{Diller-Favre:2001,Guedj:etats}).

\begin{eg}\label{eg:monomial_degree}
If $M$ is an element of $\GL(2,\Z),$  $M$ acts on
$\C^*\times \C^*$ monomially (see equation \ref{eq:monomial}). 
The dynamical degree of this monomial transformation is equal to the spectral radius
$\rho(M)$ of $M.$ If $f_M$ is the automorphism of the Cayley cubic 
$S_C$ which is induced by $M,$ the dynamical degree of $f_M$ coincides also
with $\rho(M)$ 
(see \cite{Favre:2003}, or the survey article \cite{Guedj:etats}). 
\end{eg}

%%%%%%%%%%%%%%%%%%%%%%%%%%%%%%%%%%%

\subsection{Birational extension} Let $S$ be a member of the family $\SS.$ The group $\A$ acts by polynomial automorphisms on $S$ and also  by birational transformations
of the compactification $\overline{S}$ of $S$ in $\P^3(\C).$ 
Let $\Delta$ be the triangle at infinity, $\Delta={\overline{S}}\setminus S.$
The three sides of this triangle are the lines $D_x=\{x=0,w=0\},$
$D_y=\{y=0,w=0\}$ and $D_z=\{z=0,w=0\}$; the vertices are $v_x=[1:0:0:0],$
$v_y=[0:1:0:0]$ and $v_z=[0:0:1:0].$ The ``middle points'' of the sides are respectively
$$
m_x=[0:1:1:0], \quad m_y=[1:0:1:0], \,{\text{and}} \, \, m_z=[1:1:0:0]
$$
(see  figure \ref{fig:markov} in \S \ref{par:markov}).
Let $ V $ be the subspace of $H^2(\overline{S},\Z)$ defined by
$$
V= \Z[D_x]+ \Z [D_y] + \Z [D_z],
$$
where $[D_x]$ denotes either the homology class of $D_x$ in $H_2(\overline{S},\Z)$ or
its dual in $H^2(\overline{S},\Z).$  Since $\Delta$ is $\A$-invariant, 
the  action of any element $f$ in $\A$ on $H^2(\overline{S},\Z)$  preserves the subspace $V.$ 

\begin{lem}[see \cite{El-Huti:1974} or \cite{Iwasaki-Uehara:2006}]\label{lem:IU}
The involution $s_x$ acts on the triangle $\Delta$ in the following way.
\begin{itemize}
 \item The image of the side $D_x$ is the vertex $v_x$ and the vertex $v_x$ is 
blown up onto the side $D_x.$
\item the sides $D_y$ and $D_z$ are invariant and $s_x$ permutes the vertices
and fixes the middle point of each of these sides.
\end{itemize}
\end{lem}

Of course, we have the same result for $s_y$ and $s_z,$ with the obvious 
required modifications. In the following, we shall denote by ${\sf{s}}^*_x$ (resp. 
${\sf{s}}^*_y$ or ${\sf{s}}^*_z$) the restriction of $(s_x)^*$ (resp. $(s_y)^*$ 
or $(s_z)^*$) on the subspace $V$ of $H^2({\overline{S}}, \Z).$  

\begin{rem}\label{rem:indeterminacy} 
The ``action'' of $\A$ on the triangle $\Delta$ does not 
depend on the choice of the parameters $(A,B,C,D).$ 
Let $f=w(s_x,s_y,s_z)$ be an element of $\A,$ given by a reduced word in the letters
$s_x,$ $s_y$ and $s_z.$ Since $s_x$ (resp. $s_y,$ $s_z$) does not blow down any 
curve on indeterminacy points of the other two involutions, the linear transformation
$f^*:V\to V$ 
is the composition $f^*=w^t(s_x^*,s_y^*,s_z^*),$ where $w^t$ is the transpose of $w$
(see section \ref{par:bira_nota}).
If $w$ ends with $s_x$ (resp. $s_y$ or $s_z$), then $f$
contracts the side $D_x$ (resp. $D_y$ or $D_z$). 
If $w$ starts with $s_x$ (resp. $s_y$ or $s_z$), the image of the critical set of $f$ is the vertex $v_x$ (resp. $v_y$ or $v_z$).
In particular, $\Ind(f)$ and $\Ind(f^{-1})$ are not empty if $f$ is different from
the identity.
\end{rem}

\begin{eg}\label{eg:parabolic-bira} The element $g_x=s_z\circ s_y$ preserves 
 the coordinate variable $x.$ Its
action  on $\Delta$ is the following: $g_x$ contracts
both $D_y$ and $D_z\setminus\{v_y \}$ on $v_z, $ and preserves $D_x$;  its inverse 
contracts $D_y$ and $D_z\setminus\{v_z\}$ on $v_y.$ In particular $\Ind(g_x)=v_y$ and 
$\Ind(g_x^{-1})=v_z.$ The elements $g_y$ and $g_z$ act in a similar way. In particular,
$g_x,$ $g_y$ and $g_z$ are algebraically stable.
\end{eg}

Let us now present a  nice way of describing the ``action'' of 
$\A,$ {\sl{i.e.}} of $\Gamma_2^*,$ on the triangle $\Delta.$ 
Since this action does not depend on the parameters, we choose
 $(A,B,C,D)=(0,0,0,0)$ and use what we know about
the Markov surface $\M$ (see \S \ref{par:markov}). 
The closure of $\M_+(\R)$ in $\overline{\M}$ contains a part of the triangle at infinity, 
namely the set $\Delta_+(\R)$ of points $[x:y:z:0]$ such that
$xyz=0,$ and $x,\, y,\, z \geq 0.$ This provides a compactification of $\M_+(\R)$ by 
the triangle $\Delta_+(\R).$ The conjugation 
$$
{\mathbf{c}}:\H \to \M_+(\R)
$$
between the Poincar\'e half plane and $\M_+(\R)$ described in \S \ref{par:markov}
does not extend up to the boundary of this compactification.
Nevertheless, one can ``extend'' the map ${\mathbb{c}}$ 
in the following way (see  figure \ref{fig:markov}): 
\begin{itemize}
\item the three segments $(j_y,j_z),$ $(j_z,j_x)$ and $(j_x,j_y)$ of $\partial \H$ are sent
to the three vertices $v_x,$ $v_y$ and $v_z$ of $\Delta$;
 \item the three points $j_x,$ $j_y,$ and $j_z$ are ``sent'' to the three sides $D_x,$ $D_y$ and $D_z$ of $\Delta_+(\R)$
by ${\mathbf{c}}$ (or equivalently to the middle points $m_x,$ $m_y$ and $m_z$);
 \end{itemize}
Then, if $M$ is a hyperbolic element of $\Gamma_2^*,$ the two fixed points of 
$M$ on the boundary of $\H$ are sent to the indeterminacy points of $f_M$ 
and $f_M^{-1}$: If $M$ is hyperbolic, with one attractive 
fixed point $\omega_M$ and one repulsive fixed point $\alpha_M$, 
then 
$$
\Ind(f_M)={\mathbf{c}}(\alpha_M), \quad \Ind(f_M^{-1})={\mathbf{c}}(\omega_M).
$$

\begin{rem} Let us  consider the surface obtained by blowing up the vertices
of the triangle $\Delta.$ This is a new compactification of the affine cubic $S_M$ by a
cycle of six rational curves. Then we blow up the six vertices of this hexagon, 
and so on : This defines a sequence of rational surfaces $S^i$. Let $S^\infty$ 
be the projective limit of these surfaces. The group $\Gamma_2^*$ acts
continously on this space, and we can
extend $\mathbf{c}^{-1}$ so as to obtain a semi-conjugation 
between the action on $S_M ^\infty\setminus S_M$ and the action of $\Gamma_2^*$ on the 
circle. Such a construction is presented in details in a similar context in
 \cite{Hubbard-Papadopol:preprint}, chapter 4 (see also \cite{Cantat:cremona} for a related approach).
\end{rem}

The following proposition reformulates and makes more precise, section 7 of \cite{Iwasaki-Uehara:2006}.

\begin{pro}\label{pro:IU}
Let $S$ be any member of the family $\SS$ and $f$ an element of $\A.$ 
\begin{itemize}
\item The birational transformation $f:{\overline{S}}\to {\overline{S}}$
is algebraically stable if, and only if  $f$  is a cyclically reduced composition
of the three involutions $s_x,$ $s_y$ and $s_z$ of length at least $2.$ 

\item Every hyperbolic element $f$ of $\A$ is conjugate to an algebraically stable 
element of $\A.$

 \item If $f$ is algebraically stable and hyperbolic, $\Ind(f)$ and $\Ind(f^{-1})$ are two distinct vertices of
$\Delta,$ and $f^n$ contracts the whole triangle $\Delta\setminus \Ind(f)$ onto $\Ind(f^{-1})$
as soon as $n$ is a  positive integer.
\end{itemize}
\end{pro}

\begin{proof}  If $\Ind(f)=\Ind(f^{-1})\neq \emptyset,$
$f$ is not algebraically stable. This shows, for example, that an involution with a
non empty indeterminacy set is not algebraically stable. 

Let $M$ be an element of $\Gamma_2^*\setminus\{Id\}$ and $f_M$ the corresponding element 
of $\A,$ viewed as a birational transformation of $S.$
From remark \ref{rem:indeterminacy},
we know that $\Ind(f_M)$ is  non empty, and from proposition \ref{pro:classification_isometry}
that any elliptic element of $\Gamma_2^*$ is an involution. This shows that $f_M$
is not algebraically stable if $M$ is elliptic. 

Let us now fix a non elliptic element $M$ of $\Gamma_2^*,$ 
which we write as a reduced word $w(r_x,r_y,r_z)$ in 
the generators $r_x,$ $r_y$ and $r_z$ of $\Gamma_2^*$
(see \S \ref{par:classification-eph}).

Let us first assume that $M$ is parabolic. If $f_M$ is a non trivial iterate  
of $g_x$ (resp. $g_y$
or $g_z$), we know from example \ref{eg:parabolic-bira} that $M$ is algebraically 
stable. 
If not, the unique fixed point of $M$ on $\partial H$ is different from 
$j_x,$ $j_y$ and $j_z$ 
and its image by ${\mathbf{c}}$ is a vertex of $\Delta.$ This vertex $v$ coincides with 
$\Ind(f_M)$ and $\Ind(f_M^{-1}),$ and $f_M$ is not algebraically stable. Since
$M$ is cyclically reduced if, and only if $M$ is an iterate of $g_x,$ $g_y,$ or $g_z,$ 
the result is proved in the parabolic case.  

Let us now suppose that $M$ is hyperbolic :
The fixed points $\alpha_M$ and $\omega_M$ define
two distinct elements of $\partial \H\setminus \{j_x,j_y,j_z\}$ 
and the indeterminacy sets of 
$f_M$ and $f_M^{-1}$ are  the vertices
$\Ind(f_M)={\mathbf{c}}(\alpha_M)$ and $\Ind(f_M^{-1})={\mathbf{c}}(\omega_M)$ of $\Delta.$
These vertices  are distinct if, and only if 
$\alpha_M$ and $\omega_M$ are contained in two distinct components of 
$\partial \H\setminus \{j_x,j_y,j_z\},$ if, and only if $f_M$ 
is a cyclically reduced composition of the three involutions $s_x,$ $s_y,$ $s_z$ 
(see remark \ref{rem:cyclic}). This shows that $f_M$ is not algebraically stable if $w$ is
not cyclically reduced. In the other direction, if $w$ is cyclically reduced, then ${\mathbf{c}}(\omega_M)$ 
is not an indeterminacy point of $f_M,$ $f_M$ fixes this point, and contracts the three sides of $\Delta$
on this vertex. As a consequence, the positive orbit of $\Ind(f_M^{-1})$ does 
not intersect
$\Ind(f_M),$ and $f_M$ is algebraically stable. 
\end{proof}

\begin{thm}
Let $f$ be an element of $\A$ and $M_f$ the element of $\PGL(2,\Z)$ which 
is associated to $f.$ The dynamical degree $\lambda(f)$ is equal to the 
spectral radius of $M_f.$
\end{thm}

This result is different from, but similar to, the main theorem of \cite{Iwasaki-Uehara:2006},
which provides another algorithm to compute  $\lambda(f).$ 

\begin{proof}
Let $f$ be an element of $\A.$ After conjugation inside $\A$
(this does not change the dynamical degree and the spectral radius of $M_f$), 
we can assume that $f=w(s_x,s_y,s_z)$ is a cyclically reduced word. If $f$ is periodic, then $f$ 
is one of the involutive generators and the theorem is proved. If $f$ is parabolic, 
then $f$ is conjugate to an iterate of $g_x,$ $g_y$ or $g_z,$ $f$ preserves a fibration
of $S$ into rational curves, and $\lambda(f)=1.$
If $f$ is hyperbolic, proposition \ref{pro:IU} shows that $f$ is algebraically stable.
Let $[v]=[D_x]+[D_y]+[D_z]$ be the class of the hyperplane section of $\overline{S}$
which is obtained by cutting $\overline{S}$ with the plane at infinity. We know that 
$$
\limsup_{k\to \infty} \left(\frac{1}{k} \log\Vert (f^k)^*[v] \Vert\right) = \log (\lambda(f)).
$$ 
Since the action of $f^*$ on the subspace $V$ of $H^2(X,\Z)$  does not depend
on the parameters $(A,B,C,D),$ and since $[v]$ is contained in $V,$ $\lambda(f)$ 
does not depend on $(A,B,C,D).$ 
 Consequently, to calculate $\lambda(f),$ we can choose the parameters
$(0,0,0,4)$ and work on the Cayley cubic. The conclusion now follows from example \ref{eg:monomial_degree}.
\end{proof}

%%%%%%%%%%%%%%%%%%%%%%%%%%%%%%%%%%%%%%%%%%%%%%

\subsection{Entropy of birational transformations}

Let $f$ be a hyperbolic element of  $\A$ (see section \ref{par:classification-eph}).  
Up to conjugation, the birational transformation $f:{\overline{S}}\to {\overline{S}}$ 
is algebraically stable, $\Ind(f)$ is a fixed point of $f^{-1}$  and $\Ind(f^{-1})$ is
a fixed point of $f.$ As remarked in \cite{Iwasaki-Uehara:2006}, this enables
us to apply the main results from \cite{Bedford-Diller:2005} and \cite{Dujardin:2006}. 

\begin{thm}(Bedford, Diller, Dujardin, Iwasaki, Uehara)
Let $f$ be an element of the group $\A$ and $S$ be an element of $\SS.$ The topological
 entropy of $f_M:S\to S$ is equal to the logarithm of the spectral radius $\lambda(f)$ of $M_f,$
the number of periodic (saddle) points of $f$ of period $n$ grows like $\lambda(f)^n$ and these
points equidistribute toward an ergodic measure of maximal entropy for $f.$
\end{thm}

In \cite{Cantat:BHPS}, we shall explain how the dynamics of $f$ is related to the
dynamics of H\'enon mappings, and deduce a much more precise description
of the dynamics. 

\begin{eg}
Let $M$ be an element of $\GL(2,\Z).$ Let $U$ be the unit circle in $\C^*$ and $\T$  be
the subgroup $U\times U$ of $\C^*\times \C^*.$ The monomial automorphism $M$ of $\C^*\times \C^*$ 
preserves $T$ and induces a ``linear'' automorphism on this real torus. The entropy of 
$M:\T\to \T$ is equal to the logarithm of the spectral radius of $M.$ If $(x,y)$ is 
a point of $\C^*\times \C^*,$ the orbit $M^n(x,y),$ $n\geq 0,$ converges toward $\T$ or goes to infinity.
The same property remains true for the dynamics of $f_M$ on the Cayley cubic $S_C$; the role played by $\T$ is now 
played by $\T/\eta= S_C(\R)\cap [-2,2]^3.$
\end{eg}

%%%%%%%%%%%%%%%%%%%%%%%%%%%%%%%%%%%%%%%%%%%%%%%
\section{Bounded Orbits}\label{par:sectionbounded}

%%%%%%%%%%%%%%%%%%%%%%%%%%%%%%%%%%%%%%%%%%%%%%%

\subsection{Dynamics of parabolic elements}\label{par:dynamics_parabolic}

%%%%%%%%%%%%%%%%%%%%%%%%%%
Parabolic elements will play an important role in the proof of theorem \ref{thm:main_affine}. 
In this section, we describe the dynamics of these automorphisms, on any member
$S$ of our family of cubic surfaces. Since any parabolic element is conjugate to
a power of $g_x$, $g_y$ or $g_z$, we just need to study one of these examples.

\vv

Once the parameters $A$, $B$, $C$, and $D$ have been fixed, the automorphism
$g_z$ is given by
$$
g_z\left(\begin{array}{c}
          x \\ y \\ z 
         \end{array}
 \right) 
=
\left(\begin{array}{c}
          A - x - zy \\ B-Az +zx + (z^2-1) y \\ z 
         \end{array}
 \right).
$$
This defines a global polynomial diffeomorphism of $\C^3$, that preserves
each horizontal plane $\Pi_{z_0}=\{(x,y,z_0), \, x\in \C, \, y\in \C\}$.  
On each of these planes, $g_z$ induces an affine transformation 
$$
\left(\begin{array}{c}
          x \\ y 
         \end{array}
 \right) \ 
\mapsto\ 
\left(\begin{array}{cc}
          -1 & -z_0 \\ z_0 & z_0^2 -1 
         \end{array}
 \right) \left(\begin{array}{c}
          x \\ y 
         \end{array}
 \right) + \left(\begin{array}{c}
          A \\ B-Az_0 
         \end{array}
 \right),
$$ 
which preserves the conic $S_{z_0}=S\cap \Pi_{z_0}$. The trace of the linear part of this 
affine transformation is $z_0^2-2$ while the determinant is $1$. 

\begin{pro}\label{pro:para_fiber} Let $S$ be any member of the family of cubic surfaces
$\SS$. Let $g_z$ be the automorphism of $S$ defined by
the composition $s_y\circ s_x$. On each fiber $S_{z_0}$ 
of the fibration 
$$
\pi_z:S\to \C, \quad \pi_z(x,y,z)=z,
$$
$g_z$ induces a homographic transformation ${\overline{g_{z_0}}}$, and
\begin{itemize}
 \item ${\overline{g_{z_0}}}$ is an elliptic homography if and only if
${z_0}\in (-2,2)$; this homography is periodic if and only if 
$z_0$ is of type $\pm 2 \cos(\pi \theta)$ with $\theta$ rational;
\item ${\overline{g_{z_0}}}$ is parabolic (or the identity) if and only if ${z_0}=\pm 2$;
\item ${\overline{g_{z_0}}}$ is loxodromic if and only if ${z_0}$ is not in the
interval $[-2,2]$.   
\end{itemize}
\end{pro}

If ${z_0}$ is different from $2$ and $-2$, $g_z$ has a unique fixed point inside $\Pi_{z_0}$,
the coordinate of which are $(x_0,y_0,z_0)$ where
$$
x_0= \frac{B{z_0}-2A}{{z_0}^2-4}, \quad y_0 =\frac{Az_0 -2B}{{z_0}^2-4}.
$$
This fixed point is contained in the surface $S$ if and only if ${z_0}$ 
satisfies the quartic equation $P_z({z_0})=0$ where
\begin{equation}\label{eq:pz}
P_z=z^4 -C z^3 -(D+4)z^2 + (4C-AB)z +4D + A^2 + B^2.
\end{equation}
In that case, the union of the two $g_z$-invariant lines of $\Pi_{z_0}$ which go through the fixed point 
coincides with  $S_{z_0}$; moreover, the involutions $s_x$ and $s_y$ permute
those two lines. If the fixed point is not contained in $S$, the conic $S_{z_0}$ is smooth, and the two fixed points of the (elliptic or loxodromic) homography ${\overline{g_{z_0}}}$ are at infinity.

\vv 

If ${z_0}=2$, the affine transformation induced by $g_z$ on $\Pi_{z_0}$ is 
$$
{\overline{g_{z_0}}}\left(\begin{array}{c}
          x \\ y 
         \end{array}
 \right) 
=
\left(\begin{array}{cc}
          -1 & -2 \\ 2 & 3 
         \end{array}
 \right) \left(\begin{array}{c}
          x \\ y 
         \end{array}
 \right) + \left(\begin{array}{c}
          A \\ B-2A 
         \end{array}
 \right). 
$$ 
Either ${\overline{g_{z_0}}}$ has no fixed point, or $A=B$ and there is a line
of fixed points, given by $x+y=A/2$. This line of fixed points intersects
the surface $S$ if and only if $S_{z_0}$ coincides with this (double) line. 
In that case the involutions $s_x$ and $s_y$ also fix the line pointwise.
When the line does not intersect $S$, the conic $S_{z_0}$ is smooth,
with a unique point at infinity; this point is the unique fixed point of the parabolic 
transformation ${\overline{g_{z_0}}}$. In particular, any point of $S_{z_0}$ goes to infinity 
under the action of $g_{z}$.

\vv

If ${z_0}=-2$, then 
$$
\overline{g_{z_0}}\left(\begin{array}{c}
          x \\ y 
         \end{array}
 \right) 
=
\left(\begin{array}{cc}
          -1 & 2 \\ -2 & 3 
         \end{array}
 \right) \left(\begin{array}{c}
          x \\ y 
         \end{array}
 \right) + \left(\begin{array}{c}
          A \\ B+2A 
         \end{array}
 \right). 
$$ 
Either $g_z$ does not have any fixed point in $\Pi_{z_0}$, or $A=-B$ and $g_{z_0}$ has
a line of fixed points given by $x-y=A/2$. This line intersects $S$
if and only if $S_{z_0}$ coincides with this (double) line. In that case the involutions
$s_x$ and $s_y$ fixe the line pointwise. 

\begin{lem}\label{lem:dyna_para}
With the notation that have just been introduced, the homographic transformation
${\overline{g_{z_0}}}$ induced by $g_z$ on $S_{z_0}$ has a fixed point in $S_{z_0}$ if and only if $z_0$ satisfies equation (\ref{eq:pz}). Moreover
\begin{itemize}
 \item when ${z_0}\neq 2,-2$, $S_{z_0}$ is a singular conic, 
 namely a union of two lines that are permuted by $s_x$ and $s_y$, 
 and the unique fixed point of ${\overline{g_{z_0}}}$ is the point of intersection of these two lines, with coordinates 
$$
x_0= \frac{B{z_0}-2A}{{z_0}^2-4}, \quad y_0 =\frac{Az_0 -2B}{{z_0}^2-4};
$$
\item when ${z_0}=2$,  then $A=B,$ $S_{z_0}$ is the double line $x+y=A/2$,
and this line is pointwise fixed  by ${\overline{g_{z_0}}}$, $s_x$ and $s_y$;
\item when ${z_0}=-2$, then $A=-B,$ $S_{z_0}$ is the double line $x-y=A/2$, and 
 this line is pointwise fixed  by ${\overline{g_{z_0}}}$, $s_x$ and $s_y$;
\end{itemize}
\end{lem}

The dynamics of $g_z$ on $S$ is now easily described. Let $p_0=(x_0,y_0,z_0)$ be a point of $S$. 
If ${z_0}$ is in the interval $(-2,2)$, the orbit of $p_0$ under $g_z$ is bounded, and it is periodic if, and only if, either $p_0$ is a fixed point, or ${z_0}$ is of type $\pm2\cos(\pi\theta),$ where $\theta$ is a rational number. If $z_0=\pm2$, and if $p_0$ is not a fixed point, $g^n(p_0)$ goes to infinity
when $n$ goes to $+\infty$ and $-\infty$. If $z_0$ is not contained in the interval $[-2,2]$,
for instance if the imaginary part of $z_0$ is not $0$, either $p_0$ is fixed or $g^n(p_0)$
goes to infinity when $n$ goes to $-\infty$ or $+\infty.$
Of course, the same kind of results are valid for $g_x$ and $g_y$, with the appropriate
permutation of variables and parameters. 

%%%%%%%%%%%%%%%%%%%%%%%%%%%%%%%%%%%%%%%%%%%%%%%

\subsection{Bounded Orbits}\label{par:subsection_bounded}
There is a huge literature on the classification of algebraic
solutions of Painlev\'e VI equation (see \cite{Boalch:2007} and references therein).
Such solutions give rise to periodic orbits for 
the action of $\A$ on the cubic surface $S_{(A,B,C,D)}$, where the parameters are defined
in terms of the coefficients of the Painlev\'e equation (see \S \ref{par:appendixB}). Of course,
periodic orbits are bounded. Here, we study infinite bounded  orbits.

\begin{thm}\label{Thm:BoundedOrbits}
Let $S=S_{(A,B,C,D)}$ be a surface in the family $\SS$, and $p$ be a point with 
an infinite and bounded $\Gamma_2^*$ orbit  $\Orb(p)$. 
%If $\#\Orb(p)>4$, 
%then the parameters $A$, $B$, $C$, and $D$ that define $S$ are %real numbers and the orbit of $p$ is contained in $S(\R)$. 
Then $A,$ $B,$ $C,$ and $D$ are real numbers, the orbit  is contained
in $[-2,2]^3$ and it  forms a dense subset of the unique bounded 
component of $S(\R)\setminus {\sf{Sing}}(S)$.
\end{thm}

We fix a point $p$ in one of the surfaces  $S$ 
and denote its $\Gamma^*(2)$-orbit by $\Orb(p)$. 
Let us first study orbits of small finite length.
Recall that orbits of length $1$ are singular points 
of the cubic $S$.

\begin{pro}\label{Prop:SmallOrbits}
Modulo Benedetto-Goldman symmetries (see \S \ref{par:twists}), $\Gamma_2^*$-orbits 
of length $2$ are equivalent to
$$
\{(0,0,z_1),(0,0,z_2)\}\in S_{(0,0,C,D)},\ \ \ C^2+4D\neq 0
$$
where $z_1$ and $z_2$ are the two roots of $z^2=Cz+D,$
$\Gamma_2^*$-orbits of length $3$ are equivalent to
$$
\{(0,0,1),(A,0,1),(0,A,1)\}\in S_{(A,A,2,-1)},
$$
and $\Gamma_2^*$-orbits of length $4$ are equivalent to
$$
\{(1,1,1),(A-2,1,1),(1,A-2,1),(1,1,A-2)\}\in S_{(A,A,A,4-3A)}.
$$
\end{pro}

%\begin{itemize}
%\item $A=B=0$, $x=y=0$, $z^2=Cz+D$, and $C^2+4D\neq 0$. 
%\item $B=C=0$, $y=z=0$, $x^2=Ax+D$, and $A^2+4D\neq 0$. 
%\item $C=A=0$, $z=x=0$, $y^2=By+D$, and $B^2+4D\neq 0$. 
%\end{itemize}

\begin{eg}An orbit of length $2$ is for instance provided by
the representation $\rho$ defined by 
$$\rho\ :\ (\alpha,\beta,\gamma,\delta)\mapsto(M,N,M,-N)$$
where $M,N\in\SL(2,\C)$ are any element satisfying $\Tr(MN)=0$
i.e. $(MN)^2=-I$.
Trace parameters are given by $(a,b,a,-b)$ where 
$a=\Tr(M)$ and $b=\Tr(N)$ : we get $C=a^2-b^2$, $D=(a^2-2)(b^2-2)$
and $z=a^2-2$.
The other representation in the orbit, given by $z'=2-b^2$, 
is defined by 
$$\rho'\ :\ (\alpha,\beta,\gamma,\delta)\mapsto(M,M^{-1}NM,NMN^{-1},-N).$$
To this length $2$ orbit corresponds a two-sheeted algebraic 
solution of $P_{VI}$-equation, namely 
$$q(t)=1+\sqrt{1-t},\ \ \ \text{for parameters}\ \ \ \theta=(\theta_0,\theta_1,\theta_0,-\theta_1),$$
with $a=2\cos(\pi\theta_0)$ and $b=2\cos(\pi\theta_1)$.
This representation was already considered in \cite{Previte-Xia:2003} :
for convenient choice of parameters $a$ and $b$, the image
of the representation is a dense subgroup of $\SU(2)$.

Other choice of the trace parameters are provided by 
$$(a',b',a',-b')\ \ \ \text{with}\ a'=\sqrt{4-b^2}\ \text{and}\ b'=\sqrt{4-a^2},$$
giving rise to a representation of the same kind, and
$$(a'',0,c'',0)\ \ \ \text{with}\ a'',c''=\frac{a}{2}\sqrt{4-b^2}\pm\frac{b}{2}\sqrt{4-a^2}.$$
The later one corresponds to a dihedral representation of the form
$$(\alpha,\beta,\gamma,\delta)\mapsto
\begin{pmatrix}\lambda&0\\0&\lambda^{-1}\end{pmatrix},
\begin{pmatrix}0&\mu\\ -\mu^{-1}&0\end{pmatrix},
\begin{pmatrix}\tau^{-1}&0\\0&\tau\end{pmatrix},
\begin{pmatrix}0&-\nu^{-1}\\ \nu&0\end{pmatrix})$$
with $\lambda\mu\nu\tau=1$.
\end{eg}

\begin{proof}  Let $p=(x_0,y_0,z_0)$ be a point of $S_{(A,B,C,D)}$.
Recall that $p$ is fixed if, and only if, $p$ is a singular point of $S.$
%we are in one 
%of the following cases:
%$$\left\{\begin{matrix}
%P_z(z_0)=0\\ 2x_0+y_0z_0=A\\ 2y_0+x_0z_0=B
%\end{matrix}\right.
%\ \ \ \text{or}\ \ \ 
%\left\{\begin{matrix}
%z_0=2\\ 2(x_0+y_0)=A\\ A=B
%\end{matrix}\right.
%\ \ \ \text{or}\ \ \ 
%\left\{\begin{matrix}
%z_0=-2\\ 2(x_0-y_0)=A\\ A=-B
%\end{matrix}\right.$$
On the other hand, $p$ is periodic of order $n>1$ for $g_z$ if,
and only if, 
$$z_0=2\cos(\pi \frac{k}{n}),\ \ \ k\wedge n=1$$
and at least one of the equalities $P_z(z_0)=0$, $2x_0+y_0z_0=A$, $2y_0+x_0z_0=B$ does not hold. In particular, denoting by 
$\Orb_{g_z}(p)$ the orbit of $p$ under the action of $g_z$, we have:
$$\# \Orb_{g_z}(p)=2\ \Rightarrow\ z_0=0,$$
$$\# \Orb_{g_z}(p)=3\ \Rightarrow\  z_0=\pm 1,$$
$$\# \Orb_{g_z}(p)=4\ \Rightarrow\  z_0=\pm \sqrt{2},$$
$$\# \Orb_{g_z}(p)=6\ \Rightarrow\  z_0=\pm \sqrt{3}.$$

Up to permutation of variables $x$, $y$ and $z$ 
(and correspondingly of the parameters $A$, $B$ and $C$), 
an orbit of length $2$
takes the form $\Orb(p)=\{p,s_z(p)\}$.
In this case, $p$ and $p'=s_z(p)=(x_0,y_0,z_0')$ 
are permuted by $s_z$, and thus by $g_x=s_z\circ s_y$ and $g_y=s_x\circ s_z$ ;
this implies $x_0=y_0=0$. 
On the other hand, $p$ and $p'$ are fixed by $s_x$, $s_y,$ 
and therefore $A=B=0.$ Since $p=(0,0,z_0)$ is contained in $S,$ we deduce that $z_0$ and $z_0'$  are the roots of $z^2=Cz+D$.

Up to permutation of the variables $x$, $y$ and $z$, 
an orbit of length $3$ takes the form 
$$\Orb(p)=\{p_0,p_1=s_x(p_0),p_2=s_y(p_0)\}$$
with $p_0=(x_0,y_0,z_0),$ $p_1=(x_0',y_0,z_0),$ and $p_2=(x_0,y_0',z_0).$
Since $g_x$ (resp. $g_y$) permutes $p_0$ and $p_2$ (resp. $p_0$ and $p_1$),
we get $x_0=y_0=0$. On the other hand, $g_z$ permutes cyclically $p_0\to p_2\to p_1,$
so that $z_0=\pm 1$. Changing signs if necessary by a twist (see \ref{par:twists}),
we can assume $z_0=1$. Now, studying the fixed points of $s_x$, $s_y$ and $s_z$,
amongst  $p_0$, $p_1$ and $p_2,$ we obtain:
$$2z_0+x_0y_0=C,\ \ \ 
\left\{\begin{matrix}2y_0+x_0'z_0=B\\ 2z_0+x_0'y_0=C\end{matrix}\right.\ \ \ \text{and}\ \ \ 
\left\{\begin{matrix}2x_0+y_0'z_0=A\\ 2z_0+x_0y_0'=C\end{matrix}\right.$$
and thus $C=2$, $x_0'=B$ and $y_0'=A$. We also have
$$x_0'=A-x_0-y_0z_0\ \ \ \text{and}\ \ \ y_0'=B-y_0-x_0z_0$$
(action of $s_x$ and $s_y$) yielding $A=B$. Finally, 
the fact that $p_0$ is contained in $S$ gives $1=C+D$, whence the result.

Up to symmetry, an orbit of length $4$ consists in  $p_0$, $p_1$
and $p_2$ like before ($p_1=s_x(p_0)$ and $p_2=s_y(p_0)$) 
and there are  $4$ possibilities for the fourth point $p_3$:
\begin{itemize}
\item[(1)] $p_3=(x_0',y_0',z_0)=s_y(p_1)=s_x(p_2)$,
\item[(2)]  $p_3=(x_0',y_0'',z_0)=s_y(p_1)\not=s_x(p_2)$,
\item[(3)]  $p_3=(x_0',y_0,z_0')=s_z(p_1)$,
\item[(4)]  $p_3=(x_0,y_0,z_0')=s_z(p_0)$.
\end{itemize}
The first case is impossible: Since $g_x$ and $g_y$  have order $2$
for each $p_i,$  the coordinates $x$ and $y$ vanish for each point $p_i,$ and therefore $p_0=p_1=p_2=p_3$, a contradiction. The second case is impossible for the same reason since  $g_x$ and $g_y$ have order $2$ for $p_0$ and $p_1,$ 
so that $p_0=p_1$, a contradiction. The same argument applies in third case:
$g_x$ has order $2$ for $p_0$ and $p_1$ implying $p_0=p_1$, contradiction.

For the fourth case, since $g_x$, $g_y$ and $g_z$ have order $3$ for $p_0$,
we get $p_0=(\pm1,\pm1,\pm1)$. Up to symmetry, there are two subcases:
$p_0=(1,1,1)$ or $p_0=(-1,-1,-1)$. In the first sub-case, 
conditions given by the fixed points of $s_x$, $s_y$ and $s_z$ yield 
$$A=B=C=2+x_0'=2+y_0'=2+z_0',$$
and the fact that $p_i$ is contained in $S$ gives $4=3A+D$. 
Proceeding in the same way with the second sub-case, 
conditions given by the fixed points of $s_x$, $s_y$ and $s_z$ yield 
$$A=B=C=-2-x_0'=-2-y_0'=-2-z_0'$$
and the fact that $p_i$ is in $S$ gives
$$2=-3A+D\ \ \ \text{and}\ \ \ x_0'=A$$
implying $x_0'=A=-1$ and $p_1=p_0$, a contradiction.
\end{proof}

\begin{lem}\label{Lem:Bounded+Large=Real}
If $\Orb(p)$ is bounded and $\#\Orb(p) > 4$, then $A$,  $B$, $C$, and $D$
are real and $p\in S(\R)$.
\end{lem}

\begin{proof}
Let $p_0=(x_0,y_0,z_0)$ be a point of the orbit.
If the third coordinate $z_0\not\in(-2,2)$, the homography induced 
by $g_z$ on the conic $S_{z_0}$ is parabolic or hyperbolic. 
Since the orbit of $p_0$ is bounded, this implies that $p_0$ 
is a fixed point of $g_z$, $s_x$ and $s_y$ (see lemma \ref{lem:dyna_para}).
Since  $\Orb(p_0)$ has length $>4,$ $s_z(p_0)$ is different from $p_0$,
so that $p_0$   is not fixed by $g_x$, nor by $g_y$ either ; this implies 
that $x_0,y_0\in(-2,2)$. 
Moreover, the point $p_1:=s_z(p_0)=(x_0,y_0,z_1)$ is 
not fixed by $g_z$, otherwise the orbit would have length $2,$ so that
$z_1\in(-2,2)$ and $p_1\in(-2,2)^3$. This argument shows the following: 
If one of the coordinates of $p_0$ is not contained in $(-2,2),$ then $p_0$ 
is fixed by two of the involutions $s_x,$ $s_y$ and $s_z$ while the third
one maps $p_0$ into $(-2,2)^3.$ 

Let now $p$ be a point of the orbit with coordinates in $(-2,2)^3$; 
if  the three points $s_x(p),$  $s_y(p)$ and 
$s_z(p)$ either escape from $(-2,2)^3$ or coincide with $p$, 
then the orbit reduces to $\{p, \, s_x(p), \, s_y(p),\, s_z(p)\},$ and has length $\le 4$.
From this we deduce that the orbit contains at least
two distinct points $p_1,p_2\in(-2,2)^3,$ which are, say, permuted by $s_x$.
Let $(x_i,y_1,z_1)$ be the coordinates of $p_i,$ $i=1,2.$
Then, $A=x_1+x_2+y_1z_1\in\mathbb R$. 
If $B$ and $C$ are also real, then
$p_1$ is real and satisfies the equation of $S$, so that $D$ is real
as well and $\Orb(m)=\Orb(p_1)\subset S(\R)$. 

Now, assume by contradiction that $B\not\in\R$. Then,
$q_i:=s_y(p_i)=(x_i,B-y_1-x_iz_1,z_1)\not\in(-2,2)$
and is therefore fixed by $s_x$ (otherwise the orbit
would not be bounded): We thus have 
$$2x_i+(B-y_1-x_iz_1)z_1=A.$$
Since $B$ is the unique imaginary number of this equation, $z_1$ must vanish,
and we get $x_1=x_2(=\frac{A}{2})$, a contradiction.
A similar argument shows that $C$ must be real as well.
\end{proof}

\begin{pro}\label{Prop:AlgebraicPeriodic}
If $\Orb(p)$ is finite and $\#\Orb(p) > 4$, then $A$,  $B$, $C$, and $D$
are real algebraic numbers and $p\in S(\R)$ has algebraic coordinates as well.
\end{pro}

The  proof is exactly the same, replacing $(-2,2)$ by $(-2,2)\cap 2\cos(\pi \Q)$
and thus $\R$ by $\R\cap {\overline{\Q}}.$

\begin{lem}\label{Lem:Selberg}
Let $S$ be an element of the family $\SS$ and $p$ a point of $S$. There exists a positive integer
$N$ such that, if $p'$ is a point of the orbit of $p$ with a coordinate of the form
$$
2 \cos (\pi\frac{k}{n}), \quad k\wedge n =1,
$$
then $n$ divides $N$.
\end{lem}

\begin{proof}
The point $p$ is an element of the character variety $\chi(\Sphere_4)$. Let us
choose a representation $\rho:\pi_1(\Sphere_4)\to \SL(2,\C)$ in the  conjugacy
class that is determined by $p$. Since $\pi_1(\Sphere_4)$ is finitely generated,  Selberg's lemma (see \cite{Alperin:1987})
implies the existence of  a torsion free, finite
index subgroup $G$ of $\rho(\pi_1(\Sphere_4))$. If we define $N$ to be
the cardinal of the quotient $\rho(\pi_1(\Sphere_4))/G$, then the order of any torsion element
in $\rho(\pi_1(\Sphere_4))$ divides $N$. 

If $p'$ is a point of the orbit of $p$, the coordinates of $p'$ are traces of elements
of $\rho(\pi_1(\Sphere_4))$. Assume that the trace of an element $M$ in $\rho(\pi_1(\Sphere_4))$ 
is of type 
$
2 \cos(\pi\theta).
$
If $\theta=\frac{k}{n}$ and $k$ and $n$ are relatively prime integers, then $M$ is a cyclic element 
of $\rho(\pi_1(\Sphere_4))$ of order $n$, so that $n$ divides $N$. \end{proof}

The subset of $\SU(2)$-representations always
form a connected component of $S\setminus \Sing(S)$
contained into $[-2,2]^3$; the corresponding orbits
are bounded, generally infinite. A bounded component
can also consist in $\SL(2,\R)$-repre\-sentations, depending
on the choice of $(a,b,c,d)$; for instance, in the Cayley
case, the bounded component consists in 
$\SL(2,\R)$-representations (resp. $\SU(2)$-representations)
when $(a,b,c,d)=(2,2,2,-2)$ (resp. $(0,0,0,0)$). 

\begin{pro}[Benedetto-Goldman \cite{Benedetto-Goldman:1999}]\label{Prop:BenedettoGoldman} When $A$, $B$, $C$ and $D$
are real, then $S(\R)\setminus\{\Sing(S)\}$ has at most one 
bounded connected component. In this case, 
$a$, $b$, $c$ and $d$ lie in $[-2,2]$, 
whatever the choice of $(a,b,c,d)$
corresponding to $(A,B,C,D)$.
\end{pro}

When $S(\R)$ is smooth, the converse is true: When 
$a$, $b$, $c$ and $d$ lie in $[-2,2]$, $S(\R)$ has 
a ``bounded component'' maybe degenerating to
a singular point, like in the Markov case.
It is proved in Apendix B, \S \ref{section:SU(2)}, 
that a bounded component 
always correponds to $\SU(2)$-representations for a convenient
choice of parameters $(a,b,c,d)$.

\begin{proof}[Proof of theorem \ref{Thm:BoundedOrbits}]
Let $\Orb(p)$ be an infinite and bounded $\Gamma_2^*$-orbit
in $S=S_{(A,B,C,D)}$.
Following Lemma \ref{Lem:Bounded+Large=Real},
$A$, $B$, $C$ and $D$ are real and $\Orb(p)\subset S(\R)$.
We want to prove that the closure
$\overline{\Orb(p)}$ is open in $S(\R)\setminus\{\Sing(S(\R))\}$;
since $\overline{\Orb(p)}$ is closed, it will therefore coincide 
with the (unique) bounded connected component of
$S\setminus\{\Sing(S)\}$, thus proving the theorem.

We first claim that there exists an element (actually infinitely many)
$p_0=(x_0,y_0,z_0)$ of the orbit which is contained in $(-2,2)^3$
and for which at least one of the M\"obius transformations 
$\overline{g_{x_0}}$, $\overline{g_{y_0}}$ or $\overline{g_{z_0}}$
is (elliptic) non periodic. Indeed, if a point $p_0$ of the orbit 
is such that $\overline{g_{z_0}}$ is not of the form above,
then we are in one of the following cases
\begin{itemize}
\item $P_z(z_0)=0$ and $p_0$ is a fixed point of $\overline{g_{z_0}}$,
\item $z_0=2\cos(\pi\frac{k}{n})$ with $k\wedge n=1$, $n\vert N$ and $\overline{g_{z_0}}$ is periodic of period $n$
\end{itemize}
(where $N$ is given by Lemma \ref{Lem:Selberg}).
This gives us finitely many possibilities for $z_0$; we also get 
finitely many possibilities for $x_0$ and $y_0$ and the claim follows.

Let $p_0$ be a point of $\Orb(p)$, with, say, $\overline{g_{x_0}}$
elliptic and non periodic, so that the closure $\overline{\Orb(p)}$ contains
the "circle"
$\overline{\Orb_{g_x}(p_0)}=S_{x_0}(\R).$
 Let us first prove that
$\overline{\Orb(p)}$ contains an open neighborhood of $p_0$ in $S(\R)\setminus\{\Sing(S(\R))\} $. 

Since the point $p_0$ is not fixed by $g_x=s_z\circ s_y$, 
then either $s_y$ or $s_z$ does not fix $p_0$, say $s_z$;
this means that the point $p_0$ is not a critical point
of the projection 
$$
\pi_x\times\pi_y\ :\ S(\R)\to\R^2\ ;\ (x,y,z)\mapsto(x,y).
$$
Therefore, there exists some $\varepsilon>0$ and a neighborhood
$V_\varepsilon$ of $p_0$ in $S(\R)$ such that $\pi_x\times\pi_y$ 
maps $V_\varepsilon$ diffeomorphically onto the square
$(x_0-\varepsilon,x_0+\varepsilon)\times(y_0-\varepsilon,y_0+\varepsilon)$. By construction, we have
$$
\pi_x\times\pi_y(\overline{\Orb(p)})\supset 
\pi_x\times\pi_y(\overline{\Orb_{g_x}(p_0)})\supset
\{x_0\}\times(y_0-\varepsilon,y_0+\varepsilon).
$$
For each $y_1\in(y_0-\varepsilon,y_0+\varepsilon)$ 
 of irrational type, that is to say not of the form 
$2\cos(\pi\theta)$ with $\theta$ rational,
there exists
$p_1=(x_0,y_1,z_1)\in{\overline{\Orb(p)}}$ (namely, the preimage of 
$(x_0,y_1)$ by $\pi_x\times\pi_y$) 
and 
$$
\overline{\Orb(p)}\supset\overline{\Orb_{g_y}(p_1)}=S_{y_1}(\R);
$$
in other words, for each $y_1\in(y_0-\varepsilon,y_0+\varepsilon)$ 
of irrational type, we have
$$
\pi_x\times\pi_y(\overline{\Orb(p)})\supset 
\pi_x\times\pi_y(\overline{\Orb_{g_y}(p_1)})\supset
(x_0-\varepsilon,x_0+\varepsilon)\times\{y_0\}.
$$
Since those coordinates $y_1$ of irrational type
are dense in $(y_0-\varepsilon,y_0+\varepsilon)$, 
we deduce that 
$V_\varepsilon\subset\overline{\Orb(p)}$, 
and $\overline{\Orb(p)}$ is open at $p_0$.

It remains to prove that $\overline{\Orb(p)}$ is open 
at any point $q\in\overline{\Orb(p)}$ which is not a singular
point of $S(\R).$ 
Let $q=(x_0,y_0,z_0)$ be such a point and assume that
$q\not\in\Orb(p)$ (otherwise we have already proved the
assertion).

Since $q$ is not a singular point of $S(\R)$, one of the projections, 
say $\pi_x\times\pi_y\ :\ S(\R)\to\R^2$, is regular at $q$
and we consider a neighborhood $V_\varepsilon$ like above,
$\pi_x\times\pi_y(V_\varepsilon)=(x_0-\varepsilon,x_0+\varepsilon)\times(y_0-\varepsilon,y_0+\varepsilon)$. By assumption, 
$\Orb(p)\cap V_\varepsilon$ is infinite (accumulating $q$)
and, applying once again Lemma \ref{Lem:Selberg}, 
one can find one such point 
$p_1=(x_1,y_1,z_1)\in\Orb(p)\cap V_\varepsilon$
such that either $x_1$ or $y_1$ has irrational type, say $x_1$.
Now, reasonning with $p_1$ like we did above with $p_0$, 
we eventually conclude that $V_\varepsilon\supset\overline{\Orb(p)}$,
and $\overline{\Orb(p)}$ is open at $q$.
\end{proof}

\section{Invariant geometric structures}

%%%%%%%%%%%%%%%%%%%%%%%%%%%%%%%%%%%%%%%%

In this section, we study the existence of $\A$-invariant geometric structures
on surfaces $S$ of the family $\SS.$ 
An example of such an invariant structure is given by the
area form $\Omega,$ defined in Proposition 
\ref{pro:AreaForm}. Another example occurs
for the Cayley cubic: $S_C$ is covered by $\C^*\times \C^*$ and the action of
$\A$ on $S_C$ is covered by the monomial action of $\GL(2,\Z),$ that is also covered by the linear
action of $\GL(2,\Z)$ on $\C\times \C$ if we use the covering mapping
$$
\pi:\C\times \C\to \C^*\times \C^*, \quad \pi(\theta, \phi)=(\exp(\theta),\exp(\phi));
$$
as a consequence, there is an obvious $\A$-invariant affine structure on $S_C.$

\begin{rem}
The surface $S_C$ is endowed with a natural orbifold structure, the 
analytic structure near its singular points being locally isomorphic to
the quotient of $\C^2$ near the origin by the involution $\sigma(x,y)=(-x,-y).$
The affine structure can be understood either in the orbifold language, or as an 
affine structure defined only outside the singularities (see below).
\end{rem}

%%%%%%%%%%%%%%%%%%%%%%%%%%%%%%%%%%%%%%%%%%%%%%%%%%%%%

\subsection{Invariant curves, foliations and webs}\label{par:inv_curves_webs}

%%%%%%%%%%%%%%%%%%%%%%%%%%%%%%%%%%%%%%%%%%%%%%%%%%%%%

We start with

\begin{lem}\label{Lem:NoCurve}
Whatever the choice of $S$ in the family $\SS,$ the group $\A$ does not preserve any (affine) algebraic curve on $S.$
\end{lem}

Of course, invariant curves appear if we blow up  singularities. This is important for the study 
of special (Riccati) solutions of Painlev\'e VI equation (see section \ref{par:Painleve}).

\begin{proof}
Let $C$ be an algebraic curve on $S.$ Either $C$ is contained in a fiber of $\pi_z,$ 
or the projection $\pi_z(C)$ covers $\C$ minus at most finitely many points. If $C$ 
is not contained in a fiber, we can choose $m_0=(x_0,y_0,z_0)$ in $C$ and
a neighborhood $U$ of $m_0$ such that  $z_0$ is contained in  $(0,2)$ and, 
in $U,$ $C$ intersects each fiber $S_z$ of the projection $\pi_z$ in exactly one point. 
Let $m'=(x',y',z')$ be any element of $C\cap U$ such that $z'$ is an element of $(0,2).$
Then ${\overline{g_z}}$ is an elliptic transformation of $S_{z'}$ that preserves $C\cap S_{z'}$; since the intersection
of $C$ and $S_{z'}$ contains an isolated point $m',$ this point is ${{g_z}}$ periodic. As a consequence, $z'$ is of the form $2\cos(\pi p/q)$ (see proposition \ref{pro:para_fiber}). Since any 
$z'\in (0,2)$ sufficiently close to $z_0$ should satisfy an  equation of this type, we obtain a contradiction.

Since no curve can be simultaneously contained in fibers of $\pi_x,$ $\pi_y$ and $\pi_z,$
the lemma is proved.
\end{proof}

A (singular) web on a surface $X$ is given 
by a hypersurface in the projectivized  tangent bundle $\P TX$; for each point, the web
determines a finite collection of directions tangent to $X$ through that point. 
The number of directions is constant on an open subset of $X$ but it may vary  along 
the singular locus of the web. Foliations are particular cases of webs, and any web is locally
made of a finite collection of foliations in the complement of its singular locus.

\begin{pro}\label{Prop:NoWeb}
Whatever the choice of $S$ in the family $\SS,$ the group $\A$ does not preserve any web on $S.$
\end{pro}

\begin{proof}
Let us suppose that there exists an invariant web $W$ on one of the surfaces $S.$
Let $k$ and $l$ be coprime positive integers and $m=(x,y,z)$ be a periodic point of $g_z$  of period $l,$ 
with 
$$
z=2\cos(\pi k/l).
$$
Let $L_1,$ ..., $L_d$ be the directions determined by $W$ through the point $m,$ and $C_1,$ ..., $C_d$ 
the local leaves of $W$ which are tangent to these directions. The automorphism $g_z^{s},$ with $s= l (d!),$ fixes $m,$ preserves
the web and fixes each of the directions $L_i$; it therefore preserves each of the $C_i.$ The proof 
of lemma \ref{Lem:NoCurve} now shows that $d=1$ and that the curves $C_i$ are contained in the
fiber of $\pi_z$ through $m.$ Since the set of points $m$ which are $g_z$-periodic is Zariski dense
in $S,$ this argument shows that the web is the foliation by fibers of $\pi_z.$ The same 
argument shows that the web should also coincide with the foliations by fibers of $\pi_x$ or $\pi_y,$ a contradiction.
\end{proof}

\begin{cor}
Whatever the choice of $S$ in the family $\SS,$ the group $\A$ does not preserve any holomorphic riemannian
metric on $S.$
\end{cor}

\begin{proof} Let $\sf{g}$ be an invariant holomorphic riemannian metric. At each point $m$ 
of $S,$ ${\sf{g}}$ has two isotropic lines. This determines an $\A$-invariant web, and we get 
a contradiction with the previous proposition. 
\end{proof}

%%%%%%%%%%%%%%%%%%%%%%%%%%%%%%%%%%%%%%%%%%%%%%%%%%%%%%%%%%%%%%%%%%%%%%%%%%%%%%%%%

\subsection{Invariant Affine Structures}\label{par:InvAffStruct}

%%%%%%%%%%%%%%%%%%%%%%%%%%%%%%%%%%%%%%%%%%%%%%%%%%%%%

A holomorphic affine structure on a complex surface $M$ is given by an atlas of charts $\Phi_i:U_i\to \C^2$ for which the transition functions 
$\Phi_i \circ \Phi_j^{-1}$ are affine transformations of the plane $\C^2.$ A local chart $\Phi:U\to \C^2$ is said to be affine if, for all $i,$ $\Phi\circ \Phi_i^{-1}$
is the restriction of an affine transformation of $\C^2$ to $\Phi_i(U_i)\cap \Phi(U).$  
A subgroup $G$ of $\Aut(M)$ preserves the affine
structure if elements of $G$ are given by affine transformations in  local affine charts. 

\begin{thm}\label{thm:affine}
Let $S$ be an element of $\SS.$ Let $G$ be a finite index subgroup of $\Aut(S).$
The group $G$ preserves an  affine structure on $S\setminus \Sing(S),$ if, and only if
$S$ is the Cayley cubic $S_C.$ \end{thm}

In what follows, $S$ is a cubic of the family $\SS$ and $G$ will
be a finite index subgroup of $\A$ preserving an affine structure on $S.$

Before giving the proof of this statement, we collect a few basic results concerning
affine structures. Let $X$ be a complex surface with a holomorphic  affine structure.
Let $\pi:{\widetilde{X}}\to X$ be the  universal cover of $X$; the group
of deck transformations of this covering is isomorphic to the  fundamental group
$\pi_1(X).$ Gluing together the affine local charts of $X,$ 
we get a developping map 
$$
\dev:{\widetilde{X}}\to\C^2,
$$ 
and a monodromy representation $\hol:\pi_1(X)\to \Aff(\C^2)$ such that 
$$
\dev(\gamma(m))= \hol(\gamma)(\dev(m))
$$
for all $\gamma$ in $\pi_1(X)$ and all $m$ in ${\widetilde{X}}.$
The map $\dev$ is a local diffeomorphism but, a priori, it is neither surjective,
nor a covering onto its image. 

Let $f$ be an element of  $\Aut(X)$ that preserves the affine structure
of $X.$ Let $m_0$ be a fixed point of $f,$ let  ${\widetilde{m_0}}$
be an element of the fiber $\pi^{-1}(m_0),$ and let ${\widetilde{f}}:{\widetilde{X}}\to {\widetilde{X}}$
be the lift of $f$ that fixes ${\widetilde{m_0}}.$ Since $f$ is  affine, 
there exists a unique affine automorphism $\Af{(f)}$ of $\C^2$ such that
$$
\dev \circ {\widetilde{f}} = {\Af{(f)}} \circ \dev.
$$

%%%%%%%%%%%%%%%%%%%%%%%%%%%%%

\subsection{Proof of theorem \ref{thm:affine}; step 1.} In this first step, we show that 
$S\setminus \Sing(S)$ cannot be simply connected, and deduce from this fact
that $S$ is singular. Then we study the singularities of $S$ and the  fundamental
group of $S\setminus \Sing(S).$ 

%%%%%%%
\subsubsection{Simple connectedness} Assume that $S\setminus \Sing(S)$ is simply connected. The developping map $\dev$ is therefore defined on 
$S\setminus \Sing(S)\to \C^2.$
Let $N$ be a positive integer for which $g_x^N$ is contained
in $G.$ Choose a fixed point $m_0$ of $g_x$ as a base point. Since $g_x^N$ preserves the affine structure, there exists an affine transformation
$\Af (g_x^N)$ such that 
$$
\dev \circ g_x^N = \Af(g_x^N)\circ \dev. 
$$
In particular, $\dev$ sends periodic points of $g_x^N$ to periodic points of $\Af(g_x^N).$ 
Let $m$ be a nonsingular point of $S$ with its first coordinate in the interval $(-2,2),$ 
and let $U$ be an open neighborhood of $m.$ 
Section \ref{par:dynamics_parabolic} shows that periodic points of $g_x^N$ 
form a Zariski-dense subset of $U,$ by which we mean that 
any holomorphic functions $\Psi:U\to \C$ which vanishes on the set of
periodic points of $g_x^N$ vanishes everywhere. Since $\dev$ is a local 
diffeomorphism, periodic points of $\Af(g_x^N)$ are Zariski-dense in a neighborhood
of $\dev(m),$ and therefore $\Af(g_x^N)={\text{Id}}.$ This provides a contradiction, 
and shows that $S\setminus \Sing(S)$ is not simply connected. 

Consequently, lemma \ref{Lem:TopologyCubic} 
implies that $S$ is singular and that the fundamental group of $S\setminus \Sing(S)$ is generated,
as a normal subgroup, by the local fundamental groups around the singularities. 

%%%%%%%
\subsubsection{Orbifold structure}\label{par:orbifoldstructure} We already explained in section \ref{par:sing_fp_orb} that the singularities
of $S$ are quotient singularities. If $q$ is a singular point of $S,$ $S$ is locally isomorphic
to the quotient of the unit ball $\B$  in $\C^2$ by a finite subgroup $H$ of $\SU(2).$ 

The local affine structure around $q$ can therefore be lifted into a $H$-invariant affine structure on 
 $\B\setminus \{(0,0)\},$ and then extended up to the origin by Hartogs theorem. In particular, 
$\dev$ lifts to a local diffeomorphism between $\B$ and an open subset of $\C^2.$ This remark 
shows that the affine structure is compatible with the orbifold structure of $S$ defined in section \ref{par:sing_fp_orb}. 

Let $h$ be an element of the local fundamental group $H.$ Let us lift
the affine structure on $\B$ and assume that the monodromy action of $h$ is trivial, i.e. $\dev \circ h = \dev.$
Since $\dev$ is a local diffeomorphism, the singularity is isomorphic to a quotient of $\B$ by a proper
quotient of $H,$ namely the quotient of $H$ by the smallest normal subgroup containing $h.$ This provides
a contradiction and shows that $(i)$ $H$ embeds in the global fundamental group of $S\setminus \Sing(S)$
and $(ii)$ the universal cover of $S$ in the orbifold sense is smooth (it is obtained by adding points to
 the universal cover of ${S\setminus \Sing(S)}$ above singularities of $S$). 

In what follows, we denote the orbifold universal cover by $\pi:{\widetilde{S}}\to S,$ and the developing map
by $\dev:{\widetilde{S}}\to \C^2.$

%%%%%%%
\subsubsection{Singularities}\label{par:node} Let $q$ be a singular point of $S.$ Let ${\widetilde{q}}$ be a point of the fiber $\pi^{-1}(q).$ 
Since the group $\A$ fixes all the singularities of $S,$ it fixes $q$ and one can lift the action of $\A$ on $S$ to an action 
of $\A$ on ${\widetilde{S}}$ that fixes ${\widetilde{q}}.$ If $f$ is an element of $\A,$ $\widetilde{f}$ will denote the corresponding
holomorphic diffeomorphism of ${\widetilde{S}}.$ Then we compose $\dev$ by a translation of the affine plane $\C^2$ in order
to assume that 
$$
\dev({\widetilde{q}})=(0,0).
$$
By assumption, $\dev\circ {\widetilde{g}}=\Af(g)\circ \dev$ for any element $g$ in $G,$ from which we deduce that 
the affine transformation $\Af(g)$ are in fact linear. Since $\A$ almost preserves an area form, 
$\Af(g)$ is an element of $\GL(2,\C)$ with determinant $+1$ or $-1$; passing to a subgroup of 
index $2$ in $G,$ we shall assume that the determinant is $1.$  Since $\dev$ realizes a local 
conjugation between the action of $G$ near ${\widetilde{q}}$ and the action of $\Af(G)$ 
near the origin, the morphism
$$
\left\{ \begin{array}{rcl} G & \to & \SL(2,\C) \\ g & \mapsto & \Af(g) \end{array} \right.
$$
is injective.

\vv

Since $G$ is a finite index subgroup of $\Aut(S),$ $G$ contains a 
non abelian free  group of finite index and is not virtually solvable.
Let $H$ be the finite subgroup of $\pi_1(S\setminus \Sing(S))$ that fixes the point ${\widetilde{q}}.$ 
This group is normalized by the action of $\A$ on ${\widetilde{S}}.$ Consequently, using the local 
affine chart determined by $\dev,$ the group $\Af(G)$
normalizes the monodromy group $\hol(H).$ If $\hol(H)$ is not contained in the center of $\SL(2,\C),$
 the eigenlines of the elements of  $\hol(H)$ determine a  finite, non empty, and $\Af(G)$-invariant set of lines in $\C^2,$ so that $\Af(G)$ is virtually solvable.  This would contradict the injectivity
 of $g\mapsto \Af(g).$ 
 From this we deduce that any element of $\hol(H)$ is a homothety with determinant $1.$
Since the monodromy representation is injective on $H,$ we conclude that  $H$ "coincides" with the subgroup $\{+{\text{Id}},-{\text{Id}}\}$ of $\SU(2).$

%%%%%%%
\subsubsection{Linear part of the monodromy}\label{par:lpotm}
By lemma \ref{Lem:TopologyCubic}, the fundamental group of $S\setminus \Sing(S)$ is generated, as a normal subgroup, by 
the finite fundamental groups around the singularities of $S.$ Since $\pm {\text{Id}}$ is in the center of $\GL(2,\C),$
the linear part of the monodromy $\hol(\gamma)$ of any element $\gamma$ in $\pi_1(S\setminus \Sing(S))$ is
equal to $+{\text{Id}}$ or $-{\text{Id}}.$ 

%%%%%%%%%%%%%%%%%%%%%%%%%%%%%%%%%%%%%%%%%%%%%%%%

\subsection{Proof of theorem \ref{thm:affine}; step 2}\label{par:pdftaodt} We now study the dynamics of the parabolic 
elements of $G$ near the fixed point $q.$ 

%%%%%%%
\subsubsection{Linear part of automorphisms} Let $g$ be an element of the group $G.$ Let $m$ 
be a fixed point of $g$ and ${\widetilde{m}}$ a point of the fiber $\pi^{-1}(m).$ Let ${\widetilde{g}}_{\widetilde{m}}$
be the unique lift of $g$ to ${\widetilde{S}}$ fixing ${\widetilde{m}}$ (with the notation used in step 1,
${\widetilde{g}}_{\widetilde{q}}={\widetilde{g}}).$ Since $g$ preserves the affine structure, there exists an 
affine transformation $\Af( {\widetilde{g}}_{\widetilde{m}} )$ such that 
$$
\dev \circ {\widetilde{g}}_{\widetilde{m}}= \Af( {\widetilde{g}}_{\widetilde{m}} )\circ \dev.
$$
Note that $\Af( {\widetilde{g}}_{\widetilde{m}} )$ depends on the choice of $m$ and ${\widetilde{m}},$ but that
$\Af({\widetilde{g}}_{\widetilde{m}})$ is uniquely determined by $g$ up to composition by an element
of the monodromy group $\hol(\pi_1(S\setminus \Sing(S)).$ Since the linear parts of the monodromy
are equal to $+{\text{Id}}$ or $-{\text{Id}},$ we get  a well defined morphism 
$$
\left\{ \begin{array}{rcl} G & \to & \PSL(2,\C) \\ g & \mapsto & \Lin(g) \end{array} \right. 
$$
that determines the linear part of $\Af({\widetilde{g}}_{\widetilde{m}})$ modulo $\pm{\text{Id}}$ 
for any choice of $m$ and $\widetilde{m}.$ 

%%%%%%%
\subsubsection{Parabolic elements} Since the linear part $\Lin(g)$ does not depend 
on the fixed point $m,$ it turns out that $\Lin$ preserves the type of $g$: 
We now prove and use this fact in the particular
case of the parabolic elements $g_x,$ $g_y$ and $g_z.$ 

Let $N$ be a positive integer such that $g_x^N$ is contained in $G.$ For $m,$ we 
choose a regular point of $S$ which is periodic of period $l$ for $g_x^N$ and 
which is not a critical point of the projection $\pi_x.$ 
Then $g_x^{Nl}$ fixes the fiber $S_x$ of  $\pi_x$ through $m$ pointwise. 
Since $g_x$ is not periodic and preserves the fibers of $\pi_x,$ this implies that
the differential of $g_x^{Nl}$ at $m$ is parabolic. Let $\widetilde{m}$ be a point of $\pi^{-1}(m)$ and $(\widetilde{g_x^{Nl}})_{\widetilde{m}}$ 
the lift of $g_x^{Nl}$ fixing that point. The universal cover $\pi$ provides a local conjugation between
$g_x^{Nl}$ and $(\widetilde{g_x^{Nl}})_{\widetilde{m}}$ around $m$ and
$\widetilde{m},$ and the developping map provides
a local conjugation between $(\widetilde{g_x^{Nl}})_{\widetilde{m}}$ and $\Lin(g_x^{Nl}).$ As a consequence, $\Lin(g_x^{Nl})$ is a parabolic element of $\PSL(2,\C).$ 

Since a power of $\Lin(g_x^N)$ is parabolic, $\Lin(g_x^N)$ itself is parabolic. 
In particular, the dynamics of ${\widetilde{g_x^N}}$ near ${\widetilde{q}}$ is conjugate to 
a linear upper triangular transformation of $\C^2$ with 
diagonal entries equal to $1.$ 

As a consequence, the lift ${\widetilde{g_x}}$ is locally conjugate near ${\widetilde{q}}$ to 
a linear parabolic transformation with eigenvalues $\pm 1.$ The eigenline of this transformation
corresponds to the fiber $S_z$ through $q.$ Since the local fundamental group $H$ coincides
with $\pm Id,$ this eigenline is mapped to a curve a fixed point by the covering $\pi.$ In particular, 
the fiber $S_z$ through $q$ is a curve of fixed points for $g_x.$

Of course, a similar study holds for $g_y$ and $g_z.$ 

%%%%%%%
\subsubsection{Fixed points and coordinates of the singular point}\label{par:fp_cayley_yo}
The study of fixed points of $g_x,$ $g_y$ and $g_z$ (see lemma \ref{lem:dyna_para})
now shows that the coordinates of the singular point $q$ are equal to $\pm 2.$ Let $\epsilon_x,$ $\epsilon_y$ and $\epsilon_z$ 
be the sign of the coordinates of $q,$ so that
$$
q=(2\epsilon_x, 2\epsilon_y, 2\epsilon_z).
$$
Recall from section \ref{par:sing_fp_orb} that the coefficients $A,$ $B,$ $C,$ 
and $D$ are uniquely determined by the
coordinates of any singular point of $S.$ 
If the product $\epsilon_x \epsilon_y \epsilon_z$ is positive, then, up to symmetry, $q=(2,2,2)$ and $S$ is the surface
$$
x^2 + y^2 + z^2 + xyz= 8x + 8 y + 8 z - 28;
$$
in this case, $q$ is the unique singular point of $S,$ and this singular point is not a node: The second jet of the equation near $q$ is $(x+y+ z)^2=0.$ 
This contradicts the fact that $q$ has to be a node (see section \ref{par:node}). From this we deduce that the 
product $\epsilon_x \epsilon_y \epsilon_z$ is equal to $-1,$ and that $S$ is the Cayley 
cubic.

%%%%%%%%%%%%%%%%%%%%%%%%%%%%%%%%%%%%%%%%%%%%%%%%%%%%%%%%

\section{Irreducibility of Painlev\'e VI Equation.}\label{par:Painleve}

%%%%%%%%%%%%%%%%%%%%%%%%%%%%%%%%%%%%%%%%%%%%%%%%%%%%%%%%%

The goal of this section is to apply the previous section to the irreducibility of Painlev\'e VI equation. 

%%%%%%%
\subsection{Phase space and space of initial conditions}
The naive phase space of Painlev\'e VI equation is parametrized
 by coordinates $(t,q(t),q'(t))\in(\P^1\setminus\{0,1,\infty\})\times\C^2$;
the ``good'' phase space is a convenient semi-compactification
still fibering over the three punctured sphere
$$
\mathcal M(\theta)\to\P^1\setminus\{0,1,\infty\}
$$
whose fibre $\mathcal M_{t_0}(\theta)$, at any point $t_0\in\P^1\setminus\{0,1,\infty\}$,  is the Hirzebruch surface $\mathbb F_2$
 blown-up at $8$-points minus some divisor,
a union of $5$ rational curves (see \cite{Okamoto:1979}).
The analytic type of the fibre, namely the position of the $8$ centers
 and the $5$ rational curves, only depends on Painlev\'e parameters
 $\theta=(\theta_\alpha,\theta_\beta,\theta_\gamma,\theta_\delta)\in\C^4$ and $t_0$.
This fibre bundle is analytically (but not algebraically!)
locally trivial: The local trivialization is given 
by the Painlev\'e foliation (see \cite{SaitoTakebeTerajima})
which is transversal to the fibration.
The monodromy of Painlev\'e equation is given by a representation
$$
\pi_1(\P^1\setminus\{0,1,\infty\},t_0)\to\Diff(\mathcal M_{t_0}(\theta))
$$
into the group of analytic diffeomorphisms of the fibre.

%%%%%%%
\subsection{The Riemann-Hilbert correspondance and $P_{VI}$-monodromy}
On the other hand, the space of initial conditions 
$\mathcal M_{t_0}(\theta)$ may be 
interpreted as the moduli space of rank $2$, trace free 
meromorphic connections over $\P^1$ having simple poles 
at $(p_\alpha,p_\beta,p_\gamma,p_\delta)=(0,t_0,1,\infty)$ 
with prescribed residual 
eigenvalues $\pm\frac{\theta_\alpha}{2}$, $\pm\frac{\theta_\beta}{2}$,
$\pm\frac{\theta_\gamma}{2}$ and $\pm\frac{\theta_\delta}{2}$.
The Riemann-Hilbert correspondance therefore provides an
analytic diffeomorphism 
$$
\mathcal M_{t_0}(\theta)\to \hat{S}_{(A,B,C,D)}
$$
where $\hat{S}_{(A,B,C,D)}$ is the minimal desingularization 
of $S=S_{(A,B,C,D)}$, the parameters $(A,B,C,D)$ being 
given by formulae (\ref{theta/abcd}) 
and (\ref{eq:parameters}).
From this point of view, the Painlev\'e VI foliation coincides with
the isomonodromic foliation: Leaves correspond to universal
 isomonodromic deformations of those connections.
 The monodromy of Painlev\'e VI equation correspond
 to a morphism  
$$
\pi_1(\P^1\setminus\{0,1,\infty\},t_0)\to\Aut(S_{(A,B,C,D)})
$$
and coincides with the $\Gamma_2$-action described in section \ref{par:modularformulae}.
For instance, $g_x$ (resp. $g_y$) is the Painlev\'e VI monodromy
when $t_0$ turns around $0$ (resp. $1$) in the obvious simplest way.
All this is described with much detail in \cite{IIS:2005}.

%%%%%%%
\subsection{Riccati solutions and singular points}
When $S_{(A,B,C,D)}$ is singular, the exceptional divisor
in $\hat{S}_{(A,B,C,D)}$ is a finite union of rational
curves in restriction to which $\Gamma_2$ acts by 
M\"obius transformations. To each such rational curve 
corresponds a rational hypersurface ${\mathcal{H}}$ of the phase space
$\mathcal M(\theta)$ invariant by the Painlev\'e VI foliation.
On ${\mathcal{H}}$, the projection $\mathcal M(\theta)\to\P^1\setminus\{0,1,\infty\}$ restricts to
a regular rational fibration 
and the Painlev\'e equation restricts to a Riccati equation
of hypergeometric type: We get a one parameter family of Riccati solutions.
See \cite{Watanabe:1998,SaitoTerajima,IIS:2005} for a classification
of singular points of $S_{(A,B,C,D)}$ and their link with Riccati solutions;
they occur precisely when either one of the $\theta$-parameter
is an integer, or when the sum $\sum\theta_i$ is an integer.
Since $S_{(A,B,C,D)}$ is affine, there are obviously no other
complete curve in $\mathcal M_{t_0}(\theta)$ 
(see section \ref{par:inv_curves_webs}).

%%%%%%%
\subsection{Algebraic solutions and periodic orbits}
A complete list of algebraic solutions of Painlev\'e VI equation
is still unknown. Apart from those solutions arising as special cases of 
Riccati solutions, that are well known, they correspond
to periodic $\Gamma_2$-orbits on the smooth part of $S_{(A,B,C,D)}$
(see \cite{Iwasaki:2007}). Following section 
\ref{par:subsection_bounded}, 
apart from the three well-known families of $2$, $3$ and $4$-sheeted
algebraic solutions, other algebraic solutions are countable 
and the cosines of the corresponding $\theta$-parameters are real algebraic 
numbers.
In the particular Cayley case $S_C=S_{(0,0,0,4)}$,
periodic $\Gamma_2$-orbits arise from pairs of roots
of unity $(u,v)$ on the two-fold cover $(\C^*)^2$ (see \ref{par:Cayley});
there are infinitely many periodic orbits in this case
and they are dense in the real bounded component of 
$S_C\setminus\{\Sing(S_C)\}$. The corresponding 
algebraic solutions were discovered by Picard in 1889
(before Painlev\'e discovered the general $P_{VI}$-equation !);
see \cite{Mazzocco:2001} and below.
All algebraic solutions (resp. periodic $\Gamma_2$-orbits)
have been classified in the particular case $\theta=(0,0,0,*)$
(resp. $(A,B,C,D)=(0,0,0,*)$) in \cite{Dubrovin-Mazzocco:2000,Mazzocco:2001}: Apart from Riccati and Picard algebraic solutions,
there are $5$ extra solutions up to symmetry (see also \cite{Boalch:2006} for
finite orbit coming from finite subgroups of $\SU(2)$).

Bounded $\Gamma_2$-orbits correspond to what Iwasaki calls
``tame solutions'' in \cite{Iwasaki:2002}.

%%%%%%%
\subsection{Nishioka-Umemura irreducibility}
In 1998, Watanabe
proved in \cite{Watanabe:1998}
the irreducibility of Painlev\'e VI equation in the sense
of Nishioka-Umemura for any 
parameter $\theta$: The generic solution of $P_{VI}(\theta)$
is non classical,
and classical solutions are
\begin{itemize}
\item Riccati solutions (like above),
\item algebraic solutions.
\end{itemize}
Non classical roughly means ``very transcendental'' with regards
to the XIXth century special functions: The general solution 
cannot be expressed in an algebraic way by means of solutions
of linear, or first order non linear differential equations.
A precise definition can be found in \cite{Casale:2007}.

%%%%%%%
\subsection{Malgrange irreducibility}
Another notion of irreducibility was introduced by Malgrange 
in \cite{Malgrange:2001}: He defines the Galois groupoid of
an algebraic foliation to be the smallest algebraic Lie-pseudo-group
that contains the tangent pseudo-group of the foliation (hereafter referred
to as the "pseudo-group");
this may be viewed as a kind of Zariski closure for the pseudo-group
of the foliation. Larger Galois groupoids correspond
to  more complicated  foliations. From this point of view,
it is natural to call irreducible any foliation whose Galois groupoid 
is as large as possible, {\sl{i.e.}} coincides with the basic pseudo-group.

For Painlev\'e equations, a small restriction has to be taken into
account: It has been known since Malmquist that Painlev\'e foliations
may be defined as  kernels of  closed meromorphic  $2$-forms. 
The pseudo-group, and the Galois groupoid, both preserve 
the closed $2$-form. The irreducibility
conjectured by Malgrange is that the Galois groupoid of Painlev\'e
equations coincide with the algebraic Lie-pseudo-group of 
those transformations on the phase space preserving
$\omega$. This was proved for Painlev\'e I equation by Casale
in \cite{Casale:2005}.

For a second order
polynomial differential equation $P(t,y,y',y'')=0$, 
like Painlev\'e equations, Casale
proved in \cite{Casale:2007} that Malgrange-irreducibility
implies Nishioka-Umemura-irreducibility; the converse is not
true as we shall see.

%%%%%%%
\subsection{Invariant geometric structures}
Restricting to a transversal, e.g. the space of initial conditions 
$\mathcal M_{t_0}(\theta)$ for Painlev\'e VI equations,
the Galois grou\-poid defines an algebraic geometric structure
which is invariant under monodromy transformations;
reducibility would imply the existence of an extra geometric structure
on $\mathcal M_{t_0}(\theta)$, additional to the volume form 
$\omega$, preserved by all monodromy transformations.
In that case, a well known result of Cartan, adapted to our 
algebraic setting by Casale in \cite{Casale:2005}, 
asserts that monodromy transformations 
\begin{itemize}
\item either preserve an algebraic foliation,
\item or preserve an algebraic affine structure.
\end{itemize}
Here, ``algebraic'' means that the object is defined 
over an algebraic extension of the field of rational functions,
or equivalently, becomes well-defined over the field of rational functions after some finite ramified 
cover. For instance, ``algebraic foliation'' means polynomial web.
As a corollary of proposition \ref{Prop:NoWeb} 
and Theorem \ref{thm:affine}, we shall prove the following 

\begin{thm}\label{thm:irreducibility}
The sixth Painlev\'e equation is irreducible in the sense of Malgrange, except in one of the following cases:
\begin{itemize}
\item $\theta_\omega\in\frac{1}{2}+\Z$, $\omega=\alpha,\beta,\gamma,\delta$,
\item $\theta_\omega\in \Z$, $\omega=\alpha,\beta,\gamma,\delta$, and $\sum_\omega \theta_\omega$ is even.
\end{itemize}
All these special parameters are equivalent, modulo Okamoto symmetries,
to the case $\theta=(0,0,0,1).$ The corresponding cubic surface
is the Cayley cubic.
\end{thm}

Of course, in the Cayley  case, the existence of an invariant affine structure
shows that the Painlev\'e foliation is Malgrange-reducible
(see \cite{Casale:2006}). This will be made more precise in section \ref{par:picard_reduce}.

Before proving the theorem, we need a stronger version of Lemma
\ref{Lem:NoCurve}

\begin{lem}\label{Lem:NoRiemannSurface}
Let $S$ be an element of the family $\SS.$ 
There is no $\A$-invariant curve of finite type in $S.$
\end{lem}

By "curve of finite type" we mean a complex analytic curve in $S$ with a finite number 
of irreducible components $C_i,$ such that the desingularization of each $C_i$ is a
Riemann surface of finite type. 

\begin{proof}Let $C\subset S$ be a complex analytic curve of  finite type. Since
$S$ is embedded in $\C^3,$ $C$ is not compact. In particular, $C$ is not isomorphic to
the projective line and the group of holomorphic diffeomorphisms of $C$ is virtually solvable.
Since $\A$ contains a non abelian free subgroup, there exists an element $f$ in $\A\setminus \{ { \text{Id}}\}$
which fixes $C$ pointwise. From this we deduce that $C$ is contained in the algebraic curve 
of fixed points of $f.$ This shows that the Zariski closure of $C$ is an $\A$-invariant algebraic 
curve, and we conclude by Lemma \ref{Lem:NoCurve}.
\end{proof}

%%%%%%
\subsection{Proof of theorem \ref{thm:irreducibility}}

In order to prove that Painlev\'e VI equation, for a given parameter
$\theta\in\C^4$ is irreducible, it suffices, due to \cite{Casale:2005}
and the discussion above, to prove that the space of initial
conditions $\mathcal M_{t_0}(\theta)$ does not admit any
monodromy-invariant web or algebraic affine structure.
Via the Riemann-Hilbert correspondance, such a geometric
structure will induce a similar $\Gamma_2$-invariant  structure
on the corresponding character variety $S_{(A,B,C,D)}$.
But we have to be carefull: The Riemann-Hilbert map is not 
algebraic but analytic. As a consequence, the geometric structures 
we have now
to deal with on $S_{(A,B,C,D)}$ are not rational anymore,
but meromorphic (on a finite ramified cover). 
Anyway, the proof of proposition \ref{Prop:NoWeb}
is still valid in this context and exclude the possibility
of $\Gamma_2$-invariant analytic web.

%%%%%%
\subsubsection{Multivalued affine structures.}
We now explain more precisely what is a $\Gamma_2$-invariant
multivalued meromorphic affine structure in the above sense.
First of all, a meromorphic affine structure is an affine structure
in the sense of section \ref{par:InvAffStruct} defined on the
complement of a proper analytic subset $Z$, having moderate growth 
along $Z$ in a sense that we do not need to consider here.
This structure is said to be $\Gamma_2$-invariant if both $Z$, 
and the regular affine structure induced on the complement of $Z$,
are $\Gamma_2$-invariant. Now, a multivalued meromorphic affine structure is a meromorphic structure (with polar locus $Z'$) 
defined on a finite analytic 
ramified cover $\pi':S'\to S$; the ramification locus $X$ is an analytic set. 
This structure is said  to be $\Gamma_2$-invariant if both $X$ 
and $Z=\pi'(Z')$ are invariant and, over the complement of 
$X\cup Z$, $\Gamma_2$ permutes the various regular affine 
structures induced by the various branches of $\pi'$. 

\vv

Let us  prove that the multivalued meromorphic affine structure
induced on $S$ by a reduction of Painlev\'e VI Galois groupoid
has actually no pole, and no ramification apart from singular points
of $S$. Indeed, let $C$ be the union of $Z$ and $R$;
then $C$ is analytic in $S$ but comes from an algebraic curve 
in $\mathcal M_{t_0}(\theta)$ (the initial geometric structure
is algebraic in $\mathcal M_{t_0}(\theta)$), so that the 
$1$-dimensional part of $C$
is a curve of finite type. Lemma 
\ref{Lem:NoRiemannSurface} then show that $C$ is indeed a finite
set. In other words, $C$ is contained in $\Sing(S),$ $R$ itself  is  contained in $\Sing(S)$ and  $Z$ is empty.

%%%%%%
\subsubsection{Singularities of  $S$}

Since the ramification set $R$ is contained in $\Sing(S),$ the cover $\pi'$ is
an \'etale cover in the orbifold sense (singularities of $S'$ are also 
quotient singularities).  Changing the cover $\pi':S'\to S$ if 
necessary, we may assume that $\pi'$ is a Galois cover. 

If $S$ is simply connected, then of course $\pi'$ is trivial, the affine structure is 
univalued, and theorem \ref{thm:affine} provides a contradiction. We can
therefore choose a singularity $q$ of $S,$ and a point $q'$ in
the fiber $(\pi')^{-1}(q).$ Since $\pi_1(S ; q)$ is finitely generated, the
 number of subgroups of index $\deg(\pi')$ in $\pi_1(S ; q)$ is finite.
As a consequence, there is a finite index subgroup $G$ in $\Gamma_2$ 
which lifts to $S'$ and preserves the univalued affine structure defined
on $S'.$ 

We now follow the proof of theorem \ref{thm:affine}  for $G,$ $S'$ and its affine
structure. First, we denote $\pi:{\tilde{S}} \to S'$ the universal cover of $S',$
we choose a point ${\tilde{q}}$ in the fiber $\pi^{-1}(q'),$ and
we lift the action of $G$ to an action on the universal cover ${\tilde{S}}$
fixing ${\tilde{q}}.$ Then we fix a  developping map  $\dev:{\tilde{S}}\to \C^2$ 
with $\dev({\tilde{q}})=0$; these choices imply that $\Af(g)$ is linear for any $g$ in $G.$
Section \ref{par:node} shows that the singularities of $S$ and $S'$ are simple nodes. 

Now comes the main difference with sections \ref{par:lpotm} and \ref{par:pdftaodt}: A priori, the
fundamental group of $S'$ is not generated, as a normal subgroup, by the
local fundamental groups around the singularities of $\Sing(S').$ It could be the case that
$S'$ is smooth, with an infinite fundamental group. So, we need a new argument
 to prove that $g_x$ ({\sl{resp.}} $g_y$ and $g_z$) has a curve of fixed points through 
 the singularity $q.$

%%%%%%
\subsubsection{Parabolic dynamics}
 
 Let $g=g_x^n$ be a non trivial iterate of $g_x$  that is contained in $G.$ The affine transformation
 $\Af(g)$ is linear, with determinant $1$ ; we want to show that this transformation is parabolic.

 Let ${\tilde{U}}$ be an open subset of ${\tilde{S}}$ on which both $\dev$ and the universal cover $\pi'\circ \pi$
 are local diffeomorphisms, and let $U$ be the projection of ${\tilde{U}}$ on $S$ by $\pi'\circ \pi.$  
 We choose ${\tilde{U}}$ in such a way that $U$ contains points $m=(x,y,z)$ with $x$ in the
 interval $[-2,2].$ The fibration of $U$ by fibers of the projection $\pi_x$ is mapped onto a fibration
 ${\mathcal{F}}$ of $\dev({\tilde{U}})$ by the local diffeomorphism $\dev \circ (\pi'\circ\pi)^{-1}.$ Let us prove, first, that  ${\mathcal{F}}$ is a foliation by parallel lines. 
 
 \vv
 
 Let $m$ be a point of $U$ which is $g$-periodic, of period $l.$ Then, the fiber of $\pi_x$ 
 through $m$ is a curve of fixed point for $g^l.$ If ${\tilde{m}}$ is a lift of $m$ in ${\tilde{S}},$
 one can find a lift $\gamma \circ {\tilde{g}}^l$ of $g$ to ${\tilde{S}}$ 
 ($\gamma$ in $\pi_1(S,q)=\Aut(\pi)$)  that fixes pointwise the  fiber through ${\tilde{m}}.$ 
 As a consequence, the fiber of ${\mathcal{F}}$
 through $\dev({\tilde{m}})$ coincides locally with the set of fixed points of 
 the affine transformation $\Af(g^l)\circ \hol(\gamma).$ As such, the fiber of 
 ${\mathcal{F}}$ through $\dev({\tilde{m}})$ is an affine line. 
 
This argument shows that an infinite number of leaves of ${\mathcal{F}}$
are affine lines, or more precisely coincide with the intersection of affine lines with 
$\dev ({\tilde{U}}).$ Since $g$ preserves each fiber of $\pi_x,$  the foliation ${\mathcal{F}}$ is 
leafwise $(\Af(g^l)\circ \hol(\gamma))$-invariant. 
Assume now that $L$ is a  line which coincides with a leaf of ${\mathcal{F}}$ on $\dev ({\tilde{U}}).$
If $L$ is not parallel to the line of fixed points of $\Af(g^l)\circ \hol(\gamma),$ then the affine 
transformation $\Af(g^l)\circ \hol(\gamma)$  is a linear map (since it has a fixed point),
with determinant $\pm 1,$  and with two eigenlines, one of them, the line of fixed points, corresponding to the eigenvalue $1.$  This implies that $\Af(g^l)\circ \hol(\gamma)$
has finiter order.
Since $g$ is not periodic, we conclude that $L$ is parallel to the line of fixed points of $\Af(g^l)\circ \hol(\gamma),$ and that the foliation $ {\mathcal{F}}$ is a foliation 
by parallel lines.

By holomorphic continuation, we get that the image by $\dev$ of the fibration $\pi_x\circ \pi$ is a foliation of the plane by parallel lines. 

\vv

Let us now study the dynamics of ${\tilde{g}}$ near the fixed point ${\tilde{q}}.$ 
Using the local chart $\dev,$ ${\tilde{g}}$ is conjugate to the linear transformation
$\Af(g).$ Since $g$ preserves each fiber of $\pi_x,$ $\Af(g)$ preserves each leaf
of the foliation ${\mathcal{F}}.$ Since $g$ is not periodic, $\Af(g)$ is not periodic
either, and $\Af(g)$ is a linear parabolic transformation. As a consequence, 
$g$ has a curve of fixed points through $q.$ 

\subsubsection{Conclusion}

We can now apply the arguments of section \ref{par:fp_cayley_yo} word by word to conclude that
$S$ is the Cayley cubic.

%We now prove that the cover $\pi:S'\to S$ actually splits into 
%the disjoint union of degree $1$ covers; each component
%$\pi_0:S_0'\to S$ therefore defines a uniform unramified
%and holomorphic affine structure on $S\setminus\{\Sing(S)\}$
%which is invariant by a subgroup $G$ of finite index in $\Gamma_2$
%and we may conclude the proof by Theorem \ref{thm:affine}.
%In order to see this, we consider the action of the fundamental
%group of $S\setminus\{\Sing(S)\}$ on a given fibre $\pi_0^{-1}(m_0)$
%of a given connected component of $\pi$.
%Following Lemma \ref{Lem:TopologyCubic}, the fundamental
%group is normally generated by the local fundamental group
%around the singular points of $S$; it suffices to show that this
%local action is trivial on the fibre $\pi_0^{-1}(m_0)$. 
%The developping map $\dev$ can still be defined form the 
%``irreducible'' multivalued affine structure defined by $\pi_0:S_0'\to S$.
%However, the monodromy is only defined for a subgroup of 
%finite index of the fundamental group of $S$, namely the 
%image of the fundamental group of $S_0'$.
 
 \subsection{Picard parameters of Painlev\'e VI equation and the Cayley cubic}\label{par:picard_reduce}
 
 Let us now explain in more details why the Cayley case is so special 
 with respect to Painlev\'e equations. Consider the universal cover
$$
\pi_t:\C\to\{y^2=x(x-1)(x-t)\}\ ;\ z\mapsto(x(t,z),y(t,z))
$$
of the Legendre elliptic curve
with periods $\Z+\tau\Z$ - this makes sense at least on a
neighborhood of $t_0\in\P^1\setminus\{0,1,\infty\}.$ The
functions  $\tau=\tau(t)$ and $\pi_t$  are analytic in $t$.
 
 The following theorem, obtained by Picard in 1889, shows that
 the Painlev\'e  equations corresponding to the Cayley cubic
 have (almost) classical solutions. 
 
 \begin{thm}[Picard, see \cite{Casale:2006} for example] The general solution of the Painlev\'e sixth differential
equation $P_{VI}(0,0,0,1)$ is given by
$$
t\mapsto x(t,c_1+c_2\cdot\tau(t)),\ \ \ c_1,c_2\in\C.
$$
Moreover, the solution is algebraic if, and only if
$c_1$ and $c_2$ are rational numbers.
\end{thm}

Note that $c_1,c_2\in\Q$ exactly means that 
$\pi_t(c_1+c_2\cdot\tau(t))$ is a torsion point of the elliptic curve.

Finally, $P_{VI}(0,0,0,1)$-equation can actually be integrated
by means of elliptic functions, but in a way that is non classical
with respect to Nishioka-Umemura definition. 
Coming back to 
Malgrange's point of view, the corresponding polynomial affine structure on the phase 
space $\mathcal M(0,0,0,1)$ has been computed by Casale
in \cite{Casale:2006}, thus proving the reducibility of 
$P_{VI}(0,0,0,1)$ equation (and all its birational Okamoto
symmetrics) in the sense of Malgrange.

\section{Appendix A}

This section is devoted to the proof of theorem \ref{thm:carac_cayley}, 
according to which the unique surface in the family $\SS$ 
with four singularities is the Cayley cubic $S_C.$

\begin{proof}

{\bf{I.}} The point $q=(x,y,z)$ is a singular point of $S_{(A,B,C,D)}$ if, and only if
$q$ is contained in $S_{(A,B,C,D)}$ and 
$$
2x+yz=A,\quad 2y+zx=B, \ {\text{ and }} \ 2z+xy=C.
$$
In particular, any pair of  two coordinates of $q$ determines the third coordinate.

\vv

{\bf{II.}} If $(u,v)$ is a couple of complex numbers, 
 $\kappa_{uv}(X)$ will denote  the following quadratic polynomial 
 $$
 \kappa_{uv}(X)= X^2-uvX+(u^2+v^2-4).
 $$
 This polynomial has a double root, namely $\alpha=uv/2,$ if and only if
 $\kappa_{uv}(X)=(X-uv/2)^2,$ if and only if $(u^2-4)(v^2-4)=0.$
 
 Let us now fix a set of $(a,b,c,d)$ parameters that determines $(A,B,C,D).$
 It is proved in \cite{Benedetto-Goldman:1999}  that the coordinates of a singular
 point $q$ satisfy the following properties: 
 \begin{itemize}
 \item[(i)] The $x$ coordinate satisfy one of the following conditions
 \begin{itemize}
 \item $x$ is a double root of $\kappa_{ab}(X),$
 \item $x$ is a double root of $\kappa_{cd}(X),$
 \item $x$ is a  common root of $\kappa_{ab}$ and $\kappa_{cd}(X)$ ;
 \end{itemize}
 \item[(i)] $y$ satisfies the same kind of conditions with respect to $\kappa_{ad}$ and 
 $\kappa_{bc}$;
  \item[(i)] $z$ also, with respect to $\kappa_{ac}$ and $\kappa_{bd}.$
   \end{itemize}

 This shows that the number of possible $x$ (resp. $y,$ $z$)-coordinates for $q$ is 
 bounded from above by $2.$ Together with step {\bf{I}}, this shows that 
 $S_{(A,B,C,D)}$ has at most four singularities. 
 
  When $S_{(A,B,C,D)}$ has four singuarities, there are two possibilities
 for the $x$ coordinate, and either $\kappa_{ab}$ and $\kappa_{cd}$ both have
 a double root, of $\kappa_{ab}$ and $\kappa_{cd}$ coincide and have two simple
 roots. 

 \vv
 
 {\bf{III.}} 
 Let us assume that $\kappa_{ab}$ and $\kappa_{cd}$ have a double root. 
 After a symmetry (see \S \ref{par:twists}), we may assume that $a=c=2.$
 Then, $\kappa_{ac},$ $\kappa_{ad}$ and $\kappa_{bc}$ all have a double root.
 In particular, since $S_{(A,B,C,D)}$ has four singularities, the two choices
 for the $z$-coordinate of singular points are two double roots, the root
 of $\kappa_{ac},$ and, necessarily, the double root of $\kappa_{bd}.$ 
 This implies that $b$ or $d$ is equal to $\pm 2.$ Applying a symmetry of
 the parameters, we may assume that $b=2,$ so that $(a,b,c,d)$ is now
  of type $(2,2,2,d).$
   
   Under this assumption, 
   the $x,$ $y$ and $z$ coordinates of singular points are contained in $\{2,d\}$
   (these are the possible double roots).  If $d^2\neq 4,$ the equations of step {\bf{I}} 
    show that two of the coordinates are equal to $2,$ when 
  one is equal to $d.$  This gives at most three singularities. As a consequence, 
   $d=2$ or $d=-2,$ and the conclusion follows from the fact
    that when $d=2,$ there is only one singularity, namely $(2,2,2).$
    
    \vv
    
    {\bf{IV.}} 
   The last case that we need to consider is when all polynomials
    $\kappa_{uv},$ $u,\, v \in \{a,b,c,d\},$ coincide. In that case, 
    up to symmetries, $a=b=c=d.$ Then, a similar argument shows
    that $a=0$ if $S$ has four singularities (another way to see it
     is to apply the covering $\Quad \circ \Quad$ from section 
     \ref{section:RamifiedCovers}). 
   \end{proof}

\section{Appendix B}\label{par:appendixB}
%%%%%%%%%%%%%%%%%%%%%%%%%%%%%%%%%%%%%%%%

\subsection{Painlev\'e VI parameters $(\theta_\alpha,\theta_\beta,\theta_\gamma,\theta_\delta)$ and Okamoto symmetries}\label{section:Okamoto}

%%%%%%%%%%%%%%%%%%%%%%%%%%%%%%%%%%%%%%%%

Many kinds of conjugacy classes  of representations
$\rho$ with 
$$
\chi(\rho)=(a,b,c,d,x,y,z)
$$ 
give rise to the same 
$(A,B,C,D,x,y,z)$-point ; in order to underline this phenomenon, 
we would like to understand the ramified cover
$$
\Pi\ :\ \ \ 
\left\{\begin{matrix}\C^4&\to&\C^4\\
(a,b,c,d)&\mapsto&(A,B,C,D)
\end{matrix}\right. 
%\ \ \ \text{with}\ \ \ 
%\left\{\begin{matrix}A&=&ab+cd\\ B&=&bc+ad\\ C&=&ac+bd\\
%D&=&4 - a^2 - b^2- c^2 - d^2 - abcd\end{matrix}\right.
$$
defined by equation (\ref{eq:parameters}).

%%%%%%%%%%%%%%%
\subsubsection{Degree of $\Pi$}

\begin{lem}The degree of the covering map $\Pi,$ that is 
the number of points  $(a,b,c,d)$ giving rise to a given generic 
$(A,B,C,D)$-point, is $24$.
\end{lem}

\begin{proof}We firstly assume $B\not=\pm C$ so that $a\not=\pm b$.
Then, solving $B=bc+ad$ and $C=ac+bd$ in $c$ and $d$ yields
$$
c=\frac{aC-bB}{a^2-b^2}\ \ \ \text{and}\ \ \ d=\frac{aB-bC}{a^2-b^2}.
$$
Subsituting in $A=ab+cd$ and $D=4 - a^2 - b^2- c^2 - d^2 - abcd$
gives $\{P=Q=0\}$ with
\begin{eqnarray*}
\, \, P & = & -ab(a^2-b^2)^2+A(a^2-b^2)^2+(B^2+C^2)ab-BC(a^2+b^2) \\
\text{and} \ \,  Q &  = & (a^2+b^2)(a^2-b^2)^2+(D-4)(a^2-b^2)^2 \\
\, & \, & +(B^2+C^2)(a^2-a^2b^2+b^2)+BCab(a^2+b^2-4).
\end{eqnarray*}
These two polynomials have both degree $6$ in $(a,b)$ 
and the corresponding curves must intersect
in $36$ points. However, one easily check that they intersect
along the line at infinity with multiplicity $4$ at each of the two points
 $(a:b)=(1:1)$ and $(1:-1)$; moreover, they also intersect
along the forbidden lines $a=\pm b$ at $(a,b)=(0,0)$ with multiplicity
$4$ as well, provided that $BC\not=0$. As a consequence, the number of preimages of $(A,B,C,D)$ is $36-4-4-4=24$ (counted with multiplicity).
 \end{proof}

\begin{rem}$\Pi$ is not a Galois cover: The group of deck transformations
 is the order $8$ group $Q=\langle P_1,P_2,\otimes_{(-1,-1,-1,-1)}\rangle$
(see \S \ref{par:twists}). \end{rem}

%%%%%%%%%%
\subsubsection{Okamoto symmetries}
To understand the previous remark, it is convenient to introduce the 
Painlev\'e VI parameters, which are related to $(a,b,c,d)$ by the map
$$
\begin{matrix}\C^4&\to&\C^4\\
(\theta_\alpha,\theta_\beta,\theta_\gamma,\theta_\delta)
&\mapsto&(a,b,c,d)\end{matrix}
\ \ \ \text{with}\ \ \ 
\left\{\begin{matrix}
a & = & 2\cos(\pi\theta_\alpha)\\
b & = & 2\cos(\pi\theta_\beta)\\
c & = & 2\cos(\pi\theta_\gamma)\\
d & = & 2\cos(\pi\theta_\delta)
\end{matrix}\right.
$$
The composite map $(\theta_\alpha,\theta_\beta,\theta_\gamma,\theta_\delta)\mapsto(A,B,C,D)$ has been studied in \cite{IIS:2004}:
It is an infinite Galois ramified cover whose deck transformations 
coincide with the group $G$ of so called Okamoto symmetries. Those symmetries 
are "birational transformations" of Painlev\'e VI equation; they have
been computed directly on the equation by Okamoto in \cite{Okamoto:1987} 
(see \cite{Noumi-Yamada} for a modern presentation). 
Let $\Bir(P_{VI})$ be the group of all birational symmetries of Painlev\'e sixth equation.
The Galois group $G$ is the subgroup of $\Bir(P_{VI})$ generated 
by the following four kind of affine transformations. 

\begin{enumerate}
\item Even 
translations by integers
$$
\oplus_n\ :\ \ \ 
\left\{\begin{matrix}
\theta_\alpha&\mapsto&\theta_\alpha+n_1\\
\theta_\beta&\mapsto&\theta_\beta+n_2\\
\theta_\gamma&\mapsto&\theta_\gamma+n_3\\
\theta_\delta&\mapsto&\theta_\delta+n_4
\end{matrix}\right.
\ \ \ \text{with}\ \ \ 
\left\{\begin{matrix}n=(n_1,n_2,n_3,n_4)\in\Z^4,\\ 
n_1+n_2+n_3+n_4\in2\Z.
\end{matrix}\right.
$$
Those symmetries also act on the space of initial conditions
of $P_{VI}$ in a non trivial way, but the corresponding action 
on $(x,y,z)$ is very simple: We recover the twist symmetries
$\otimes_\epsilon$ of section \ref{par:twists} by considering 
$n$ modulo $2\Z^4$.

\item An action of $\Sym_4$ permuting 
$(\theta_\alpha,\theta_\beta,\theta_\gamma,\theta_\delta).$ This
corresponds to the action of $\Sym_4$ on $(a,b,c,d,x,y,z)$
permuting $(a,b,c,d)$ in the same way. This group is generated by 
the four permutations $T_1$, $T_2$, $P_1$ and $P_2$ (see sections
\ref{par:toruscover} and \ref{par:modularformulae}).

\item Twist symmetries on Painlev\'e parameters
$$
\otimes_\epsilon\ :\ \ \ 
\left\{\begin{matrix}
\theta_\alpha&\mapsto&\epsilon_1\theta_\alpha\\
\theta_\beta&\mapsto&\epsilon_2\theta_\beta\\
\theta_\gamma&\mapsto&\epsilon_3\theta_\gamma\\
\theta_\delta&\mapsto&\epsilon_4\theta_\delta
\end{matrix}\right.
\ \ \ \text{with}\ \ \ 
\epsilon=(\epsilon_1,\epsilon_2,\epsilon_3,\epsilon_4)\in\{\pm1\}^4.
$$
The corresponding action on $(a,b,c,d,x,y,z)$ is trivial.

\item The special Okamoto symmetry 
(called $s_2$ in  \cite{Noumi-Yamada}) 
$$
\Ok\ :\ \ \ 
\left\{\begin{matrix}
\theta_\alpha&\mapsto&\frac{\theta_\alpha-\theta_\beta-\theta_\gamma-\theta_\delta}{2}+1\\
\theta_\beta&\mapsto&\frac{-\theta_\alpha+\theta_\beta-\theta_\gamma-\theta_\delta}{2}+1\\
\theta_\gamma&\mapsto&\frac{-\theta_\alpha-\theta_\beta+\theta_\gamma-\theta_\delta}{2}+1\\
\theta_\delta&\mapsto&\frac{-\theta_\alpha-\theta_\beta-\theta_\gamma+\theta_\delta}{2}+1
\end{matrix}\right.
$$
The corresponding action on $(A,B,C,D,x,y,z)$ is trivial (see \cite{IIS:2004}), but the action on $(a,b,c,d)$ is rather subbtle, as we shall see.
\end{enumerate}

The ramified cover 
$(\theta_\alpha,\theta_\beta,\theta_\gamma,\theta_\delta)\mapsto
(a,b,c,d)
$
is also a Galois cover: Its Galois group $K$ is the subgroup of $G$ 
generated by those translations $\oplus_n$ with $n\in(2\Z)^4$ 
and the twists $\otimes_\epsilon$. 
%The ramified cover 
%$$(\theta_\alpha,\theta_\beta,\theta_\gamma,\theta_\delta)\mapsto
%(A,B,C,D)
%$$
%is Galois as well, with Galois group
%$G\subset\Bir(P_{VI})$ generated by $K$, $\oplus_{(-1,-1,-1,-1)}$, 
%the permutations
%$(\theta_\delta,\theta_\gamma,\theta_\beta,\theta_\alpha)$
%and $(\theta_\beta,\theta_\alpha,\theta_\delta,\theta_\gamma)$
%(corresponding to $P_1$ and $P_2$) and $\Ok$. 
One can check that $[G:K]=24$ but $K$ is not a normal subgroup 
of $G$: It is not $\Ok$-invariant. In fact, $K$ is normal in the subgroup $G'\subset G$ where we omit the generator $\Ok$ and $Q=G'/K$
coincides with the order $8$ group of symmetries fixing $(A,B,C,D)$.
Therefore, $G/K$ may be viewed as the disjoint union of left cosets
$$
G/K=Q\cup \Ok\cdot Q\cup \widetilde{\Ok}\cdot Q
$$
where $\widetilde{\Ok}$ is the following symmetry (called $s_1s_2s_1$
 in  \cite{Noumi-Yamada}) 
$$
\widetilde{\Ok}\ :\ \ \ 
\left\{\begin{matrix}
\theta_\alpha&\mapsto&\frac{\theta_\alpha-\theta_\beta-\theta_\gamma+\theta_\delta}{2}\\
\theta_\beta&\mapsto&\frac{-\theta_\alpha+\theta_\beta-\theta_\gamma+\theta_\delta}{2}\\
\theta_\gamma&\mapsto&\frac{-\theta_\alpha-\theta_\beta+\theta_\gamma+\theta_\delta}{2}\\
\theta_\delta&\mapsto&\frac{\theta_\alpha+\theta_\beta+\theta_\gamma+\theta_\delta}{2}
\end{matrix}\right.
$$

%%%%%%%%%%%%%%%%%
\subsubsection{From $(A,B,C,D)$ to $(a,b,c,d)$} Now, given a $(a,b,c,d)$-point, we 
would like to describe explicitly
all other parameters $(a',b',c',d')$ in the same $\Pi$-fibre,
i.e. giving rise to the same parameter $(A,B,C,D)$. 
We already know that the $Q$-orbit
$$
\left\{\begin{matrix}
(a,b,c,d)&(-a,-b,-c,-d)&(d,c,b,a)&(-d,-c,-b,-a)\\
(b,a,d,c)&(-b,-a,-d,-c)&(c,d,a,b)&(-c,-d,-a,-b)
\end{matrix}\right\},
$$
, which generically is of length $8,$ is contained in the fibre.
In order to describe the remaining part of the fibre, 
let us choose $(a_\epsilon,b_\epsilon,c_\epsilon,d_\epsilon)\in\C^4$,
$\epsilon=0,1$,
such that 
$$
\left\{\begin{matrix}
a_0&=&\frac{\sqrt{2+a}}{2}\\
b_0&=&\frac{\sqrt{2+b}}{2}\\
c_0&=&\frac{\sqrt{2+c}}{2}\\
d_0&=&\frac{\sqrt{2+d}}{2}
\end{matrix}\right.
\ \ \ \text{and}\ \ \ 
\left\{\begin{matrix}
a_1&=&\frac{\sqrt{2-a}}{2}\\
b_1&=&\frac{\sqrt{2-b}}{2}\\
c_1&=&\frac{\sqrt{2-c}}{2}\\
d_1&=&\frac{\sqrt{2-d}}{2}
\end{matrix}\right.
$$
If $\theta_\alpha$ is such that $(a_0,a_1)=
(\cos(\pi\frac{\theta_\alpha}{2}),\sin(\pi\frac{\theta_\alpha}{2}))$, 
then $a=2\cos(\pi\theta_\alpha)$; therefore, the choice of 
$(a_0,a_1)$ is equivalent to the choice of a $P_{VI}$-parameter
$\theta_\alpha$ modulo $2\Z$, i.e. of $\frac{\theta_\alpha}{2}$
modulo $\Z$. Then, looking at the action of the special
Okamoto symmetry $\Ok$ on Painlev\'e parameters 
$(\theta_\alpha,\theta_\beta,\theta_\gamma,\theta_\delta)$,
we derive the following new point $(a',b',c',d')$ in the $\Pi$-fibre 
$$
\left\{\begin{matrix}
a'&=&-2\sum(-1)^{\frac{\sum\epsilon_i}{2}+\epsilon_1}a_{\epsilon_1}b_{\epsilon_2}c_{\epsilon_3}d_{\epsilon_4}\\
b'&=&-2\sum(-1)^{\frac{\sum\epsilon_i}{2}+\epsilon_2}a_{\epsilon_1}b_{\epsilon_2}c_{\epsilon_3}d_{\epsilon_4}\\
c'&=&-2\sum(-1)^{\frac{\sum\epsilon_i}{2}+\epsilon_3}a_{\epsilon_1}b_{\epsilon_2}c_{\epsilon_3}d_{\epsilon_4}\\
d'&=&-2\sum(-1)^{\frac{\sum\epsilon_i}{2}+\epsilon_4}a_{\epsilon_1}b_{\epsilon_2}c_{\epsilon_3}d_{\epsilon_4}
\end{matrix}\right.
$$
where the sum is taken over all $\epsilon=(\epsilon_1,\epsilon_2,\epsilon_3,\epsilon_4)\in(\{0,1\})^4$ for which $\sum_{i=1}^4\epsilon_i$
is even. One can check that the different choices for 
$(a_0,b_0,c_0,d_0)$ and $(a_1,b_1,c_1,d_1)$
lead to $16$ distinct possible $(a',b',c',d')$, namely
$2$ distinct $Q$-orbits, which together with the $Q$-orbit
of $(a,b,c,d)$ above provide the whole $\Pi$-fibre.

\begin{eg}\label{eg:OkamotoSymetrics}
When $(a,b,c,d)=(0,0,0,d)$, we have $(A,B,C,D)=(0,0,0,D)$
with $D=4-d^2$. The $\Pi$-fibre is given by the $Q$-orbits 
of the $3$ points
$$
(0,0,0,d)\ \ \ \text{and}\ \ \ (\tilde d,\tilde d,\tilde d,-\tilde d)\ \ \ 
\text{where}\ \ \ \tilde d=\sqrt{2\pm\sqrt{4-d^2}}
$$
(only the sign of the square root inside is relevant up to $Q$).
The fibre has length $24$ except in the Cayley case $d=0$
where it has length $9$, consisting of the two $Q$-orbits of 
$$
(0,0,0,0)\ \ \ \text{and}\ \ \ (2,2,2,-2)
$$
(note that $(0,0,0,0)$ is $Q$-invariant) and in the Markov case 
$d=2$ where it has length $16$, consisting of the two $Q$-orbits of 
$$
(0,0,0,2)\ \ \ \text{and}\ \ \ (\sqrt{2},\sqrt{2},\sqrt{2},-\sqrt{2}).
$$
\end{eg}

%%%%%%%%%%%%%%%%%%%%%%%%%%%%%%%%%

\subsection{Reducible representations versus singularities.}\label{section:ReducibleSingular}

%%%%%%%%%%%%%%%%%%%%%%%%%%%%%%%%%%

\begin{thm}[\cite{Benedetto-Goldman:1999,IIS:2005}]\label{Thm:ReducibleSingular}
The surface
$S_{(A,B,C,D)}$ is singular if, and only if, we are in one of the following cases
\begin{itemize}
\item $\Delta(a,b,c,d)=0$ where
$$
\Delta=(2(a^2+b^2+c^2+d^2)-abcd-16)^2-(4-a^2)(4-b^2)(4-c^2)(4-d^2),
$$
\item at least one of the parameters $a$, $b$, $c$ or $d$ 
equals $\pm2$.
\end{itemize}
More precisely, a representation $\rho$ is sent to a singular point if, and only if, we are in one of the following cases :
\begin{itemize}
\item the representation $\rho$ is reducible and then $\Delta=0$,
\item one of the generators $\rho(\alpha)$, $\rho(\beta)$, $\rho(\gamma)$ or $\rho(\delta)$ equals $\pm I$ (the corresponding trace parameter is then equal to $\pm2$).
\end{itemize}
\end{thm}

In fact, it is proved in \cite{Benedetto-Goldman:1999} 
that the set $Z$ of parameters 
$(A,B,C,D)$ for which $S_{(A,B,C,D)}$ is singular is defined by 
$\delta=0$ where $\delta$ is the discriminant of the 
polynomial 
$$P_z=z^4 -C z^3 -(D+4)z^2 + (4C-AB)z +4D + A^2 + B^2$$ 
defined in section \ref{par:dynamics_parabolic}: $P_z$ 
has a multiple root at each singular point.
Now, consider the ramified cover 
$$\Pi:\C^4\to\C^4;(a,b,c,d)\mapsto(A,B,C,D)$$
defined by (\ref{eq:parameters}).
One can check by direct computation that 
$$
\delta\circ\Pi=\frac{1}{16}(a^2-4)(b^2-4)(c^2-4)(d^2-4)\Delta^2.
$$
One also easily verifies that the locus of reducible representations
is also the ramification locus of $\Pi$:
$$\mathrm{Jac}(\Pi)=-\frac{1}{2}\Delta.$$

It is a well known fact (see \cite{IIS:2005}) that Okamoto symmetries
permute the two kinds of degenerate representations 
given by Theorem \ref{Thm:ReducibleSingular}.
For instance,
a singular point is defined by the following equations:
$$
A=2x+yz,\ \ \ B=2y+xz,\ \ \ C=2z+xy
$$
$$
\text{and}\ \ \ x^2+y^2+z^2+xyz=Ax+By+Cz+D.
$$
Now, a compatible choice of parameters $(a,b,c,d)$ is provided by
$$
(a,b,c,d)=(y,z,x,2)
$$
and one easily check that the corresponding representations
satisfy $\rho(\delta)=I$.

%%%%%%%%%%%%%%%%%%%%%%%%%%%%%%%%%

\subsection{$\SU(2)$-representations versus bounded components.}\label{section:SU(2)}

%%%%%%%%%%%%%%%%%%%%%%%%%%%%%%%%%%

When $a$, $b$, $c$, and $d$ are real numbers, $A$, $B$, $C$, and $D$ are real as well. In that case, the real part $S_{(A,B,C,D)}(\R)$ stands for $\SU(2)$ and $\SL(2,\R)$-representations; precisely, each connected component of 
the smooth part of $S_{(A,B,C,D)}(\R)$ is either purely $\SU(2)$, 
or purely $\SL(2,\R)$, depending
on the choice of $(a,b,c,d)$ fitting to $(A,B,C,D).$

Moving into the parameter space $\{(a,b,c,d)\},$ when we pass from 
$\SU(2)$ to $\SL(2,\R)$-representations,
we must go through a representation of the group 
$\SU(2)\cap\SL(2,\R)=\SO(2,\R).$ Since representations into $\SO(2,\R)$ are reducible, they  correspond to singular points of the cubic surface 
(see \S \ref{section:ReducibleSingular}). In other words, 
any bifurcation between $\SU(2)$ and $\SL(2,\R)$-represen\-tations creates a real singular point of $S_{(A,B,C,D)}$.

Since $\SU(2)$-representations are contained in 
the cube $[-2,2]^3$, they always form a bounded 
component of the smooth part of $S_{(A,B,C,D)}(\R)$:
Unbounded components always correspond to 
$\SL(2,\R)$-representations, whatever the choice 
of parameters $(a,b,c,d)$ is.

The topology of $S_{(A,B,C,D)}(\R)$ is studied in
\cite{Benedetto-Goldman:1999} when $(a,b,c,d)$ are real numbers. 
There are at most four singular points, 
and the smooth part has at most one bounded 
and at most four unbounded components. 
On the other hand, if $A$, $B$, $C$, and $D$ are  real numbers, then
$a$, $b$, $c$, and $d$ are not necessarily real. 
 
 \begin{eg} If $a$, $b$, $c$, and $d$  are purely imaginary numbers, then
$A$, $B$, $C$, and $D$ are  real numbers. In this specific example, there are representations $\rho:\pi_1(\Sphere_4)\to \SL(2,\C)$ with trace parameters 
$$
(a,b,c,d,x,y,z)\in (i\R)^4\times (\R)^3,
$$ 
the image of which are Zariski dense in the (real) Lie group $\SL(2,\C).$ Such 
a representation correspond to a point $(x,y,z)$ on $\S_{(A,B,C,D)(\R)}$ which is
not realized by a representation into $\SL(2,\R).$
\end{eg}

The goal of this section is to prove the following theorem, which partly extends
the above mentionned results of Benedetto and Goldman \cite{Benedetto-Goldman:1999}.

\begin{thm}\label{Thm:boundedSU2}
Let $A$, $B$, $C$, and $D$ be real numbers, for 
which the smooth part of $S_{(A,B,C,D)}(\R)$ 
has a bounded component.
Then for any choice of parameters $(a,b,c,d)$ fitting to $(A,B,C,D)$,
the numbers $a$, $b$, $c$, and $d$ are real, contained in $(-2,2)$
and the bounded component stands for $\SU(2)$ or
$\SL(2,\R)$-representations. 
Moreover, for any such parameter $(A,B,C,D)$, we can choose 
between $\SU(2)$ and $\SL(2,\R)$ by conveniently choosing
$(a,b,c,d)$: The two cases both occur.
\end{thm}

In particular, bounded components of real surfaces 
$S_{(A,B,C,D)}(\R)$ always arise from $\SU(2)$-representations\footnote{This strengthens the results of \cite{Previte-Xia:2005}
where the bounded component was assumed to arise from $\SU(2)$-representations.}.

Denote by $Z\subset\R^4$ the subset of those parameters 
$(A,B,C,D)$ for which the corresponding surface $S_{(A,B,C,D)}(\R)$
is singular (see section \ref{section:ReducibleSingular}). 
Over each connected component of $\R^4\setminus Z$,
the surface $S_{(A,B,C,D)}(\R)$ is smooth and has constant
topological type. Let $\mathcal B$ be the union of connected components of $\R^4\setminus Z$ over which the smooth surface 
has a bounded component. 

The ramified cover $\Pi:\C^4\to\C^4;(a,b,c,d)\mapsto(A,B,C,D)$
has degree $24$; Okamoto correspondences, defined
in section \ref{section:Okamoto}, ``act'' transitively on fibers
(recall that $\Pi$ is not Galois). 
Because of their real nature, these correspondences permute
real parameters $(a,b,c,d)$: Therefore, $\Pi$ restricts as a degree $24$
ramified cover $\Pi\vert_{\R^4}:\R^4\to\Pi(\R^4)$. 
Following \cite{Benedetto-Goldman:1999}, we have
$$
\Pi^{-1}(\mathcal B)\cap\mathbb R^4=(-2,2)^4\setminus\{\Delta=0\}.
$$
Using again that $\SU(2)\cap\SL(2,\R)=\mathrm{SO}(2)$ is abelian,
and therefore corresponds to reducible representations,
we promptly deduce that, along each connected component of 
$(-2,2)^4\setminus\{\Delta=0\}$, the bounded component of the
corresponding surface $S_{(A,B,C,D)}(\R)$ constantly stands 
either for $\SU(2)$-representations, or for $\SL(2,\R)$-representations.
We shall denote by $\mathcal B^{\SU(2)}$ and $\mathcal B^{\SL(2,\R)}$
the corresponding components of $\mathcal B$. 
Theorem \ref{Thm:boundedSU2} may now be rephrased as the following
equalities:
$$
\mathcal B=\mathcal B^{\SU(2)}=\mathcal B^{\SL(2,\R)}.
$$
To prove these equalities, we first note that
t $\mathcal B^{\SU(2)}\cup\mathcal B^{\SL(2,\R)}\subset\Pi([-2,2]^4)$ is obviously bounded by $-8\le A,B,C\le 8$ and $-20\le D\le 28$ (this bound is not sharp !).

\begin{lem}The set $\mathcal B$ is bounded, contained into
$-8\le A,B,C\le 8$ and $-56\le D\le 68$.
\end{lem}

\begin{proof}The orbit of any point $p$ belonging to a bounded
component of $S_{(A,B,C,D)}(\R)$ is bounded. Applying 
the tools involved in section \ref{par:sectionbounded},
we deduce that the bounded component is contained into 
$[-2,2]^3$. Therefore, for any $p=(x,y,z)$ and $s_x(p)=(x',y,z)$ 
belonging to the bounded component,
we get $A=x+x'+yz$ and then $-8\le A\le 8$.
Using $s_y$ and $s_z$, we get the same bounds for $B$ and $C$.
Since $p$ is in the surface, we also get $D=x^2+y^2+z^2+xyz-Ax-By-Cz$.
\end{proof}

The order $24$ group of Benedetto-Goldman symmetries act on the parameters $(A,B,C,D)$ by freely permutting the triple $(A,B,C)$,
and freely changing sign for two of them. This group acts on 
the set of connected components of $\R^4\setminus Z$, 
$\mathcal B$, $\mathcal B^{\SU(2)}$ and $\mathcal B^{\SL(2,\R)}$.
The crucial Lemma is

\begin{lem}Up to Benedetto-Goldman symmetries,
$\R^4\setminus Z$ has only one bounded component.
\end{lem}

\begin{proof}Up to Benedetto-Goldman symmetries,
one can always assume $0\le A\le B\le C$.
This fact is easily checked by looking at the action 
of symmetries on the projective coordinates
$[A:B:C]=[X:Y:1]$: the triangle $T=\{0\le X\le Y\le 1\}$
happens to be a fundamental domain for this group action. 
We shall show that 
$\R^4\setminus Z$ has at most one bounded component 
over the cone
$$\mathcal C=\{(A,B,C)\ ;\ 0\le A\le B\le C\}$$
with respect to the projection $(A,B,C,D)\mapsto(A,B,C)$.

The discriminant of $\delta$ with respect to $D$ reads
$$
\mathrm{disc}(\delta)=
-65536\, \left( B-C \right) ^{2} \left( B+C \right) ^{2} \left( A-C
 \right) ^{2} \left( A+C \right) ^{2} \left( A-B \right) ^{2} \left( A
+B \right) ^{2}\delta_1^3
$$
where $\delta_1$ is the following polynomial (with $(X,Y)=(\frac{A}{C},\frac{B}{C})$)
$$
\delta_1=-{C}^{9}{X}^{3}{Y}^{3}+ \left( 27\,{Y}^{4}+27\,{X}^{4}{Y}^{4}
-6\,{X}^{2}{Y}^{4}-6\,{X}^{4}{Y}^{2}+27\,{X}^{4}-6\,{X}^{2}{Y}^{2}
 \right) {C}^{8}
 $$
 $$
 + \left( -768\,{X}^{5}Y+192\,{Y}^{3}X-768\,XY+192\,{X}
^{3}Y-768\,{Y}^{5}X+192\,{X}^{3}{Y}^{3} \right) {C}^{7}
 $$
 $$
+ \left( 4096\,
{Y}^{6}-1536\,{Y}^{2}+4096+23808\,{X}^{2}{Y}^{2}-1536\,{X}^{4}-1536\,{
X}^{2}{Y}^{4}
\right. $$
 $$\left.
 -1536\,{X}^{4}{Y}^{2}+4096\,{X}^{6}-1536\,{X}^{2}-1536\,{
Y}^{4} \right) {C}^{6}
 $$
 $$
+ \left( -86016\,{X}^{3}Y-86016\,{Y}^{3}X-86016
\,XY \right) {C}^{5}
 + \left( 712704\,{X}^{2}{Y}^{2}
 \right. $$
 $$\left.
 -196608\,{Y}^{4}-
196608-196608\,{X}^{4}+712704\,{X}^{2}+712704\,{Y}^{2} \right) {C}^{4}
 $$
 $$
-5505024\,{C}^{3}XY+ \left( 3145728\,{X}^{2}+3145728+3145728\,{Y}^{2}
 \right) {C}^{2}
 -16777216
$$
First, we want to show that $\mathcal C\setminus\{\mathrm{disc}(\delta)=0\}$ 
has $5$ connected components, only two of which are bounded.
The polynomial $\delta_1$ has degree $9$ in $C$ in restriction 
to any line $L_{X,Y}=\{A=XC,B=YC\}\subset\mathcal C$ 
with $0<X<Y<1$; we claim that it has constantly $3$ simple
real roots (and $6$ non real ones) 
$$c_1(X,Y)<0<c_2(X,Y)<c_3(X,Y).$$
In order to check this, let us verify that the discriminant of $\delta_1$
with respect to $C$ does not vanish in the interior of the triangle $T$.
After computations, we find
$$
\mathrm{disc}(\delta_1)=k(X^2-Y^2)^8(X^2-1)^8(Y^2-1)^8(Y\delta_2)^2
$$
where $k$ is a huge constant and $\delta_2$ is given, setting
$X=tY$, by
$$
\delta_2=\left( 22272\,{t}^{8}+40337\,{t}^{6}+16384\,{t}^{10}+16384\,{t}^{2}+
22272\,{t}^{4} \right) {Y}^{10}
$$
$$
+ \left( -59233\,{t}^{4}+16384\,{t}^{10
}-59233\,{t}^{6}+40448\,{t}^{8}+16384+40448\,{t}^{2} \right) {Y}^{8}
$$
$$
+
 \left( 22272+22272\,{t}^{8}-59233\,{t}^{2}-59233\,{t}^{6}-118893\,{t}
^{4} \right) {Y}^{6}
$$
$$
+ \left( 40337\,{t}^{6}+40337-59233\,{t}^{2}-59233
\,{t}^{4} \right) {Y}^{4}
$$
$$
+\left( 22272+22272\,{t}^{4}+40448\,{t}^{2}
 \right) {Y}^{2}+16384+16384\,{t}^{2}.
$$
This later polynomial has non vanishing discriminant 
with respect to $Y$ for $0<t<1$ and has a non real root, for instance,
when $t=1/2$: Thus $\mathrm{disc}(\delta_1)$ does not vanish
in the interior of the triangle $T$. Therefore, the polynomial
$\delta_1$ has always the same number of real roots
when $(X,Y)$ lie inside the triangle $T$ and one can 
easily check that $0$ is never a root, and by specializing $(X,Y)$, 
that there are indeed $3$ roots, one of them being negative.
The claim is proved.

The cone $\mathcal C$ is cutted off by $\mathrm{disc}(\delta)=0$ into 
$5$ components, namely
$$
\mathcal C_1=\{C<c_1(X,Y)\},\ \ \ \mathcal C_2=\{c_1(X,Y)<C<0\},\ \ \ \mathcal C_3=\{0<C<c_2(X,Y)\},
$$
$$
\mathcal C_4=\{c_2(X,Y)<C<c_3(X,Y)\}\ \ \ \text{and}\ \ \ \mathcal C_5=\{c_3(X,Y)<C\}.
$$
But $\delta_1$ has degree $8$ when $X=0$ and one of the roots
$c_i(X,Y)$ tends to infinity when $X\to0$. One can check
that $c_3\to\infty$ and only $\mathcal C_2$ and $\mathcal C_3$
are bounded. 

We now study the possible bounded components of 
$\R^4\setminus Z$ over the cona $\mathcal C$;
they necessarily project onto 
$\mathcal C_2$, $\mathcal C_3$ or the union 
(together with $(A,B,C)=0$).
The polynomial $\delta$ defining $Z$ has degree $5$ in $D$. 
After several numerical specializations, we obtain the following
picture:
\begin{itemize}
\item the polynomial $\delta$ has $5$ real roots 
$d_1<d_2<d_3<d_4<d_5$ over $\mathcal C_2$ and $\mathcal C_3$,
$d_i=d_i(A,B,C)$ for $i=1,\ldots,5$, 
\item over $C=c_1$ or $C=c_2$, $0<A<B<C$, the $5$ roots extend continuously, satisfying $d_1=d_2<d_3<d_4<d_5$
\item over $(A,B,C)=0$, the $5$ roots extend continuously as $0=d_1<d_2=d_3=d_4=4$.
\end{itemize}
Among the $6$ connected components of $\R^4\setminus Z$
over $\mathcal C_2$ (resp. $\mathcal C_3$), only that one 
defined by $\{d_1(A,B,C)<D<d_2(A,B,C)\}$ does not 
``extend'' over the unbounded component 
$\mathcal C_1$ (resp. $\mathcal C_4$).
The unique bounded component of $\R^4\setminus Z$ 
over the cona $\mathcal C$ is therefore defined over
$\mathcal C_2\cup\{A=B=C=0\}\cup\mathcal C_3$
by $\{d_1(A,B,C)<D<d_2(A,B,C)\}$. The corresponding
connected component of $\R^4\setminus Z$ must be
bounded as well, since there is at least one bounded 
component, given by $\mathcal B^{\SU(2)}$, 
or $\mathcal B^{\SL(2,\R}$.
\end{proof}

We thus conclude that $\mathcal B=\mathcal B^{\SU(2)}=\mathcal B^{\SL(2,\R}$ and Theorem \ref{Thm:boundedSU2} is proved 
in the case the real surface $S_{(A,B,C,D)}(\R)$ is smooth.
The general case follows from the following lemma, the proof of which
is left to the reader.

\begin{lem}Let $(A,B,C,D)$ be real parameters such that
the smooth part of the surface $S_{(A,B,C,D)}(\R)$ 
has a bounded component. Then, there exist an arbitrary 
small real perturbation of $(A,B,C,D)$ such that the corresponding surface
is smooth and has a bounded component.
\end{lem}

We would like now to show that there is actually
only one bounded component in $\R^4\setminus Z$
(up to nothing).

Inside $[-2,2]^4$, the equation $\Delta$ splits into
the following two equations
$$
2(a^2+b^2+c^2+d^2)-abcd-16=\pm\sqrt{(4-a^2)(4-b^2)(4-c^2)(4-d^2)}.
$$
Those two equations cut-off the parameter space $[-2,2]^4$
into many connected components and we have\footnote{
In \cite{Benedetto-Goldman:1999}, the connected components 
of $[-2,2]^4$ standing for $\SL(2,\R)$-representations are 
equivalently defined by $\Delta>0$ and $2(a^2+b^2+c^2+d^2)-abcd-16>0$.}

\begin{thm}[Benedetto-Goldman \cite{Benedetto-Goldman:1999}]\label{Thm:Benedetto-Goldman}
When $a$, $b$, $c$ and $d$ are real and $S_{(A,B,C,D)}(\R)$
is smooth, then $S_{(A,B,C,D)}(\R)$ has a bounded component
if, and only if, $a$, $b$, $c$ and $d$ both lie in $(-2,2)$.
In this case, the bounded component 
corresponds to $\SL(2,\R)$-representations if, and only if, 
$$
2(a^2+b^2+c^2+d^2)-abcd-16>\sqrt{(4-a^2)(4-b^2)(4-c^2)(4-d^2)}.
$$
\end{thm}

When we cross the boundary 
$$
2(a^2+b^2+c^2+d^2)-abcd-16=\sqrt{(4-a^2)(4-b^2)(4-c^2)(4-d^2)}
$$ 
inside $(-2,2)^4$, we pass from 
$\SL(2,\R)$ to $\SU(2)$-representations: At the boundary, 
the bounded component must degenerate down to a singular point.

We now prove the 

\begin{pro}\label{Prop:SU(2)SL(2)}
The set $(-2,2)^4\setminus\{\Delta=0\}$ has $24$ connected
components, $8$ of them corresponding to $\SL(2,\R)$-representations.
Okamoto correspondence permute transitively
those components.
\end{pro} 

Recall that the group of cover transformations $Q$
has order $8$ and does not change the nature of
the representation: The image $\rho(\pi_1(\Sphere_4))$
remains unchanged in $\PGL(2,\mathbb C)$.
Therefore, up to this tame action, Okamoto correspondence
provides, to any smooth point $(A,B,C,D,x,y,z)$ of the character
variety, exactly $3$ essentially distinct representations,
two of them in $\SU(2)$, and the third one in $\SL(2,\R)$.
It may happens (see \cite{Previte-Xia:2003}) that one 
of the two $\SU(2)$-representations is dihedral, while the 
other one is dense!

\begin{proof} We shall prove that the $\SL(2,\R)$-locus, i.e. the real
semi-algebraic set 
$X$ of $[-2,2]^4$ defined by 
$$
2(a^2+b^2+c^2+d^2)-abcd-16>\sqrt{(4-a^2)(4-b^2)(4-c^2)(4-d^2)},
$$
consist in connected neighborhoods
of those $8$ vertices corresponding to the Cayley surface
$$(a,b,c,d)=(\epsilon_1\cdot 2,\epsilon_2\cdot 2,\epsilon_3\cdot 2,\epsilon_4\cdot 2),\ \ \ \epsilon_i=\pm1,\ \epsilon_1\epsilon_2\epsilon_3\epsilon_4=1.$$
Benedetto-Goldman symmetries act transitively on those components.
On the other hand, the Cayley surface also arise for 
$(a,b,c,d)=(0,0,0,0)$, which is in the $\SU(2)$-locus:
the Okamoto correspondence therefore sends any of the $8$
components above into the $\SU(2)$-locus, 
thus proving the theorem.

By abuse of notation, still denote by $Z$ the discriminant locus 
defined by $\{\Delta=0\}\subset(-2,2)^4$.
The restriction $Z_{a,b}$ of $Z$ to the slice
$$
\Pi_{a,b}=\{(a,b,c,d)\ ;\ c,d\in(-2,2)\}, \quad (a,b)\in(-2,2)^2,
$$
is the union of two ellipses, namely those defined by
$$
c^2+d^2-\delta cd+\delta^2-4=0,\ \ \ \text{where}\ \delta=\frac{1}{2}\left(ab\pm\sqrt{(4-a^2)(4-b^2)}\right).
$$

\begin{figure}[h]\label{fig:DeltaSlice}
\input{DeltaSlice.pstex_t}
\caption{ {\sf{$Z$ restricted to the slice $\Pi_{a,b}$. }}}
\end{figure}
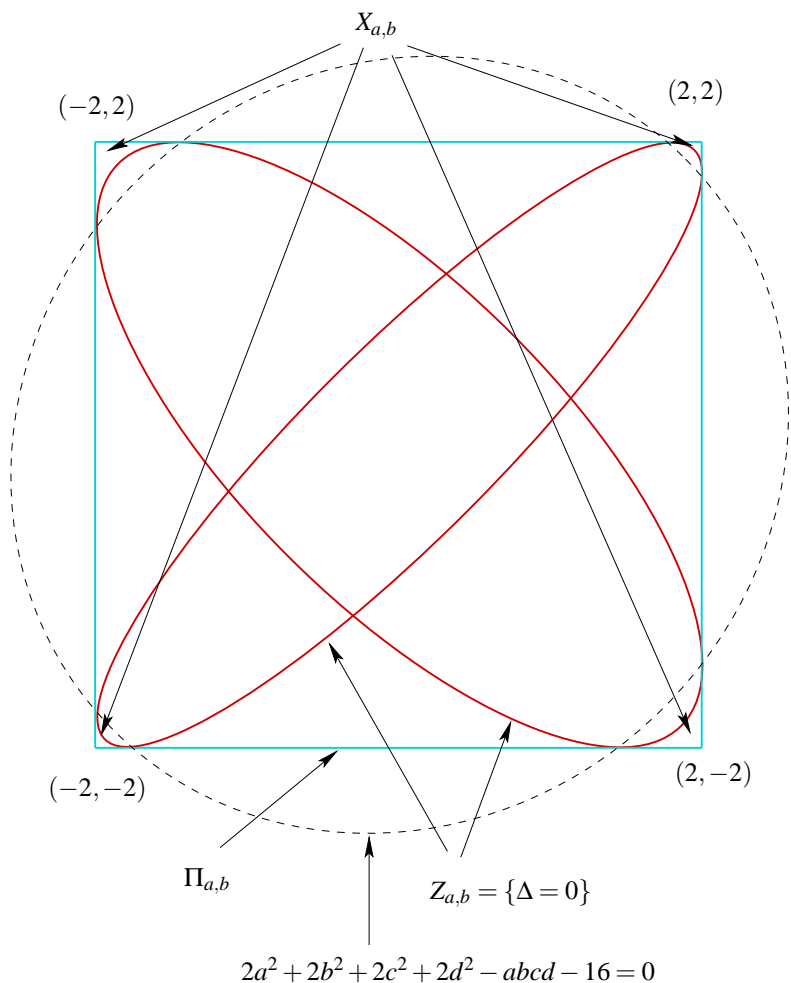

Those two ellipses are circumscribed into the square 
$\Pi_{a,b}$ (see figure \ref{fig:DeltaSlice}) and, for generic parameters $a$ and $b$,
cut the square into $13$ connected components.
One easily verify that $\SL(2,\R)$-components 
(namely those connected components of $X_{a,b}=X\cap\Pi_{a,b}$ 
defined by 
the inequality of the previous theorem) are those $4$ neighborhoods
of the vertices of the square.

This picture degenerates precisely when $a=\pm2$, $b=\pm2$ 
or $a=\pm b$. We do not need to consider the first two cases, 
since they are on the boudary of $(-2,2)^4$. Anyway, in these cases,
the two ellipses coincide; they moreover degenerate to a double line when $a=\pm b$.

In the last case $a=\pm b,$ the picture bifurcates.
When $a=b$, one of the ellipses degenerates to the double line
$c=d$, and the two components of $X_{a,b}$ near the vertices
$(2,2)$ and $(-2,-2)$ collapse. When $a=-b$, the
components of $X_{a,b}$ near the two other vertices collapse as well.
This means that each component of $X_{a,b}$ stands for exactly
two components of $X$: We finally obtain $8$ connected components
for the $\SL(2,\R)$-locus $X\subset(-2,2)^4$. One easily verify
that there are sixteen $\SU(2)$-components in $(-2,2)^4\setminus Z$.
\end{proof}

%%%%%%%%%%%%%%%%%%%%%%%%%%%%%%%%%

\subsection{Ramified covers}\label{section:RamifiedCovers}

%%%%%%%%%%%%%%%%%%%%%%%%%%%%%%%%%%

Here, we would like to describe other kinds of correspondences
between surfaces $S_{(A,B,C,D)},$ that arise by lifting representations
along a ramified cover of $\Sphere_4.$ Let $\rho\in\Rep(\Sphere_4)$ be a representation
 with $a=d=0$, so that $\rho(\alpha)^2=\rho(\beta)^2=-I$, and consider 
 the two-fold cover $\pi:\Sphere\to\Sphere$ ramifying over $p_\alpha$
 and $p_\delta$.

\begin{figure}[h]\label{fig:covering}
\input{cover.pstex_t}
\caption{ {\sf{The two-fold cover. }}{{(the point $p_\delta$ is at infinity)}}}
\end{figure}
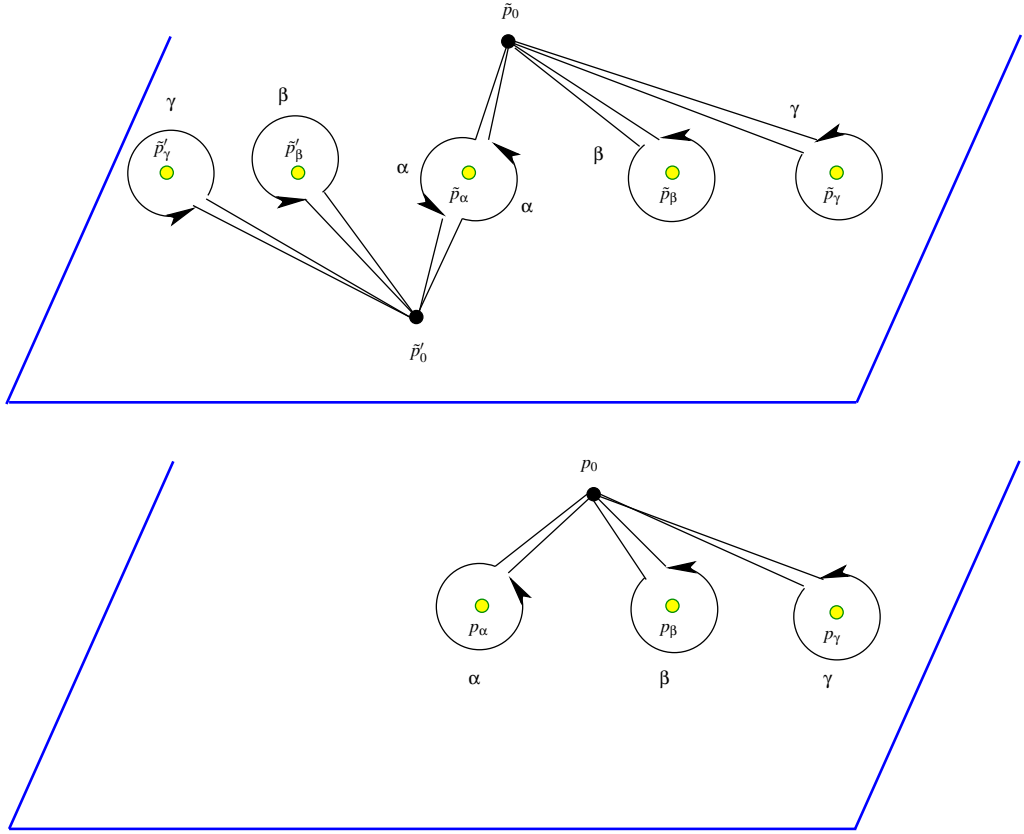

The four punctures lift-up as six punctures labelled in the obvious way
$$
\pi\ :\ \left\{\begin{matrix}
\tilde p_\alpha&\mapsto&p_\alpha\\
\tilde p_\beta,\tilde p_\beta'&\mapsto&p_\beta\\
\tilde p_\gamma,\tilde p_\gamma'&\mapsto&p_\gamma\\
\tilde p_\delta&\mapsto&p_\delta
\end{matrix}\right.
$$
After twisting the lifted representation $\rho\circ\pi$ by $-I$ at $\tilde p_\alpha$
and $\tilde p_\delta$, we get a new representation 
$\tilde\rho\in\Rep(\Sphere_4)$ ; the new punctures
are respectively $\tilde p_\gamma'$, $\tilde p_\gamma$, 
$\tilde p_\beta$ and $\tilde p_\beta'$ and the new generators
for the fundamental group are given by $\alpha\beta\gamma\beta^{-1}\alpha^{-1},$
 $\alpha\beta\alpha^{-1},$  $\beta,$ and $\gamma.$
After computation, we get a map 
$$
\left\{\begin{matrix}
0&\mapsto&c\\
b&\mapsto&b\\
c&\mapsto&b\\
0&\mapsto&c
\end{matrix}\right.
\ \ \ \text{and}\ \ \ 
\left\{\begin{matrix}
x&\mapsto&y\\
y&\mapsto&2-x^2\\
z&\mapsto&x^2y+xz-y+bc
\end{matrix}\right.
$$
defining a two-fold cover 
$$
\Quad\ :\ S_{(0,B,0,D)}\to S_{(2B,4-D,2B,2D-B^2-4)}
$$
where $B=bc$ and $D=4-b^2-c^2$. This map corresponds to
the so-called ``quadratic transformation'' of Painlev\'e VI equation. 

When moreover $c=0$,
we can iterate twice this transformation and we deduce a $4$-fold cover
$$
\Quad\circ\Quad\ :\ S_{(0,0,0,D)}\to S_{(8-2D,8-2D,8-2D,-28+12D-D^2)}
$$
$$
\left\{\begin{matrix}
0&\mapsto&b\\
b&\mapsto&b\\
0&\mapsto&b\\
0&\mapsto&b
\end{matrix}\right.
\ \ \ \text{and}\ \ \ 
\left\{\begin{matrix}
x&\mapsto&2-x^2\\
y&\mapsto&2-y^2\\
z&\mapsto&2-z^2
\end{matrix}\right.
$$
For instance, when $D=0$, we get a covering $S_{(0,0,0,0)}\to S_{(8,8,8,-28)}$. Another particular case arise when $D=4$ where
$\Quad$ defines an endomorphism of the Cayley cubic surface
$S_{(0,0,0,4)}\to S_{(0,0,0,4)}$, namely that one induced by the regular cover
$$
\C^*\times\C^*\to\C^*\times\C^*\ :\ (u,v)\mapsto(v,u^2).
$$
\begin{eg}
By the way, we note that, up to the action of $Q$, the following
traces data are related :
$$
\begin{matrix}
(0,0,0,d)&\leftrightarrow&(d'',d'',d'',-d'')&\rightarrow&S_{(0,0,0,4-d^2)}\\
\downarrow&&&&\downarrow\\
(0,0,d,d)&\leftrightarrow&(2,2,d',-d')&\rightarrow&S_{(d^2,0,0,4-2d^2)}\\
\downarrow&&&&\downarrow\\
(d,d,d,d)&\leftrightarrow&(2,2,2,d^2-2)&\rightarrow&S_{(2d^2,2d^2,2d^2,4-4d^2-d^4)}
\end{matrix}
$$
where $d'=\sqrt{4-d^2}$ and $d''=\sqrt{2+d'}$. In the previous diagram, horizontal correspondences arise from Okamoto symmetries, while
vertical arrows, from quadratic transformation $Q$.

More generally, we have related
$$
\begin{matrix}
(0,0,c,d)&\leftrightarrow&(c'',c'',d'',-d'')&\rightarrow&S_{(cd,0,0,4-c^2-d^2)}\\
\downarrow&&&&\downarrow\\
(c,c,d,d)&\leftrightarrow&(2,2,c',d')&\rightarrow&S_{(c^2+d^2,2cd,2cd,4-2c^2-2d^2-c^2d^2)}\\
\end{matrix}
$$
where $c'=\frac{cd+\sqrt{(c^2-4)(d^2-4)}}{2}$, $d'=\frac{cd-\sqrt{(c^2-4)(d^2-4)}}{2}$, $c''=\sqrt{2+c'}$ and $d''=\sqrt{2-c'}$.
\end{eg}

\begin{rem} One can check by direct computations that the quadratic tranformation
$\Quad$ is equivariant, up to finite index, with respect to the 
$\Gamma_2^*$-actions. Precisely, we have
$$
\left\{\begin{matrix}
\Quad\circ B_1^2&=&B_2^{-1}\circ\Quad,\\
\Quad\circ B_2&=&B_1^{-2}\circ\Quad,\\
\Quad\circ s_z&=&s_z\circ\Quad.
\end{matrix}\right.
$$
The group generated by 
$B_1^2$, $B_2$ and $s_z$,
acting on both sides,
contains $\Gamma_2^*$ as an index $2$ subgroup 
(recall that $B_1^2=g_x=s_z\circ s_y$ and $B_2^2=g_y=s_x\circ s_z$).
Therefore, if $q=\Quad(p)$ (for some parameters $(0,B,0,D)$),
then $p$ is $\Gamma_2^*$-periodic (resp. bounded) if, and only if,
$q$ is.
\end{rem}

%
%%%%%%%%%%%%%%%%%%%%%%%%%%%%%%%%%%%%%%%%%%%%%%%%%%%%%%%%%%%%%%%%%%
%

\bibliographystyle{plain}
\bibliography{references}

%@article {surf,
%    AUTHOR = {Endrass, Stephan and Huelf, Hans and Oertel, Ruediger and Schmitt, %Ralf and Schneider, Kai and Beigel, Johannes},
 %    TITLE = {surf version 1.0.4 - visualizing algebraic curves and algebraic surfaces},
  % JOURNAL = {http://surf.sourceforge.net/index.shtml},
    %  YEAR = {2003},
%}

\end{document}

%% file: sphere-ter.pstex_t
\begin{picture}(0,0)%
\includegraphics{sphere-ter.pstex}%
\end{picture}%
\setlength{\unitlength}{2171sp}%
\begingroup\makeatletter\ifx\SetFigFont\undefined%
\gdef\SetFigFont#1#2#3#4#5{%
  \reset@font\fontsize{#1}{#2pt}%
  \fontfamily{#3}\fontseries{#4}\fontshape{#5}%
  \selectfont}%
\fi\endgroup%
\begin{picture}(5880,7375)(3211,-7472)
\put(3301,-7336){\makebox(0,0)[lb]{\smash{{\SetFigFont{12}{14.4}{\familydefault}{\mddefault}{\updefault}$p_\delta$}}}}
\put(5326,-3511){\makebox(0,0)[lb]{\smash{{\SetFigFont{12}{14.4}{\familydefault}{\mddefault}{\updefault}$\gamma$}}}}
\put(3226,-961){\makebox(0,0)[lb]{\smash{{\SetFigFont{12}{14.4}{\familydefault}{\mddefault}{\updefault}$p_\gamma$}}}}
\put(9076,-436){\makebox(0,0)[lb]{\smash{{\SetFigFont{12}{14.4}{\familydefault}{\mddefault}{\updefault}$p_\beta$}}}}
\put(9076,-6811){\makebox(0,0)[lb]{\smash{{\SetFigFont{12}{14.4}{\familydefault}{\mddefault}{\updefault}$p_\alpha$}}}}
\put(7051,-4861){\makebox(0,0)[lb]{\smash{{\SetFigFont{12}{14.4}{\familydefault}{\mddefault}{\updefault}$\alpha$}}}}
\put(7276,-3136){\makebox(0,0)[lb]{\smash{{\SetFigFont{12}{14.4}{\familydefault}{\mddefault}{\updefault}$\beta$}}}}
\put(5026,-4786){\makebox(0,0)[lb]{\smash{{\SetFigFont{12}{14.4}{\familydefault}{\mddefault}{\updefault}$\delta$}}}}
\end{picture}%

%% file: ch-surf.pstex_t
\begin{picture}(0,0)%
\includegraphics{ch-surf.pstex}%
\end{picture}%
\setlength{\unitlength}{1697sp}%
\begingroup\makeatletter\ifx\SetFigFont\undefined%
\gdef\SetFigFont#1#2#3#4#5{%
  \reset@font\fontsize{#1}{#2pt}%
  \fontfamily{#3}\fontseries{#4}\fontshape{#5}%
  \selectfont}%
\fi\endgroup%
\begin{picture}(13266,13623)(1168,-13951)
\put(11176,-6736){\makebox(0,0)[lb]{\smash{{\SetFigFont{9}{10.8}{\familydefault}{\mddefault}{\updefault}II}}}}
\put(3901,-13936){\makebox(0,0)[lb]{\smash{{\SetFigFont{9}{10.8}{\familydefault}{\mddefault}{\updefault}III}}}}
\put(11251,-13936){\makebox(0,0)[lb]{\smash{{\SetFigFont{9}{10.8}{\familydefault}{\mddefault}{\updefault}IV}}}}
\put(4051,-6736){\makebox(0,0)[lb]{\smash{{\SetFigFont{9}{10.8}{\familydefault}{\mddefault}{\updefault}I}}}}
\end{picture}%

%% file: torus.pstex_t
\begin{picture}(0,0)%
\includegraphics{torus.pstex}%
\end{picture}%
\setlength{\unitlength}{2171sp}%
\begingroup\makeatletter\ifx\SetFigFont\undefined%
\gdef\SetFigFont#1#2#3#4#5{%
  \reset@font\fontsize{#1}{#2pt}%
  \fontfamily{#3}\fontseries{#4}\fontshape{#5}%
  \selectfont}%
\fi\endgroup%
\begin{picture}(7341,7266)(-107,-7594)
\put(4276,-1936){\makebox(0,0)[lb]{\smash{{\SetFigFont{12}{14.4}{\familydefault}{\mddefault}{\updefault}$\omega_2$}}}}
\put(4201,-4636){\makebox(0,0)[lb]{\smash{{\SetFigFont{12}{14.4}{\familydefault}{\mddefault}{\updefault}$\alpha$}}}}
\put(5401,-3811){\makebox(0,0)[lb]{\smash{{\SetFigFont{12}{14.4}{\familydefault}{\mddefault}{\updefault}$\alpha$}}}}
\put(6151,-4336){\makebox(0,0)[lb]{\smash{{\SetFigFont{12}{14.4}{\familydefault}{\mddefault}{\updefault}$1$}}}}
\put(6826,-3661){\makebox(0,0)[lb]{\smash{{\SetFigFont{12}{14.4}{\familydefault}{\mddefault}{\updefault}$\tilde p_0$}}}}
\put(5476,-2611){\makebox(0,0)[lb]{\smash{{\SetFigFont{12}{14.4}{\familydefault}{\mddefault}{\updefault}$\beta$}}}}
\put(5476,-3286){\makebox(0,0)[lb]{\smash{{\SetFigFont{12}{14.4}{\familydefault}{\mddefault}{\updefault}$\omega_1$}}}}
\put(4726,-4336){\makebox(0,0)[lb]{\smash{{\SetFigFont{12}{14.4}{\familydefault}{\mddefault}{\updefault}$1/2$}}}}
\put(5476,-4936){\makebox(0,0)[lb]{\smash{{\SetFigFont{12}{14.4}{\familydefault}{\mddefault}{\updefault}$\tilde p_0'$}}}}
\put(3226,-4336){\makebox(0,0)[lb]{\smash{{\SetFigFont{12}{14.4}{\familydefault}{\mddefault}{\updefault}$(0,0)$}}}}
\put(2851,-3661){\makebox(0,0)[lb]{\smash{{\SetFigFont{12}{14.4}{\familydefault}{\mddefault}{\updefault}$\delta$}}}}
\put(3001,-3061){\makebox(0,0)[lb]{\smash{{\SetFigFont{12}{14.4}{\familydefault}{\mddefault}{\updefault}$1/2$}}}}
\put(3676,-2161){\makebox(0,0)[lb]{\smash{{\SetFigFont{12}{14.4}{\familydefault}{\mddefault}{\updefault}$\gamma$}}}}
\put(3301,-1486){\makebox(0,0)[lb]{\smash{{\SetFigFont{12}{14.4}{\familydefault}{\mddefault}{\updefault}$1$}}}}
\end{picture}%

%% file: conj_c.pstex_t
\begin{picture}(0,0)%
\includegraphics{conj_c.pstex}%
\end{picture}%
\setlength{\unitlength}{1579sp}%
\begingroup\makeatletter\ifx\SetFigFont\undefined%
\gdef\SetFigFont#1#2#3#4#5{%
  \reset@font\fontsize{#1}{#2pt}%
  \fontfamily{#3}\fontseries{#4}\fontshape{#5}%
  \selectfont}%
\fi\endgroup%
\begin{picture}(16208,6814)(1050,-7787)
\put(3451,-6736){\makebox(0,0)[lb]{\smash{{\SetFigFont{10}{12.0}{\familydefault}{\mddefault}{\updefault}$j_x$}}}}
\put(5701,-6736){\makebox(0,0)[lb]{\smash{{\SetFigFont{10}{12.0}{\familydefault}{\mddefault}{\updefault}$j_y$}}}}
\put(3901,-1261){\makebox(0,0)[lb]{\smash{{\SetFigFont{10}{12.0}{\familydefault}{\mddefault}{\updefault}$j_z=\infty$}}}}
\put(11326,-1336){\makebox(0,0)[lb]{\smash{{\SetFigFont{10}{12.0}{\familydefault}{\mddefault}{\updefault}$v_x$}}}}
\put(13951,-6811){\makebox(0,0)[lb]{\smash{{\SetFigFont{10}{12.0}{\familydefault}{\mddefault}{\updefault}$v_z$}}}}
\put(16351,-1336){\makebox(0,0)[lb]{\smash{{\SetFigFont{10}{12.0}{\familydefault}{\mddefault}{\updefault}$v_y$}}}}
\put(10651,-4411){\makebox(0,0)[lb]{\smash{{\SetFigFont{10}{12.0}{\familydefault}{\mddefault}{\updefault}${\mathbf{c}}$}}}}
\put(13876,-1336){\makebox(0,0)[lb]{\smash{{\SetFigFont{10}{12.0}{\familydefault}{\mddefault}{\updefault}${\mathbf{c}}(j_z)$}}}}
\put(15676,-4411){\makebox(0,0)[lb]{\smash{{\SetFigFont{10}{12.0}{\familydefault}{\mddefault}{\updefault}${\mathbf{c}}(j_y)$}}}}
\put(12676,-4411){\makebox(0,0)[lb]{\smash{{\SetFigFont{10}{12.0}{\familydefault}{\mddefault}{\updefault}${\mathbf{c}}(j_x)$}}}}
\end{picture}%

%% file: DeltaSlice.pstex_t
\begin{picture}(0,0)%
\includegraphics{DeltaSlice.pstex}%
\end{picture}%
\setlength{\unitlength}{3158sp}%
\begingroup\makeatletter\ifx\SetFigFont\undefined%
\gdef\SetFigFont#1#2#3#4#5{%
  \reset@font\fontsize{#1}{#2pt}%
  \fontfamily{#3}\fontseries{#4}\fontshape{#5}%
  \selectfont}%
\fi\endgroup%
\begin{picture}(6122,7707)(1365,-6655)
\put(3181,-6586){\makebox(0,0)[lb]{\smash{{\SetFigFont{10}{12.0}{\familydefault}{\mddefault}{\updefault}$2a^2+2b^2+2c^2+2d^2-abcd-16=0$}}}}
\put(1666,-5146){\makebox(0,0)[lb]{\smash{{\SetFigFont{10}{12.0}{\familydefault}{\mddefault}{\updefault}$(-2,-2)$}}}}
\put(1741,209){\makebox(0,0)[lb]{\smash{{\SetFigFont{10}{12.0}{\familydefault}{\mddefault}{\updefault}$(-2,2)$}}}}
\put(6586,-5026){\makebox(0,0)[lb]{\smash{{\SetFigFont{10}{12.0}{\familydefault}{\mddefault}{\updefault}$(2,-2)$}}}}
\put(6526,329){\makebox(0,0)[lb]{\smash{{\SetFigFont{10}{12.0}{\familydefault}{\mddefault}{\updefault}$(2,2)$}}}}
\put(2731,-5851){\makebox(0,0)[lb]{\smash{{\SetFigFont{10}{12.0}{\familydefault}{\mddefault}{\updefault}$\Pi_{a,b}$}}}}
\put(4666,-5956){\makebox(0,0)[lb]{\smash{{\SetFigFont{10}{12.0}{\familydefault}{\mddefault}{\updefault}$Z_{a,b}=\{\Delta=0\}$}}}}
\put(4081,869){\makebox(0,0)[lb]{\smash{{\SetFigFont{10}{12.0}{\familydefault}{\mddefault}{\updefault}$X_{a,b}$}}}}
\end{picture}%

%% file: cover.pstex_t
\begin{picture}(0,0)%
\includegraphics{cover.pstex}%
\end{picture}%
\setlength{\unitlength}{2171sp}%
\begingroup\makeatletter\ifx\SetFigFont\undefined%
\gdef\SetFigFont#1#2#3#4#5{%
  \reset@font\fontsize{#1}{#2pt}%
  \fontfamily{#3}\fontseries{#4}\fontshape{#5}%
  \selectfont}%
\fi\endgroup%
\begin{picture}(11651,9516)(387,-8194)
\put(6301,-1111){\makebox(0,0)[lb]{\smash{{\SetFigFont{7}{8.4}{\familydefault}{\mddefault}{\updefault}$\alpha$}}}}
\put(2101,-436){\makebox(0,0)[lb]{\smash{{\SetFigFont{7}{8.4}{\familydefault}{\mddefault}{\updefault}$\tilde p_\gamma'$}}}}
\put(7126,-511){\makebox(0,0)[lb]{\smash{{\SetFigFont{7}{8.4}{\familydefault}{\mddefault}{\updefault}$\beta$}}}}
\put(9376, 14){\makebox(0,0)[lb]{\smash{{\SetFigFont{7}{8.4}{\familydefault}{\mddefault}{\updefault}$\gamma$}}}}
\put(3526,164){\makebox(0,0)[lb]{\smash{{\SetFigFont{7}{8.4}{\familydefault}{\mddefault}{\updefault}$\beta$}}}}
\put(2251, 89){\makebox(0,0)[lb]{\smash{{\SetFigFont{7}{8.4}{\familydefault}{\mddefault}{\updefault}$\gamma$}}}}
\put(9751,-961){\makebox(0,0)[lb]{\smash{{\SetFigFont{7}{8.4}{\familydefault}{\mddefault}{\updefault}$\tilde p_\gamma$}}}}
\put(4876,-661){\makebox(0,0)[lb]{\smash{{\SetFigFont{7}{8.4}{\familydefault}{\mddefault}{\updefault}$\alpha$}}}}
\put(5701,-6511){\makebox(0,0)[lb]{\smash{{\SetFigFont{7}{8.4}{\familydefault}{\mddefault}{\updefault}$\alpha$}}}}
\put(7876,-6511){\makebox(0,0)[lb]{\smash{{\SetFigFont{7}{8.4}{\familydefault}{\mddefault}{\updefault}$\beta$}}}}
\put(9751,-6511){\makebox(0,0)[lb]{\smash{{\SetFigFont{7}{8.4}{\familydefault}{\mddefault}{\updefault}$\gamma$}}}}
\put(5701,-5911){\makebox(0,0)[lb]{\smash{{\SetFigFont{7}{8.4}{\familydefault}{\mddefault}{\updefault}$p_\alpha$}}}}
\put(7876,-5911){\makebox(0,0)[lb]{\smash{{\SetFigFont{7}{8.4}{\familydefault}{\mddefault}{\updefault}$p_\beta$}}}}
\put(9751,-5986){\makebox(0,0)[lb]{\smash{{\SetFigFont{7}{8.4}{\familydefault}{\mddefault}{\updefault}$p_\gamma$}}}}
\put(6976,-4036){\makebox(0,0)[lb]{\smash{{\SetFigFont{7}{8.4}{\familydefault}{\mddefault}{\updefault}$p_0$}}}}
\put(6076,1139){\makebox(0,0)[lb]{\smash{{\SetFigFont{7}{8.4}{\familydefault}{\mddefault}{\updefault}$\tilde p_0$}}}}
\put(5476,-961){\makebox(0,0)[lb]{\smash{{\SetFigFont{7}{8.4}{\familydefault}{\mddefault}{\updefault}$\tilde p_\alpha$}}}}
\put(7876,-961){\makebox(0,0)[lb]{\smash{{\SetFigFont{7}{8.4}{\familydefault}{\mddefault}{\updefault}$\tilde p_\beta$}}}}
\put(5026,-2761){\makebox(0,0)[lb]{\smash{{\SetFigFont{7}{8.4}{\familydefault}{\mddefault}{\updefault}$\tilde p_0'$}}}}
\put(3601,-436){\makebox(0,0)[lb]{\smash{{\SetFigFont{7}{8.4}{\familydefault}{\mddefault}{\updefault}$\tilde p_\beta'$}}}}
\end{picture}%